\input amstex
\documentstyle{amsppt}
\magnification1200 
\tolerance=2000
\NoRunningHeads
\def\n#1{\Bbb #1}

\def\Pr{\hbox{Pr }}
\def\cl{\hbox{cl }}

\def\rank{\hbox{rank }}
\def\Sing{\hbox{Sing }}
\def\lim{\hbox{lim }}
\def\im{\hbox{im }}
\def\Pl{\hbox{Pl }}

\def\End{\hbox{End }} 

\def\Ker{\hbox{Ker }}
\def\Tor{\hbox{Tor }}
\def\Coker{\hbox{Coker }}

\def\Pic{\hbox{Pic }}

\def\pr{\hbox{pr}}
\def\dim{\hbox{dim }}
\def\opr{\hbox{def}}
\def\det{\hbox{det }}
\def\Div{\hbox{Div }}
\def\deg{\hbox{deg }}
\def\codim{\hbox{codim }}
\def\Codim{\hbox{Codim }}
\def\const{\hbox{const }}
\def\Alb{\hbox{Alb }}

\topmatter
\title 
Fano threefolds of genus 6
\endtitle
\author
D. Logachev
\endauthor
\address 
logachev\@ usb.ve
Yaroslavl State Pedagogical Institute. Yaroslavl, Russia. 
\endaddress
\date November, 1982 \enddate
\thanks I am very grateful to my advisor A.S.Tikhomirov for the 
statement of the problem and numerous consultations. I wanted also to thank 
A.N.Tyurin and V.A.Iskovskih for some important remarks. \endthanks
\abstract Ideas and methods of {\it Clemens C.H., Griffiths Ph. The 
intermediate Jacobian of a cubic threefold} are applied to a Fano threefold 
$X$ of genus 6 --- intersection of $G(2,5) \subset P^9$ with $P^7$ and a 
quadric. We prove: 
\medskip
1. The Fano surface $F(X)$ of $X$ is smooth and irreducible. Hodge numbers 
and some other invariants of $F(X)$ are calculated. 

2. Tangent bundle theorem for $X$, and its geometric 
interpretation. It is shown that $F(X)$ defines $X$ uniquely. 

3. The Abel - Jacobi map $\Phi: \Alb F(X) \to J^3(X)$ is an isogeny. 

4. As a necessary step of calculation of $h^{1,0}(F(X))$ we describe a 
special intersection of 3 quadrics in $P^6$ (having 1 double point) 
whose Hesse curve is a smooth plane curve of degree 6. 

5. $\im \Phi(F(X)) \subset J^3(X)$ is algebraically equivalent to 
$\frac{2\Theta^8}{8!}$ where $\Theta \subset J^3(X)$ is a Poincar\'e 
divisor (a sketch of the proof). 
\endabstract
\endtopmatter
\document
{\bf \centerline{Contents}}
\medskip
\settabs 15 \columns
\+0. Introduction &&&&&&&&&&&&&&2\cr
\medskip
\+1. Intermediate Jacobian of algebraic varieties and Abel -
Jacobi map &&&&&&&&&&&&&&5\cr
\medskip
\+2. Fano threefolds: a brief survey &&&&&&&&&&&&&&7\cr
\medskip
\+3. Fano surface as a determinantal variety &&&&&&&&&&&&&&12\cr
\medskip
\+4. Tangent bundle theorem &&&&&&&&&&&&&&23\cr
\medskip
\+5. Special intersection of 3 quadrics in $P^6$ &&&&&&&&&&&&&&31\cr
\medskip
\+6. Abel - Jacobi map of $F(X)$ is an isogeny &&&&&&&&&&&&&&38\cr
\medskip
\+7. Geometric interpretation of the tangent bundle theorem\cr 
\+and recovering of $X$ by its Fano surface &&&&&&&&&&&&&&43\cr
\medskip
\+\phantom{8. }References&&&&&&&&&&&&&&48\cr
\bigskip
\bigskip
{\bf 0. Introduction. }
\bigskip
The main object of the present paper is a Fano threefold $X$ of genus 6 of 
the principal series which is the intersection of the Grassmann variety 
$G(2,5) \subset P^9$ with $P^7$ and a quadric. The surface that parametrizes 
conics on $X$ is called its Fano 
surface and is denoted by $F(X)$. The middle Jacobian of $X$ is denoted by 
$J^3(X)$. We investigate the Abel - Jacobi map 
$\Phi: F(X) \to J^3(X)$. This investigation is the main component of 
solution of problems of rationality and of global Torelli theorem for $X$. 

Main results obtained in the present paper are the following: 
\medskip
1. $F(X)$ is smooth and irreducible. Hodge numbers and some other invariants
of $F(X)$ are calculated; it is shown that they coincide with the 
corresponding invariants of the two-sheeted covering of $P^3$ ramified 
at a quartic. 
\medskip
2. Tangent bundle theorem for $X$ is proved, and its geometric 
interpretation is given. It is shown that $F(X)$ defines $X$ uniquely 
(under some restriction imposed on $X$, see (7.1)). 
\medskip
3. The Abel - Jacobi map $\Phi: \Alb F(X) \to J^3(X)$ is an isogeny. 
\medskip
4. As a necessary step of calculation of $h^{1,0}(F(X))$ we describe a 
special intersection of 3 quadrics in $P^6$ whose Hesse curve is a smooth 
plane curve of degree 6 (instead of 7 in the general case) and main 
theorems concerning a general intersection of 3 quadrics are proved for 
this case. 
\medskip
5. $\im \Phi(F(X)) \subset J^3(X)$ is algebraically equivalent to 
$\frac{2\Theta^8}{8!}$ where $\Theta \subset J^3(X)$ is a Poincar\'e 
divisor (a sketch of the proof). 
\medskip
Now we give a more detailed description of the present paper. Sections
1 and 2 are of the survey nature. Section 1 contains a definition and the
main properties of intermediate Jacobian of a variety, and of an Abel - 
Jacobi map. Section 2 contains a definition of Fano threefolds, a brief 
description of methods and results of proofs of their non-rationality. Main 
steps of proofs of non-rationality of a cubic threefold and of a two-sheeted 
covering of $P^3$ ramified at a quartic (according Clemens - Griffiths) are 
given. 

In Sections 3 and 4 we find some simple properties of $F(X)$, calculate 
some of its invariants and prove the tangent bundle theorem. Namely, let $V$ 
be a vector space of dimension 5, $G=G(2,V)$ a Grassmann variety, it is 
included into $P(\lambda ^2(V))$ via Pl\"ucker embedding. Let $H_1$, $H_2$ 
be hyperplanes in  $P(\lambda ^2(V))$ and $\Omega$ a quadric in  
$P(\lambda ^2(V))$ such that $X= G \cap H_1 \cap H_2 \cap \Omega$ 
is a smooth threefold. $X$ is a Fano threefold of the principal series that 
will be studied in the present paper. We consider only those $X$ that 
satisfy some conditions (3.6), (3.7) (which are fulfilled for a generic 
$X$), although most likely all results of the paper are true for all $X$. 
Fano surface $F_c=F_c(X)$ is a set of conics on $X$. 

Main results of Section 3 are the following: 

{\bf Proposition 0.1.} $F_c$ is a smooth surface containing a distinguished
point $c_{\Omega}$ and a distinguished straight line $l_2$ \ \ 
($c_{\Omega} \notin l_2)$. Let 
$r_F: F \to F_c$ be the blowing up of $F_c$ at $c_{\Omega}$ and 
$l_1 = r_F^{-1}(c_{\Omega})$. There exists an involution $i_F$ on $F$ 
without stable points such that $i_F \vert _{l_1} : l_1 \to l_2$ is an 
isomorphism. We denote the two-sheeted covering that corresponds to $i_F$ 
by $p_F: F \to F_0$, and we denote $l = p_F(l_1)$. $l_1$ and $l_2$ are 
exceptional lines on $F$, i.e. there are blowings down maps $r: F \to F_m$ 
(resp. $r_0: F_0 \to F_{0m}$) that send lines $l_1 \cup l_2$ (resp. $l$) 
to points $c_{\Omega} \cup c_{\Omega}'$ (resp. $c_{0\Omega}$), as well as 
an involution map $i_{F_m}: F_m \to F_m$ and the corresponding two-sheeted 
covering $p_{F_m}: F_m \to F_{0m}$, so we have a commutative diagram: 
$$\matrix & & F & & \overset{i_F}\to{\to} & & F \\
& \overset{r_F}\to{\swarrow} & & \overset{p_F}\to{\searrow} & & 
\overset{p_F}\to{\swarrow} & \\ F_c & & \downarrow & & F_0 & & \downarrow \\
& \searrow & & & \downarrow & & \\
& & F_m & & \overset{i_{F_m}}\to{\to} & & F_m \\
& & & \overset{p_{F_m}}\to{\searrow} & & \overset{p_{F_m}}\to{\swarrow} & \\
& & & & F_{0m} & & \endmatrix $$
where vertical maps are respectively $r$, $r_0$, $r$. 

We denote $Z=P_{P(V^*)}(O\oplus O(1))$, there are projections
$\eta: Z \to P(V^*)$ and $r_Z: Z \to P^5$ where $r_Z$ is the blowing up
map of a point $t \in P^5$. 
\medskip
{\bf Proposition 0.2.} There are canonical inclusions 
$\bar \phi: F_0 \hookrightarrow Z$, 
$\bar \phi _m : F_{0m} \hookrightarrow P^5$ making the following diagram 
commutative (vertical maps are $r_0$, $r_Z$): 
$$ \matrix F_0 & \overset{\bar \phi}\to{\hookrightarrow} & Z \\
\downarrow & & \downarrow \\
F_{0m} & \overset{\bar \phi_m}\to{\hookrightarrow} & P^5  \endmatrix $$
such that $ \bar \phi_m(c_{0\Omega}) =t$. Further, 
$P^5$ can be identified with a space which is dual to the space of 
quadrics in $P^7$ which contain $X$. The point $t \in P^5$ corresponds to 
a hyperplane of quadrics which contain $G \cap H_1
 \cap H_2$. 
\medskip
{\bf Proposition 0.3.} There exists a locally free sheaf $M$ of rank 4 
on $P(V^*)$ and a bundle of quadrics $O_{\eta}(-1) \to S^2(\eta^*(M^*))$ 
on $Z$ such that $F_0 \hookrightarrow Z$ is its rank 2 determinantal variety. 
\medskip
{\bf Proposition 0.4.} Let $\sigma _m \in \Pic F_{0m}$ and 
$\sigma = r^*(\sigma_m) \in \Pic F_{0}$ are 
sheaves of order 2 that 
correspond to two-sheeted coverings $p_{F_m}$, $p_F$  respectively. There is 
an explicit formula on $F_0$ for $\sigma$ (see (3.8)). The canonical sheaf of
$F_{0m}$ (resp. $F_m$) is $\omega_{F_{0m}}=\bar \phi _m^*(O_{P^5}(3))
\otimes \sigma_m$ (resp. $\omega_{F_{m}}= p_{F_m}^*
( \bar \phi _m^*(O_{P^5}(3))) $ ). 
\medskip
{\bf Proposition 0.5.} $\deg \bar \phi _m (F_{0m}) = 40$, 
$c_1(\Omega_{F_m})^2=720$, $c_2(\Omega_{F_m})=384$. Taking into consideration 
that $h^{0,0}(F) = h^{0,0}(F_0) = 1$ and that $h^{1,0}(F) = 10$, 
$ h^{1,0}(F_0) = 0$ (see Theorem 0.13) we get the following table of Hodge 
numbers $h^{i,j}(F_m)$: 
$$ \matrix 1 & 10 & 101 \\ 10 & 220 & 10 \\ 101 & 10 & 1 \endmatrix 
\eqno{(0.6)} $$
and their decomposition on $i_{F_m}$-invariant part 
$h^{i,j}(F_m)^+ =h^{i,j}(F_{0m})$:
$$ \matrix 1 & 0 & 45 \\ 0 & 100 & 0 \\ 45 & 0 & 1 \endmatrix \eqno{(0.7)} $$
and on $i_{F_m}$-antiinvariant part $h^{i,j}(F_m)^- =h^{i}(F_{0m}, 
\lambda^j \Omega_{F_{0m}} \otimes \sigma _m)$:
$$ \matrix 0 & 10 & 56 \\ 10 & 120 & 10 \\ 56 & 10 & 0 \endmatrix 
\eqno{(0.8)} $$
\medskip
{\bf Proposition 0.9.} A conic $c$ on $X$ and a conic $i_F(c)$ which is 
involutory to $c$, meet at 2 (possibly coinciding) points $\gamma _1(c)$, 
$\gamma _2(c)$. We define a surface $W=\cup _{c \in F} (\gamma _1(c) \cup
\gamma _2(c))$. Then for a generic
point $x \in W$ there exists only one (up to $i_F$) conic $c$ such that 
$x \in \{ \gamma _1(c), \gamma _2(c)\}$.
\medskip
{\bf Proposition 0.10.} $O_X(W) = O_{P^7}(21)\vert _X$. There are 39 conics 
passing through a fixed generic point on $X$. 
\medskip
{\bf Remark 0.11.} 0.4 -- 0.8 remain true for a Fano surface $F(X')$ of $X'$
--- a two-sheeted covering of $P^3$ ramified at a quartic $W$, while 0.9 and 
0.10 are not true for $X'$. The end of Section 3 contains a description of a 
special Fano threefold of genus 6 defined in [1] --- a two-sheeted covering 
of a section of $G(2,5)$ by 3 hyperplanes ramified at the intersection with a 
quadric --- whose properties are similar to the ones of $X'$. 
\medskip
Section 4 of the present paper contains a proof of the tangent bundle theorem 
for $F(X)$. Let us give the necessary definitions. There exists a locally 
free sheaf $\tau_{2,M}$ of rank 2 on $F$ such that for $c \in F$ \ 
\ $P(\tau_{2,M}(c))$ is the linear envelope of $\gamma_1(c), \gamma_2(c)$. 
Let $i: l_1 \cup l_2 \to F$  be the natural inclusion. 
\medskip
{\bf Theorem 0.12.} (Tangent bundle theorem). There exists an exact sequence 
of sheaves on $F$:
$$0 \to \tau_{2,M}^* \to r^*(\Omega_{F_m}) \to i_*(2O_{l_1 \cup l_2}) \to 0$$
\medskip
The main resilt of Section 6 is a the following 
\medskip
{\bf Theorem 0.13.} $h^{1,0}(F)=10$, $h^{1,0}(F_0)=0$,  $h^{0,0}(F)=1$, 
$h^{0,0}(F_0)=1$. 
\medskip
The method of the proof is to consider a degeneration of $X$ in a family of 
Fano threefolds $f_X: X_4 \to \Delta$ whose base $\Delta$ is a complex 
analytic neighbourhood of 0. A generic fibre of $f_X$ is a smooth Fano 
threefold, the special fibre $X_0=f_X^{-1}(0)$ has one double point. $X_0$ 
is birationally equivalent to a special intersection of 3 quadrics in $P^6$ 
denoted by $X_0'$.  

Let us consider main objects related to a generic intersection of 3 quadrics 
in $P^{n-1}$ (denoted by $X_8^{n-4}$) for an odd $n$ ([28], [47] for any $n$, 
[21], [50] for $n=7$). Namely,  attached to $X_8^{n-4}$ is a Hesse curve 
$H_n$ which is a plane curve of degree $n$, and its two-sheeted non-ramified 
covering $p_{H_n}: \hat H_n \to H_n$. We define $S=S(p_{H_n})$ as a connected 
component of the set of effective divisors of degree $n$ on $\hat H_n$ such 
that $d \in S \iff O_{H_n}({p_{H_n}}_ *(d)) = O_{P^2}(1) \vert _{H_n}$. We 
define the Fano surface $F(X_8^{n-4})$ as the set of $\frac{n-5}2$-dimensional
quadrics on $X_8^{n-4}$. Then:
\medskip
(1) There exists an isomorphism $\Pr \hat H_n / H_n \to J^{n-4}(X_8^{n-4})$ 
([21]), here Pr is the Prym variety; 
\medskip
(2)  There exists an isogeny $\Alb S \to \Pr \hat H_n / H_n$ ([50]);
\medskip
(3) There exists an isomorphism $F(X_8^{n-4}) \to S$ ([47]). 
\medskip
Section 5 contains a proof of the following 
\medskip
{\bf Theorem 0.14.} An analog of a Hesse curve for $X_0'$ is a smooth plane 
curve $H_6$ of degree 6. $\Sing (F(X_0))$ is a double curve which is 
isomorphic to $H_6$. Analogs of the above (1) - (3) are true for $X_0'$ after 
replacement of $F(X_0)$ and $X_0'$ by their desingularizations 
$\widetilde{F(X_0)}$ and $\tilde X_0'$ respectively. Particularly, 
$h^{1,0}(\widetilde{F(X_0)}) = \dim \Pr (\hat H^6 /H^6) = 9$. 
\medskip
To finish a proof of Theorem 0.13, we use a Clemens - Schmidt exact sequence 
([40]) for a family $f_F: F_3 \to \Delta$ whose fibres are Fano surfaces of 
fibres of $f_X$. 

The next result of Section 6 is
\medskip
{\bf Theorem 0.15.} The Abel - Jacobi map $\Phi: \Alb F(X) \to J^3(X)$ is an 
isogeny. 
\medskip
Finally, a sketch of a proof of the following theorem is given: 
\medskip
{\bf Theorem 0.16.} $\Phi (F) \subset J^3(X)$ is algebraically equivalent to 
$\frac{2\Theta^8}{8!}$ where $\Theta \subset J^3(X)$ is a Poincar\'e divisor.
\medskip
The idea of the proof is the following. A result of Beauville ([28]) implies 
that an analog of this result is true for the image of $S(p_{H_6})$ in $\Pr 
\hat H^6 /H^6$ (i.e. that it is equivalent to $\frac{2\Theta^7}{7!}$). 
Investigation of topology of degeneration of abelian varieties given in [25], 
[31], permits to deduce (0.16) from this property. 
\medskip
Section 7 contains proofs of 2 theorems: 
\medskip
{\bf Theorem 0.17.} (Geometric interpretation of the tangent bundle theorem). 
Let us denote $V_{10} = H^0(F, \Omega_F)^*$, and let $B_{10}: F_m \to G(2, 
V_{10})$ be a map defined by $\Omega_{F_m}$. Let $B_8: F \to G(2, V_8)$ be 
a map defined as follows: $
B_8(f)= <\gamma_1(f), \gamma_2(f)>$ (here $V_8 = H_1 \cap H_2 \subset
\lambda^2(V)$). Then there exists a natural projection $p_{10,8}: V_{10} 
\to V_8$ such that 
$$B_8 = G(2, p_{10,8}) \circ B_{10} \circ r$$
(recall that $r:F \to F_m$ is a blowing down map). If $B_{10}$ is regular 
at $c_{0\Omega} \in F_m$ then $p_{10,8}$ is a 
projection from a plane $B_{10}(c_{0\Omega}) \subset V_{10}$. 
\medskip
{\bf Theorem 0.18.} Let $X$ satisfies the property: $B_{10}$ is regular at 
$c_{\Omega} \in F_m$.  Then, if for another Fano threefold $X'$ of the same 
type we have $F_c(X)=F_c(X')$ then $X=X'$.
\medskip
The proof of this theorem is based on the above geometric interpretation of 
the tangent bundle theorem that permits to recover uniquely a map $B_8: F \to 
G(2, V_8)$ by a given $F_c(X)$, and then to recover uniquely $X$ by  $B_8 $. 
\bigskip
{\bf Section 1. Intermediate Jacobian of algebraic varieties and Abel -
Jacobi map.}
\bigskip
Recall that this section is of survey nature. Let $Y$ be a smooth 
$n$-dimensional algebraic variety over $\n C$. We consider 
the Hodge decomposition of its Betti cohomology of an odd dimension $k$: 
$$H^k(Y, \n C)=\sum_{i+j=k}H^{i,j}(Y, \n C)$$
and we denote
$$H^+=H^k(Y, \n C)^+ = \sum_{i+j=k, i>j}H^{i,j}(Y, \n C)$$
$$H^-=H^k(Y, \n C)^- = \sum_{i+j=k, i<j}H^{i,j}(Y, \n C)$$
so we have 
$$H^+ = \overline{H^-}; \; \; H^+ \cap H^-=0\eqno{(1.1)}$$
(1.1) implies that the image of $H^k(Y, \n Z) $ in $H^-$ under the projection 
along $H^+$ will be a lattice. A complex torus $H^k(Y, \n C)^-/ \im 
H^k(Y, \n Z)$ is called the $k$-th intermediate Jacobian of $Y$, it is 
denoted by $J^k(Y)$. For $k=n$ it is called also the middle Jacobian. 

To define a structure of a polarized abelian variety on a complex torus $V/R$, 
it is necessary to define an hermitian form $H$ on $V$ satisfying conditions: 

(i) $H$ is positively defined; 

(ii) $\im H$ restricted on $R$ takes integer values. 
\medskip
There exists a Hodge real bilinear form $Q$ on $H^k(Y, \n C)$ (see, for
example, [39]) which is defined (for $k \le n$) as follows:
$$Q(\xi, \eta) = \xi \wedge \eta \wedge L^{n-k} \in H^{2n}(Y, \n C)$$
Here $H^{2n}(Y, \n C)$ is identified with $\n C$, $L \in H^{1,1}(Y, \n C) 
\cap H^2(Y, \n Z)$ is the fundamental class of a hyperplane section of $Y$ 
or, which is the same, the cohomology class that corresponds to a K\"ahler 
form on $Y$. Now we define $H$ on $H^-$ as follows: 
$$H(\xi, \eta) = i Q(\xi, \bar \eta)$$
$H$ is hermitian, because $Q$ is skew symmetric on odd-dimensional cohomology. 
A map 
$$\wedge L^{n-k}: H^k(Y, \n Z) \to H^{2n-k}(Y, \n Z)$$
is an inclusion, and a pairing $ H^k(Y, \n Z) \otimes H^{2n-k}(Y, \n Z) \to 
\n Z$ is unimodular, according a Poincar\'e theorem. This means that (ii) is 
satisfied for $H$, and moreover if $\wedge L^{n-k}$ is an isomorphism then 
$B(H) \overset{\opr}\to{=} 
\im H$ is unimodular. Nevertheless, (i) is not satisfied in general. 
We denote $k=2l+1$. The Hodge - 
Riemann bilinear relations show that a sufficient condition for $J^k(Y)$ to 
be an abelian variety is the following:
$H^-=H^{l,l+1}(Y, \n C)$ and $H^{l-1,l}(Y, \n C)=0$ (i.e. all elements of 
$ H^{l,l+1}(Y, \n C)$ are primitive). These conditions are always true for 
$k=1$, and $J^1(Y) = \Pic^0(Y)$. 

For $k=n$ \ \ \ $\wedge L^{n-k}$ is an identical operator, hence $B(H)$ is
unimodular. Complex torus $V/R$ having an hermitian form $H$ on $V$ satisfying 
(ii) such that $B(H)$ is unimodular, is called a principal torus. So, a middle 
Jacobian $J^n(X)$ is a principal torus (Griffiths, [38]). 

For a principal torus $T$ (like for principally polarized abelian varieties) 
there exists a notion of the dual torus which is isomorphic to $T$. For 
$J^n(Y)= H^-/\im H^n(Y, \n Z)$ the dual torus is $(H^+)^*/ \im H_n(Y, \n Z)$. 
This is a convenient way to define the Abel - Jacobi map (it is possible to 
define it for all $J^k(Y)$  without consideration of the dual torus, see for 
example [42], but we do not need it). 

Let $F$ be a family of $l$-dimensional algebraic cycles on $Y$. The 
Abel - Jacobi 
map $ \Phi: F \to J^n(Y)$ (or $ \Phi: \Alb F \to J^n(Y)$) is a generalization 
for $n \ge 1$ of a map $Y \to \Alb Y$ for $n=1$. A family of cycles $F$ is 
defined by its graph $\Gamma_F \overset{i}\to{\hookrightarrow} F \times Y$ 
such that for $f \in F$ $(p_F \circ i)^{-1}(f)$ is a cycle that corresponds 
to $f$, here $p_F$ is a projection from $F \times Y$ to $F$. Then there are 
cylinder maps 
$$\Phi_{\n Z}: H_1(F, \n Z) \to H_n(Y, \n Z)$$
$$\matrix \Phi_{\n C}: & H^1(F, \n C)^* & \to & H^n(Y, \n C)^* \\
& \parallel & & \parallel \\
& H_1(F, \n C) & \to & H_n(Y, \n C) \endmatrix $$
Since $F$ is an algebraic family, $\Phi_{\n C}$ behaves good with respect to 
the Hodge decomposition: 
$$ \Phi_{\n C}(H^{i,j}(F, \n C)^*) \subset H^{i+l,j+l}(Y, \n C)^*$$
For $(i,j)=(1,0)$ we have $ \Phi_{\n C}(H^{1,0}(F, \n C)^*) \subset 
H^{l+1,l}(Y, \n C)^*$ and hence $\Phi_{\n C}$ and $\Phi_{\n Z}$ define a map 
$\Phi: \Alb F \to J^n(Y)$ which is called the Abel - Jacobi map. 

Further we shall consider only the middle Jacobian of threefolds. To have a 
possibility to distinguish rational and non-rational threefolds we should 
know what is the behaviour of the Jacobian under the simplest rational maps. 
\medskip
{\bf Proposition 1.2.} ([32]; see also [20]). Let $Y$ be a threefold and 
$\tilde Y$ be the 
blowing up of a point (resp. of a curve $C$) on $Y$. Then $J^3(\tilde Y) = 
J^3(Y)$ (resp. $J^3(\tilde Y) = J^3(Y) \oplus J(C)$). $\square$
\medskip
According Hironaka ([41]) any birational map of algebraic varieties over 
$\n C$ is a composition of blowings up and down. This implies a general 
theorem which is used to apply the method of middle Jacobian to proofs of 
non-rationality of threefolds: 
 \medskip
{\bf Theorem 1.3.} ([32]) Middle Jacobian of a rational threefold is an abelian
variety which is isomorphic to the direct sum of Jacobians of curves. $\square$
\medskip
To apply this theorem we need to know properties of Jacobians of curves. 
Here is one of them: 
\medskip
{\bf Theorem 1.4.} (Andreotti - Mayer criterion, [24]). Let $(X, \Theta)$ be 
a principally polarized abelian variety. If $(X, \Theta)$ is the Jacobian of 
a curve $C$ then $\Codim _X \Sing \Theta =4$ (resp. 3) if $C$ is not a 
hyperelliptic curve (resp. if $C$ is a hyperelliptic curve). $\square$
\bigskip
{\bf Section 2. Fano threefolds: a brief survey.}
\bigskip
A threefold $X$ is called a Fano threefold if its anticanonical sheaf 
$K_X^{-1}$ is ample on $X$. A detailed survey on Fano threefolds can be 
found in [2], 
[3], [16]. According a Kodaira criterion, a Fano threefold $X$ has an 
important property $H^i(X, O_X)=0$ for
 $i >0$, i.e. $X$ satisfies an obvious necessary condition of rationality. 
$J^3(X)$ is an abelian variety whose dimension is $h^{2,1}(X)$. 

$X$ is called a Fano threefold of genus 1 if $\Pic (X) = \n Z$. We shall
consider below only these varieties. Main invariants of $X$ are the 
following: an index $r$ is defined by the equality $K_X^{-1} = O_X(r)$ where 
$O_X(1)$ is the positive generator of
 $\Pic (X)$; a degree $d$ is the degree of the image of $X$ under the map 
defined by the linear series $\vert K_X^{-1} \vert$, and a genus $g$ is the 
genus of the generic curve which is a section of the image of $X$ under this 
map by a subspace of codimension 2. 

Fano threefolds of genus one are completely classified, a table is given in 
[2]. All these threefolds --- except $P^3$ and a quadric threefold --- have 
index 2 or 1. It is clear that a hyperplane section of the image of $X$ under 
the map of the linear series
 $\vert O_X(1) \vert$ is a Del Pezzo surface (resp. K3 surface) for $X$ of 
index 2 (resp. 1). It is known that the degree of a Del Pezzo surface can be 
any number between 1 and 9, and the degree of a K3 surface can be any even 
number. But the degree of Fano threefolds of index 2 (resp. 1) can be any 
number between 1 and 5 (resp. any even number from 2 until 22, except 20). It 
is clear that for Fano threefolds of index 1 we have $d=2g-2$, and the map 
defined by the linear series $\vert K_X^{-1} \vert$ sends a
 threefold of genus $g$ to $P^{g+1}$. 

Fano threefolds are counterexamples to the Luroth problem. 
Unirationality is proved for most types of Fano threefolds (see, for 
example, [16]). There exists a method of a proof of non-rationality of
a Fano threefold $X$ based on investigation of the group of birational 
automorphisms of $X$. The origin of this method is a paper of M. Noether 
about the Cremona group of 
$P^2$. Fano used it in [35] giving a non-strict proof of non-rationality 
of a quartic 
threefold $X^3_4$. A strict proof was given by Iskovskih - Manin in [6]. 
Later the same method was used to prove non-rationality of 
\medskip
(a) $X^3_2$ (a two-sheeted covering of $P^3$ of index 1) ([4]); 
\medskip
(b) The complete intersection of a quadric and a cubic in $P^5$ ([5]); 
\medskip
(c) A "double Veronese cone" ([22]). 
\medskip
The lesser is the degree of a threefold, the easier is the proof of 
non-rationality by this method. For example, $X^3_4$ and $X^3_2$ have the 
simplest group of birational automorphisms --- all their automorphisms are 
regular. 
\medskip
At the same time Clemens and Griffiths in [32] proved non-rationality of and 
the global Torelli theorem for a 
cubic threefold by a method of intermediate Jacobian. This result stimulated 
analogous investigations of other Fano threefolds. Non-rationality was 
proved for an intersection of 3 quadrics in $P^6$ ([21]) and the two-sheeted 
covering of $P^3$ ramified at a quartic ([12]). 

According [32], to prove non-rationality of a threefold $X$ it is sufficient 
to show that its intermediate Jacobian $J^3(X)$ is not the Jacobian of a 
curve. This can be done by investigation of Prym varieties ([27], [43]) or 
singularities of Poincar\'e divisor ([24]). 

Practically the only way of investigation of $J^3(X)$ is to study its 
Abel - Jacobi map $\Phi: F \to J^3(X)$ defined for any family of curves 
$F$ on $X$. It is a generalization of Albanese map for curves. Clemens 
and Griffiths in [32] considered the simplest family of curves 
on a cubic threefold $X$ --- its Fano surface $F(X)$ which is the family 
of all straight lines on $X$. The main technical result that is used for a
proof of non-rationality and global Torelli theorem for a cubic
threefold --- the description of Poincar\'e divisor $\Theta(J^3(X))$ 
on $J^3(X)$ --- is obtained in terms of its Abel - Jacobi map. 

A cubic threefold $X$ and related varieties $J^3(X)$, $\Theta(J^3(X))$, 
$F(X)$ were studied also in [17], [18], [23], [29], [44] --- [46], in 
earlier Fano papers [34], [36] (this list is not exhaustive). Global 
Torelli theorem and the description of singularities of Poincar\'e divisor 
for a cubic threefold was proved independently in [18], [19]. 

Let us consider now the main steps of [32] in more details. Firstly 
we give some general definitions. 

Let $\phi: Y \to A$ be a map of a $k$-dimensional variety $Y$ to an abelian 
variety $A$. For any $y \in Y$ a map of tangent spaces is defined: 
$T_{\phi}(y): T_Y(y) \to T_A(\phi(y))$. Let $\phi$ satisfy a property that 
$T_{\phi}$ is an inclusion at a generic point. Then there exists a Gauss map 
which is a rational map $G(\phi): Y \to G(k,T_A(0))$ where $G(k,T_A(0))$ is 
the Grassmann variety. The Gauss map is defined as follows. For 
$y \in Y$ \ \ \ $T_{\phi}(y)( T_Y(y))$ is an element of $G(k, T_A(\phi(y)))$.
The group law on an abelian variety defines a canonical isomorphism of
tangent spaces at any its points, and hence an isomorphism 
$G(k, T_A(\phi(y))) \to G(k, T_A(0))$. An image of the above element under 
this isomorphism is the desired element $G(\phi)(y)$. 

Let $\Phi: Y \to \Alb Y$ be the Albanese map. If its Gauss map exists then it 
coincide with the map $Y \to G(k, H^0(Y, \Omega_Y)^*)$ associated with the 
sheaf $\Omega_Y$ (since $T_{\Alb Y}(0)$ is canonically isomorphic to 
$ H^0(Y, \Omega_Y)$). 

Let now $V$ be a vector space of dimension 5 and $X \subset P(V)$ a cubic 
threefold. By definition, its Fano surface $F=F(X)$ is the set of straight 
lines on $X$. 
It is a subvariety of $G(2, V)$, we denote the corresponding inclusion 
$F \to G(2,V)$ by $i_F$. 
\medskip
{\bf Theorem 2.1.} $\dim \Alb F = h^{1,0}(F) = 5$. 
\medskip
There exist different ways of proof of this theorem. Altman and Kleiman 
([23]) use that $F$ is a set of zeros of a section of a locally free sheaf 
$E$ of rank 4 on $G(2, V)$. This permits to get a resolvent 
$$0 \to \lambda^4(E^*) \to \lambda^3(E^*) \to \lambda^2(E^*) \to 
\lambda^1(E^*) \to O_G \to O_F \to 0 $$ of $O_F$ on $G(2, V)$ and to use it
for an explicit calculation of $H^1(F, O_F)$. In principle, this method can
be applied to all families of subvarieties on a variety, but for more 
complicated Fano threefolds calculations are too long. 

A proof of this theorem given in [32] uses a technique of degeneration of 
$X$ to a threefold $X_0$ having one double point $x_0$. A set of straight 
lines on 
$X_0$ which contain $x_0$ is a smooth curve $D_0$ on $F_0=F(X_0)$; it is easy
 to calculate that its genus is 4 and that $\Sing F_0 = D_0$. 
The tangent cone to $X_0$ at $x_0$ has 2 rulings. This implies: if 
$r: \tilde F_0 \to F_0$ is a desingularization of $F_0$ then $r^{-1}(D_0)=
D_1 \cup D_2$ is a pair of curves which are isomorphic to $D_0$. Further, we 
have: $\tilde F_0 = S^2(D_0)$, $h^{1,0}(\tilde F_0)=h^{1,0}(S^2(D_0)) = 
g(D_0) =4$. Finally, degeneration theory for $X$ and $F$ implies that 
$h^{1,0}(F) = h^{1,0}(\tilde F_0) + 1 =5$. The situation for the Fano 
threefolds considered in the present paper is the same (see Proposition 5.12). 
\medskip
Theorem 2.1 implies the following 
\medskip
{\bf Corollary 2.2.} Abel - Jacobi map $\Phi: \Alb F \to J^3(X)$ is an 
isogeny. $\square$
\medskip
Really, $\Phi$ is an isomorphism, but we do not need it.

The next step is a proof of the tangent bundle theorem for $F$. It 
establishes a connection between $\Omega_F$ (or $T_F$) and a sheaf defined 
on $X$ "geometrically". For the case of a cubic threefold its statement is 
very simple: 
\medskip
{\bf Theorem 2.3.} There is an isomorphism of sheaves on $F$: 
$\Omega_F=i_F^*(\tau_{2,V}^*)$. 
\medskip
A proof of this theorem for the Fano surface of a cubic threefold is also 
very easy (comparatively to proofs of analogous theorems for other Fano 
threefolds). The main idea of the proof is a calculation of $T_F$ by means 
of deformation theory. 

Theorems 2.1 and 2.3 imply
\medskip
{\bf Corollary 2.4.}(Geometric interpretation of the tangent bundle theorem). 
 There exists an isomorphism $H^0(F, \Omega_F)^* \to V$ such that the 
canonical map $F \to G(2, H^0(F, \Omega_F)^* )$ coimcide with the inclusion 
$i_F: F \to G(2,V)$. $\square$
\medskip
Corollaries 2.2 and 2.4 imply that the Gauss map of $\Phi$ (i.e. a map
$G(\Phi): F \to T_{J^3(X)}(0)$ ) coincides with $i_F$.
\medskip
{\bf Theorem 2.5.} (Analog of a Riemann theorem). A set of points of type 
$\Phi(f_1) - \Phi(f_2)$ where $f_1, f_2 \in F$ is (up to a shift) a 
Poincar\'e divisor $\Theta \overset{i_{\Theta}}\to{\hookrightarrow} J^3(X)$. 
$\square$
\medskip
Theorem 2.5 and Corollary 2.4 are used for a description of the Gauss map of 
$i_{\Theta}$, i.e. a map $G(i_{\Theta}): \Theta \to G(4,V)=P(V^*)$. It is 
clear that if $f_1, f_2 \in F$ and the corresponding straight lines in $X$ 
do not meet one another, then (2.4) implies that  
$G(i_{\Theta}) (\Phi(f_1) - \Phi(f_2))$ 
is a space $P^3 \in P(V^*)$ spanned on $f_1$, $f_2$ considered as straight 
lines in 
$P(V)$. It is possible to prove that if $t$ is an intersection point of 
$f_1$, $f_2$ considered as straight lines in $P(V)$ then 
$G(i_{\Theta}) (\Phi(f_1) - 
\Phi(f_2))$ is the tangent space of $X$ at $t$. 

This means that $G(i_{\Theta}): \Theta \to G(4,V)=P(V^*)$ is a $n$-fold 
covering where $n=\# (G(i_{\Theta})^*(P^3))$ for a $P^3 \in P(V^*)$ is equal 
to the number of pairs of non-meeting straight lines on a smooth cubic surface 
$P^3 \cap X$. We have $n= 216$. This is sufficient to make a conclusion that 
$J^3(X)$ is not a Jacobian of a curve $C$ of genus 5. Really, for $J(C)$ 
the degree of covering of the Gauss map $G(i_{\Theta}): \Theta \to G(4,V)$ 
is clearly $\binom 84 = 70$. Another way to see that $J^3(X) \ne J(C)$ is to
investigate the ramification divisor of $G(i_{\Theta})$. For $J^3(X)$ (resp.
$J(C)$) it is the set of $P^3$ which are tangent to $X$ (resp. to $C$). 
Finally, the above recovering of $X$ by ramification divisor of 
$G(i_{\Theta})$ 
means that the global Torelli theorem for it is proved: 
\medskip
{\bf Theorem 2.6.} The polarized abelian variety $(J^3(X), \Theta)$ defines 
$X$ uniquely. $\square$ 
\medskip
Analogs of results of [32] for cubic threefolds are obtained now (1982) 
for the following Fano threefolds: 

1. Two-sheeted covering of $P^3$ with ramification in a quartic ([11], [12], 
[14], [49], [51]). Since many properties of $X^3_{10}$ --- a Fano threefold 
of genus 6 (the main object of the present paper) are analogous to the ones 
of a two-sheeted covering of $P^3$ with ramification in a quartic, a brief 
survey of them will be given below. 

2. $X_8^3$ --- an intersection of 3 quadrics in $P^6$ ([21], [28], [50], 
[51]). Main theorems concerning $X^3_8$ are given in Section 5, because they 
are necessary for a calculation of irregularity of $X^3_{10}$. Here there is 
no complete analogy: proofs of non-rationality and of global Torelli theorem 
for $X^3_8$ use a fact that $J^3(X^3_8)$ is a Prym variety. Tangent bundle
theorem for $X^3_8$ is unknown.

3. "Double cone Veronese" $X^3_1$ ([13]). Here the situation is more 
complicated than for other investigated types of Fano varieties. Namely, 
$F(X^3_1)$ is not reduced, and $X^3_1$ is the only known Fano threefold for 
which $ h^{1,0} (F(X^3_1)) \ne h^{2,1} (X^3_1) $ (here we have $ h^{1,0} 
(F(X^3_1)) = 2 h^{2,1} (X^3_1) $  ). 

An attempt to investigate the Fano surface of a quartic threefold $X^3_4$ was 
made in [33]. This subject is very difficult: for example, computer 
calculation was used to find a number of conics passing through a fixed 
generic point of $X^3_4$ (for Fano threefolds considered above this number 
can be found easily by hand calculations). 
\medskip
Let us consider now main steps of a proof of non-rationality of a two-sheeted 
covering of $P^3 \subset P(V_4)$ with ramification in a quartic $W \subset 
P^3$. Let $v: P^3 \to P^9=P(S^2(V_4))$ be a Veronese map (i.e. defined by the 
complete linear series $O_{P^3}(2)$), $t$ a point in a $P^{10} - P^9$, \ $K$ 
a cone over $v(P^3)$ and a vertex $t$, $p: K \to v(P^3)$ a projection from 
$t$, $\Omega$ a quadric in $P^{10}$ such that $t \notin \Omega$. Then "a 
double $P^3$" $X= K \cap \Omega$. Without loss of generality we can suppose 
that $P^9$ is the polar hyperplane for $\Omega$ with respect to $t$. Then
the ramification quartic $W = v^{-1}(\Omega \cap v(P^3))$.

Smoothness of $X$ implies smoothness of $W$. We consider below only the case 
when $W$ does not contain straight lines (this restriction is rejected in 
[15]). It 
is easy to see that a condition that $W$ does not contain straight lines is 
analogous 
to a condition of the lemma 3.7 of the present paper for $X^3_{10}$. 

If $c$ is a conic on $X$ then $p(c)$ is a conic on $v(P^3)$ and $v^{-1}(p(c))$ 
is a straight line in $P(V_4)$ which is tangential to $W$ at 2 points 
(a bitangent of 
$W$). We define a Fano surface $F=F(X)$ of $X$ as a set of conics on $X$, and 
we denote a set of bitangents of $W$ in $P(V_4)$ by $F_0$. A map $c \mapsto 
v^{-1}(p(c))$ is a two-sheeted non-ramified covering $p_F: F \to F_0$, we 
denote the corresponding involution on $F$ by $i_F$. All conics on $X$ are 
smooth. Involutory conics $c$ and $i_F(c)$ meet at 2 points $\gamma_1(c), 
\gamma_2(c) \in v(W)$, and $v^{-1}(\gamma_i(c))$ are points of tangency of 
$p_F(c)$ of $W$. There are 12 conics passing through a fixed generic point 
of $X$ and 6 pairs of involutory conics passing through a fixed generic point 
of $v(W)$. 

By definition, $F_0$ is included in $G(2, V_4) \subset P^5=P(\lambda^2(V_4))$. 
It was indicated above (in (0.11)) that (0.4) and (0.5) are true for this 
inclusion and hence Hodge numbers of $F$ and $F_0$ are given in (0.6) - (0.8).

Main theorems for a "double $P^3$" are the following: 
\medskip
{\bf Theorem 2.7.} (Tangent bundle theorem). There exists a locally free sheaf 
$\tau_2$ on $F$ of rank 2 such that for $c \in F$ $P(\tau_2(c))$ is a linear 
envelope of $\gamma_1(c)$, $\gamma_2(c)$. Then $\Omega_F=\tau_2^*$. $\square$ 
\medskip
{\bf Theorem 2.8.} $h^{1,0}(F)=10$. Moreover, Abel - Jacobi map 
$\Phi: \Alb F \to J^3(X)$ is an isomorphism. $\square$ 
\medskip
{\bf Corollary 2.9.} (Geometric interpretation of the tangent bundle theorem). 
Let us consider a map $\rho_1: F \to G(2, H^0(F, \Omega_F)^*)$ given by 
$\Omega_F$, and a map $\rho_2: F \to G(2, S^2(V_4))$ which sends $c \in F$ to 
a linear envelope of $\gamma_1(c)$, $\gamma_2(c) \subset P(S^2(V_4))$. Then 
there exists a natural isomorphism $\iota: H^0(F, \Omega_F)^* \to S^2(V_4)$ 
such that $\iota \circ \rho_1 = \rho_2$. $\square$ 
\medskip
{\bf Theorem 2.10.} $F(X)$ defines $X$ uniquely. $\square$ 
\medskip
{\bf Theorem 2.11.}(Analog of a Riemann theorem). Let us consider a map 
$\Phi_5: \prod^5_{i=1} F \to J^3(X)$ defined as follows: $\Phi_5 (f_1, 
\dots , f_5)= \sum^5_{i=1} \Phi(f_i) - \Phi(i_F(f_i))$. Then its ramification
divisor is a Poincar\'e divisor $\Theta$ on $J^3(X)$. $\square$ 
\medskip
{\bf Theorem 2.12.} $\Theta$ is irreducible, and $\codim _{J^3(X)} \Sing 
\Theta=2 $. $\square$ 
\medskip
So, according Andreotti - Meier criterion, $J^3(X)$ is not the Jacobian of 
a curve, and hence $X$ is not rational. 
\medskip
{\bf Theorem 2.13.} (Global Torelli theorem). $(J^3(X), \Theta)$ defines 
$X$ uniquely. $\square$ 
\medskip
Let us consider now a class of threefolds for which it is easy to describe 
their middle jacobian, namely bundles of conics. Let $X$ be any threefold 
immersed in a projective space. By definition, a structure (denoted by $B$) 
of a bundle of 
conics on $X$ is a map 
$p=p(B): X \to P^2$ such that each fibre of $p$ is isomorphic to a conic. 
Attached to $B$ is a degeneration curve $C=C(B) \subset P^2$ defined as 
follows: $x \in C \iff p^{-1}(x)$ is a reduced conic. All singularities of 
$C$ are ordinary double points, and we have: 
$$x \notin \Sing C \iff p^{-1}(x) \hbox{ is a pair of straight lines}$$
$$x \in \Sing C \iff p^{-1}(x) \hbox{ is a double straight line}$$
Attached is also a two-sheeted covering $p_C=p_C(B): \hat C \to C$ defined
as follows: for $x \in C$ $p_C^{-1}(x) $ is a set of irreducible components
of $p^{-1}(x)$. It has the following property: ramification points of its 
normalization $p_{C_N}: \hat C_N \to C_N$ are exactly inverse images of 
singular points of $C$ under the normalization map $N: C_N \to C$. 

Coverings which satisfy this property are called two-sheeted pseudocoverings. 
Beauville ([26]) defined a notion of Prym variety for the case of 
pseudocoverings. There is a theorem: 
\medskip
{\bf Theorem 2.14.} $J^3(X)= \Pr (\hat C/C)$ ([43] for a smooth $C$, [26] 
for any $C$). $\square$ 
\medskip
There are sufficient conditions that a Prym variety is not the Jacobian of a 
curve, for example 
\medskip
{\bf Proposition 2.15.} Let $B$ be a bundle of conics such that $C(B)$ is 
non-singular and $\deg C(B) \ge 6$. Then $\Pr (\hat C/C)$ is not the Jacobian 
of a curve ([43]). $\square$ 
\medskip
For all Fano threefolds considered above, existence of a structure of bundle 
of conics is known only on threefolds which are birationally equivalent to 
cubics and intersections of three quadrics. If $X$ is a cubic threefold then 
$\deg (C)=5$, and in order to prove that $\Pr (\hat C/C)$ is not the Jacobian 
of a
curve we need to study the sheaf of second order $\sigma$ on $C$ that defines 
$p_C$. If $X$ is an intersection of 3 quadrics then $C$ is the Hesse 
curve $H$ (see Section 5). It is proved in [12] that there is no structure of 
bundle of conics on "double $P^3$". 

For other types of Fano threefolds existence of a structure of a bundle of 
conics is known only for degenerated Fano threefolds. Nevertheless, there 
exists a method of proof of non-rationality of a "generic" Fano threefold if 
there exists its degeneration to
 a non-rational bundle of conics (see [16] for details). This method was used 
in [26] for a proof of non-rationality of a "generic" Fano threefold for all 
types of Fano threefolds of genus one whose non-rationality was not known 
earlier. 
\bigskip
\bigskip
{\bf Section 3. Fano surface as a determinantal variety.}
\bigskip
Here we introduce some notations that will be used throughout the paper.

A sheaf of type $O_{\alpha}(a)\otimes O_{\beta}(b)\otimes \dots$ will 
be denoted by $O(_{\alpha}a, _{\beta}b \dots)$. The inverse image of a sheaf
--- if the map is clear --- will be denoted often by the same symbol as the
sheaf itself. The blowing up map of a variety $Y$ along $X \subset Y$ will 
be denoted 
by $r$ with some subscript, and for $Z \subset Y$ the notation $r^{-1}(Z)$ 
will mean the total inverse image of $Z$, if $Z \subset X$, and the proper 
inverse image if $Z \not\subset X$. $V_{\bullet}$ means a vector space, and 
the 
subscript is its dimension. If $E$ is a locally free sheaf on a variety $Y$, 
then the fibre of the corresponding vector bundle at a point $y \in Y$ will 
be denoted by $E(y)$. $P_Y(E) \to Y$ and $G_Y(k,E) \to Y$ mean the geometric 
projectivization and grassmannization respectively, i.e. 
$P_Y(E) = \hbox{Proj}(\sum _k S^*(E^*))$. Members of tautological exact 
sequences on $G_Y(k,E)$ are denoted as follows: 
\medskip
$$0 \to \tau _{k,E} \to E \to \tau ^* _{n-k,E} \to 0$$ 
$$0 \to \tau _{n-k,E^*} \to E^* \to \tau ^* _{k,E} \to 0$$
(here $n = \hbox{dim }E$). For $t \in G_Y(k,V_n)$    $t_V$ means a 
$k$-dimensional subspace in $V_n$ that corresponds to $t$. In some cases 
a vector space and its projectivization will be denoted by the same symbol. 
\medskip
Here we give without proofs some results about $G(2,5)$ and its hyperplane 
sections (see, for example, [48]).

We identify $G(2,V_n)$ and its image in $P(\lambda ^2(V_n))$ under the 
Pl\"ucker 
embedding. For $n=4$  $G(2,4)$ is a quadric hypersurface 
in $P(\lambda ^2(V_4))$, we shall denote the corresponding element in 
$P(S^2(\lambda ^2(V_4^*)))$ by $Pl(V_4)$. Further we shall consider a space 
$V=V_5$ and the Grassmannian $G=G(2,V_5)$. 

For a conic line $c$ we denote by $\pi (c)$ the plane spanned on $c$. Let 
$c$ be a conic line on $G$. If $\pi (c) \not\subset G$ or if 
$\pi (c) \subset G$ as an $\alpha$-plane (i.e. as a Schubert cycle $\Omega 
(1,4)$), then there exists the only $V_4 \subset V$ such that $c \subset 
G(2,V_4)$. If $\pi (c) \subset G$ as a $\beta$-plane (i.e. $P(V_3)$ 
for some $V_3 \subset V$), then all spaces $V_4$ which contain $V_3$ 
satisfy the condition $c \subset G(2,V_4)$.

Let $H_1$, $H_2$ be hypersurfaces in $\lambda ^2(V)$, i.e. elements of 
$P(\lambda ^2(V^*))$, and $E_2$ a plane that they generate in 
$\lambda ^2(V^*)$. We denote $V_8 = H_1 \cap H_2 \subset \lambda ^2(V)$ and 
$G_4 = G(2,V) \cap H_1 \cap H_2 \subset P(V_8)$. We can associate to an
element $H \in \lambda ^2(V^*)$ a skew-symmetric map $S_H: V \to V^*$, 
and we can associate to an inclusion $E_2 \hookrightarrow \lambda ^2(V^*)$
a map $S: E_2 \otimes V \to V^*$ such that for $H \in E_2$, $v \in V$ we
have $S(H \otimes v)=S_H(v)$. 

For generic $H_1$, $H_2$ \ \ \ $G_4$ is smooth. The condition of smoothness 
of $G_4$ is equivalent to the condition that $\forall H \in P(E_2)$ 
the rank of $S_H$ is 4. Further, taking into consideration (3.5), we 
shall suppose that the pair $H_1$, $H_2$ satisfies this condition. 

We shall often consider the space $P(V^*)$, and we shall consider an element 
$V_4 \in P(V^*)$ as a hyperplane in $V$. The tautological exact sequence on 
$P(V^*)$ is the following: 
$$0 \to \tau_4 \overset{i_V}\to{\hookrightarrow} V \otimes O \to O(1) \to 0$$

For $V_4 \in P(V^*)$ we denote 
$$M(V_4)= \lambda ^2(V_4) \cap H_1 \cap H_2 \subset V_8,$$
$$Q_G(V_4)=G(2,V_4) \cap P(V_8) = G(2,V_4) \cap P(M(V_4))$$

It is clear that $Q_G(V_4)$ is a quadric hypersurface in $P(M(V_4))$. 
\medskip
{\bf Lemma 3.1.} $\forall V_4 \in P(V^*)$ $\dim M(V_4) =4$. 
\medskip
{\bf Proof.} $\dim \lambda ^2(V_4)=6$. If $\dim M(V_4) > 4$, then
there exists an element $H \in P(E_2)$ such that $H \supset
\lambda ^2 (V_4)$. This means that the composition map 
$$ V_4 \hookrightarrow V \overset{S_H} \to{\to} V^* \twoheadrightarrow V^*_4$$
is 0, hence $S_H(V_4) \subset (V/V_4)^*$, i.e. $\rank (S_H) =2$ --- a 
contradiction. $\square$
\medskip
Let us define a sheaf $M$ on $P(V^*)$ as the kernel of the composition 

$$\lambda ^2(\tau_4) \overset{\lambda ^2 (i_V)}\to{\hookrightarrow} 
\lambda ^2(V) \otimes O \to E^*_2 \otimes O$$
where the epimorphism $\lambda ^2(V) \to E^*_2$ is dual to the inclusion 
$E_2 \to \lambda ^2 (V^*)$. It is clear that the fibre of $M$ at $V_4 \in 
P(V^*)$ is $M(V_4)$. Lemma 3.1 implies that $M$ is a locally free sheaf 
of rank 4. 

Let a locally free sheaf $E$ on a variety $Y$ be given. A bundle 
of quadrics $Q$ on $(Y,E)$ is a map that associates to any $y \in Y$ 
a quadric (a fibre of $Q$ in $y$) $Q(y) \subset P(E(y))$ such that 
they form an algebraic family. This means that there exists an 
invertible sheaf $L(Q)$ on $Y$ and a map of sheaves $i(Q): L(Q) \to
S^2(E^*)$ such that
$$\forall y \in Y \; \; \; i(Q)(y) = Q(y) \in P(S^2(E(y))^*)$$

So, there exists a bundle of quadrics $Q_G$ on $(P(V^*), M)$ whose 
fibre at $V_4 \in P(V^*)$ is $Q_G(V_4)$. 
\medskip
{\bf Lemma 3.2.} Let $E$ be a locally free sheaf of rank 4 on $Y$. There 
exists a bundle of quadrics $\Pl (E)$ on $(Y, \lambda ^2(E))$ having 
$\Pl(E)(y)= \Pl(E(y))$. 
Then $L(\Pl(E))=(\det E)^{-1}$. 
\medskip
{\bf Proof.} Let $U_1$, $U_2 \subset Y$ be affine open subsets such that 
$E \vert _{U_i}$ is a free $O_{U_i}$-module. Let $e_i = (e_{i1}, \dots 
, e_{i4})$ $(i = 1,2)$ be its basis, and $e^*_{i1}, \dots , e^*_{i4}$ a basis 
of $E^* \vert _{U_i}$ which is dual to $(e_i)$. Further, let $e^*_{ikl}= 
e^*_{ik} \wedge e^*_{il}$. We have: 
$$\Pl(e_i)=
e^*_{i12} \circ e^*_{i34}
-e^*_{i13} \circ e^*_{i24}
+e^*_{i14} \circ e^*_{i23} \in 
S^2(\lambda ^2 (E^*)) 
\vert _{U_i}$$
is an equation of the Pl\"ucker quadric. Let us consider a map of sheaves
$O \to S^2(\lambda ^2 (E^*)) $ on $U_i$ defined by the element $Pl(e_i)$.
Let on $U_1 \cap U_2$ we have $e_1 = A(e_2)$ for an $A \in GL(E \vert _{
U_1 \cap U_2})$. Then it is easy to check that $\Pl(e_1) = (\det A)^{-1} 
Pl(e_2)$ on $U_1 \cap U_2$. This implies the desired. $\square$
\medskip
Let us apply this lemma to $\tau _4$ on $P(V^*)$. Since $\det (\tau _4)= 
O(-1)$, we have: the bundle $\Pl(\tau _4)$ can be given by a map 
$$i(\Pl(\tau _4)): O(1) \to S^2 (\lambda ^2(\tau _4 ^*))$$
Restriction of $\Pl (\tau _4)$ on $M \subset \lambda ^2(\tau _4)  $ gives us 
the bundle $Q_G$. This means that $L(Q_G)=O(1)$. 

Let us describe now the set of planes on $G_4$. We define an inclusion 
$\psi _1: P(E_2) \to P(V)$ as follows: $\psi _1(H) = \Ker S_H$. Then 
$\im \psi _1$ is a conic line, we denote it by $c_u$ and its plane 
$\pi(c_u)$ by $U = P(U_3)$, where $U_3$ is a subspace of $V$. We have: 
$U^* = P(U^*_3) \subset G_4$ is the only $\beta$-plane on $G_4$. 
A straight line in $P(V^*)$ which is dual to $U$ will be denoted by $l$, 
i.e. $l = P((V/U_3)^*) \subset P(V^*)$, and a point of $l$ is a subspace 
$V_4$ such that $U_3 \subset V_4 \subset V$. There exists an isomorphism 
$\psi: c_u \to l$ such that $\psi (v) =S(E_2 \otimes v) \in P(V^*)$, 
and we denote by $\alpha (v)$ a Schubert cycle $\Omega (v, \psi (v))$.
Then $\forall v \in c_u$ \ \ \ $\alpha (v) \subset G_4$, and all
$\alpha$-planes on $G_4$ are $\alpha (v)$ for some $v \in c_u$. Further,
we have $\alpha (v) \cap U^* \subset G_4$ is a straight line, we 
denote it by $l_g(v)$. It is dual to the point $v$ with respect to the 
plane $U$. All straight lines $l_g(v)$ are tangent to the conic line 
$c^*_u \subset 
U^*$ (the dual to $c_u$). It is clear that $\forall V_4 \in l$ we have 
$Q_G(V_4)= U^* \cup \alpha (\psi ^{-1}(V_4))$, i.e. is a pair of planes. 
\medskip
{\bf (3.3).} The inverse is also true: if for $V_4 \in P(V^*)$ \ \ \ 
$Q_G(V_4)$ is a pair of planes, then it follows easily from the description 
of the set of planes on $G_4$ that $V_4 \in l$. 
\medskip
{\bf (3.4).} There exists a map $b: G_4 \to P^4$ having the 
following properties: it is a birational isomorphism; it is the restriction 
to $G_4$ of the linear projection 
$p_{85}$ whose center is $U^*$, from $P(V_8)$ to $P^4 = P(V_8/U^*_3)$. 
There exist: (a) the unique normcubic $c_3 \subset P^4$ whose linear 
envelope is a $P^3 \subset P^4$, (b) an isomorphism $\psi _2: c_u \to 
c_3$ and (c) an isomorphism $\tilde b: 
(\tilde G_4)_{U^*} 
\to 
(\tilde P^4)_{c_3}$ 
- a desingularization of $b$ such that the following diagram is commutative:
$$\matrix
(\tilde G_4)_{U^*}     &
\overset{\tilde b}\to{\to}& 
(\tilde P^4)_{c_3}\\
\downarrow & & \downarrow \\ 
G_4 & \overset{b}\to{\to} & P^4
\endmatrix $$
where left and right vertical maps (denoted by $r_G$ and $r_P$ 
respectively) are blowings up of $G_4$ and $P^4$ respectively. Further, for 
$v \in c_u$ we have: 
$r^*_G(\alpha(v))$ 
is a plane, and 
$$\tilde b(r^*_G(\alpha(v)))=r^*_P(\psi_2(v))$$ 
For $t \in U^* \subset G$ we denote $\{v_1, v_2\}=t_V\cap c_u$. Then 
$<\psi _2(v_1), \psi _2(v_2)>$ 
is a bisecante of $c_3$ in $P^3$, and 
$$\tilde b(r^*_G(t))=r^*_P(<\psi _2(v_1), \psi _2(v_2)>)$$ 
For $V_4 \in P(V^*)-l$ the map 
$b \vert _{Q_G(V_4)}
: Q_G(V_4) \to P^2$ 
can be described as follows. $Q_G(V_4) \cap U^* = \{t\}$ (one point); 
it is clear that $t_V = V_4 \cap U_3$. Then
$b \vert _{Q_G(V_4)}   $ is the blowing up of $t$ into to a bisecante of
$c_3$ and the blowing down of $t_1$, $t_2$ into intersection points of this
bisecante with $c_3$, where $t_1$, $t_2$ are straight lines on $Q_G(V_4)$, 
which contain $t$. 

Let us describe now conic lines on $X$. Let $X = G_4 \cap \Omega$, where 
$\Omega$ is a quadric hypersurface on $P(V_8)$. 
\medskip
{\bf Lemma 3.5.} $G_4$ is smooth.
\medskip
{\bf Proof.} If $G_4$ is not smooth, then $\exists H \in P(E_2)$ such that 
$\rank S_H = 2$. Then $\Sing (G \cap H)= G(2, \Ker S_H)= P^2$, hence 
$\Sing X \supset P^2 \cap H_1 \cap \Omega$. This means that $\Sing X$ 
is not empty - a contradiction. $\square$
\medskip
We denote by $F_c = F_c(X)$ the set of conics on $X$; it is an algebraic 
variety which is called the Fano surface of $X$ (later we shall see that 
$F_c$ is really a surface). Taking into consideration that for a generic 
conic $c \subset G$ $\exists ! V_4 \in P(V^*)$ such that 
$c \subset G(2,V_4)$
, we define a surface $F=F(X)$ (which we shall call the Fano 
surface as well) as a set of pairs $(c \in F_c, V_4 \in P(V^*))$ such that
$c \subset G(2,V_4)$.
There are natural projections 
$$r_F: F \to F_c, \ \ \ \tilde \phi: F \to P(V^*)$$
There exists a bundle of quadrics $Q_{\Omega}$ on $(P(V^*), M)$ whose fibre at 
$V_4 \in P(V^*)$ is $Q_{\Omega} (V_4)=P(M(V_4))\cap \Omega$. It is 
clear that 
\medskip
(a) $L(Q_{\Omega})=O$; 
\medskip
(b) $X \cap G(2, V_4)=Q_G(V_4) \cap Q_{\Omega} (V_4)$. 
\medskip
Further, we have $\forall V_4 \in P(V^*)$ \ \ \ 
$\dim Q_G(V_4)\cap Q_{\Omega}(V_4)=1$, 
because smoothness of $X$ implies that $\Pic (X)$ is generated 
by $O_{P(V_8)}(1)\vert _X$, so the degree of any surface on $X$ is 
a multiple of 10 $(= \deg X)$. But if 
$Q_G(V_4) \cap Q_{\Omega} (V_4)$   contains a surface, then its degree 
is 1 or 2. 

It is clear that $(c, V_4) \in F \iff c \subset 
Q_G(V_4) \cap Q_{\Omega} (V_4)$. Since 
$Q_G(V_4) \cap Q_{\Omega} (V_4)$ is a curve of degree (2,2) on a quadric 
surface, and $c$ is a curve of degree (1,1), then 
$Q_G(V_4) \cap Q_{\Omega} (V_4) -c$ is also a conic line. We denote it 
by $c'$. So, there exists an involution $i_F$ on $F$ defined as follows: 
$i_F(c,V_4)= (c',V_4)$. Let $F_0$ be the quotient surface $F/i_F$ and
$p_F: F \to F_0$ the corresponding two-sheeted covering. A point
$f_0 \in F_0$ will be denoted as follows:
$f_0 = (c,c',V_4)$, where $f_0=p_F(c,V_4)=p_F(c',V_4)$. It is clear that 
$c$ and $c'$ have 2 (possible coinciding) intersection points which we 
denote by 
$\gamma _1(f_0)$ and $\gamma _2(f_0)$ (or 
$\gamma _1(c)$ and $\gamma _2(c)$). Their linear envelope 
(i.e. the straight line 
$<\gamma _1(c), \gamma _2(c)>$) will be denoted by $l_8(c)$ (or $l_8(f_0)$). 
\medskip
{\bf (3.6).} There exist smooth threefolds $X$ which contain a pair of 
involutory conics 
$c$ and $c'$ such that both of them is a pair of straight lines such that 
$l_8(c)$ is their common component. Considerations analogous to ones 
used in the proof of Lemma 3.7 show that the set of $X$ having this property 
has codimension 1 in the set of all possible $X$. These $X$ will not be 
considered in the present paper, although most likely main results which we 
have got for a "generic" $X$, are true for them as well 
(for example, smoothness of $F$). 

It is clear that there exists a map $\phi: F_0 \to P(V^*)$ satisfying a 
condition $\phi \circ p_F= \tilde \phi$. Further, it is clear that 
$$(c, V_4) \in F \iff \pi(c) \cup \pi(c') \in <
Q_G(V_4), Q_{\Omega}(V_4)>_{P(S^2(M(V_4))^*)}$$
i.e. some linear combination of
$Q_G(V_4)$, $Q_{\Omega}(V_4)$ is a pair of planes.  So, let us define $Z$ 
as a variety of pairs $(V_4, Q)$ where $V_4 \in P(V^*)$, $Q \in 
<Q_G(V_4), Q_{\Omega}(V_4)>$ together with a natural projection 
$\eta: Z \to P(V^*)$. Since $L(Q_G)=O(1)$, $L(Q_{\Omega})=O$, we have 
$$Z=P_{P(V^*)}(O+O(1))$$
There is a tautological exact sequence on $Z$: 
$$0 \to O_{\eta}(-1) \overset{i_{\eta}}\to{\to} O+O(1) 
\to O(1, _{\eta}1) \to 0$$
as well as the dual sequence. Further, $Z$ is the blowing up of $P^5=P(V_6)$ 
at a point $t \in P^5$; we denote the corresponding map by $r_Z: Z \to P^5$. 
It is clear that $V_6/t=V^*$ and $r_Z^*(O_{P^5}(1))=O(1,_{\eta}1)$. 
There are divisors $D_G$, $D_{\Omega}$ on $Z$ which are defined as follows: 
a point of $Z$ can be considered as a pair $(V_4, Q)$; we have: 
$$(V_4, Q) \in D_G \iff Q = Q_G(V_4);$$
$$(V_4, Q) \in D_{\Omega} \iff Q = Q_{\Omega}(V_4)$$
It is clear that $D_G = r^*_Z(t)$ is the exceptional divisor on $Z$, and 
$r_Z(D_{\Omega})$ is a hypersurface in $P^5$. 

Let $(V_4, Q) \in Z$; we associate it the quadric $Q \subset P(M(V_4))$. 
This gives us a bundle of quadrics (denoted by $T$) on $(Z,M)$. It is clear
that $L(T)=O_{\eta}(-1)$ and
$$ i(T)=(i(Q_G) \oplus i(Q_{\Omega})) \circ i_{\eta} : O_{\eta}(-1) \to
S^2(M^*)$$

There exists an inclusion $\bar \phi :F_0 \to Z$ defined as follows: 
$\bar \phi (c,c',V_4) = (V_4, \pi(c) \cup \pi(c'))$. It is clear that 
$z \in \im \bar \phi \iff T(z)$ is a pair of planes, i.e. a rank 2 
quadric surface. This means that $F_0 \hookrightarrow Z$ is the 
second determinantal of $T$.

Let us recall the definition and main properties of 
determinantal varieties. Let $E$ be a locally free sheaf of rank $n$ 
on an algebraic variety $Y$, $Q$ a bundle of quadrics on $(Y,E)$ and 
$i(Q): L^* \to S^2(E^*)$ the corresponding map. For any $y \in Y$ 
we can consider the fibre $i(Q)(y)$ as a symmetric map from $E(y) $ 
to $E(y)^*$. The determinantal $D_k=D_k(Q)$ of $Q$ is defined as 
a set of points $y \in Y$ such that $\dim \im i(Q)(y) \le k$. We define also 
a variety $\tilde D_k = \tilde D_k(Q)$ as a set of pairs $(y, V_k)$ such that 
$y \in D_k$ and $\im i(Q)(y) \subset V_k \subset E(y)^*$. The natural 
projection $\pi _k: \tilde D_k \to D_k$ is an isomorphism outside of 
$\pi_k^{-1}(D_{k-1})$. There is also a tautological inclusion 
$\tilde D_k \to G_y(k,E^*)$, and we have a commutative diagram:

$$\matrix \tilde D_k & \hookrightarrow & G_Y(k, E^*) \\
\downarrow & & \downarrow \\
D_k & \hookrightarrow & Y \endmatrix $$

\medskip
There is a sheaf $C=C(k)=\Coker (S^2(\tau_{k,E^*}) \to S^2(E^*))$ on 
$G_Y(k,E^*)$. 
The sheaf $L \otimes C$ has a distinguished section coming from the map 
$$O \overset{i(Q) \otimes L}\to{\to} L \otimes S^2(E^*) \twoheadrightarrow L 
\otimes C $$
It is clear that $\tilde D_k$ is the set of zeros of this section. $Q$ is 
called $k$-regular, if for $i=k$ and $i=k-1$ 
$$\Codim_{G_Y(i,E^*)}\tilde D_i = \dim L \otimes C(i)$$
(if $D_i \ne \emptyset$). If so, we have $\Codim _Y D_k = 
\frac{(n-k)(n-k+1)}{2}$. Further, for a $k$-regular bundle $Q$ we have: 
$\bar D_{k-1}=\pi^{-1}_k(D_{k-1})$ is a divisor in $\tilde D_k$, which is 
the determinantal of rank $k-1$ of a bundle of quadrics 
$L^* \to S^2(\tau_k)$, which is defined naturally on $\tilde D_k$. 
Porteus formula ([7]) implies that 
$$O_{\tilde D_k}(\bar D_{k-1})=L^k\otimes O_G(-2)$$
Particularly, if $D_{k-1}=\emptyset$, then 
$\bar D_{k-1}=\emptyset$ as well, and $L^k = O_G(2)$ in $\Pic (D_k)$.

Let us consider the case of even $k$.  For a quadric $Q$ of rank $k$
the set of projective spaces of the maximal possible dimension which are 
contained in $Q$, has 2 connected components. We denote this set of 
connected components by $\Gamma (Q)$, and we define a two-sheeted covering 
$p_{D_k}: \hat {\tilde D_k} \to \tilde D_k$ as follows: $p_{D_k}^{-1}(y)=
\Gamma (Q(y))$ for any $y \in  \tilde D_k - \bar D_{k-1}$. We have: 
$\bar D_{k-1}$ is the ramification divisor of $p_{D_k}$, and a sheaf 
that corresponds to $p_{D_k}$ is  $L^{\frac{k}{2}}\otimes O_G(-1)$.

Let us consider the case $k=2$. For $y \in D_2$ \ \ \ $Q(y)$ is a pair 
of hyperplanes in $P(E(y))$, and a choice of a point in $p_{D_2}^{-1}(y)$ 
is a choice of one of these hyperplanes. So, there exists an inclusion 
$\hat {\tilde D_2} \to P_Y(E)$, and there exist sheaves $\tau_{n-1,E}$ and 
$O_{\pi}(1)$ on $\hat {\tilde D_2}$ (here $\pi : P_Y(E) \to Y$ is a 
projection). It is clear that $\hat {\tilde D_2}$ is the set of zeros 
of a section of $L \otimes S^2(\tau^*_{n-1,E} )$ (here $L \otimes 
S^2(\tau^*_{n-1,E} )$ is a sheaf on $P_Y(E^*)$ ). 

\medskip
Let us return to study of the Fano surface. As we have seen, $F_0=D_2(T)$. 
\medskip
{\bf Lemma 3.7.} For a generic $X$ we have $D_1(T)=\emptyset$.
\medskip
{\bf Proof.} It is clear that 
$$D_1(T) \ne \emptyset \iff \exists V_4 \in P(V^*), \; \exists V_3 \subset 
M(V_4)$$ such that the double plane $2 P(V_3) \in 
<Q_G(V_4), Q_{\Omega}(V_4)>$. Let us fix $G_4$. For given $V_4 \in P(V^*)$, 
$ V_3 \subset M(V_4)$ \ \ \ $Q_G(V_4)$ is defined uniquely, and the set of 
quadrics $Q_{\Omega}(V_4)$ such that 
$2 P(V_3) \in <Q_G(V_4), Q_{\Omega}(V_4)>$ has codimension 8 in the space 
of all quadrics in $P(V_8)$. This means that the set of 
quadrics $\Omega$ such that 
$2 P(V_3) \in <Q_G(V_4), Q_{\Omega}(V_4)>$ has codimension 8 in the space 
of all quadrics in $P(V_8)$. But the set of pairs $\{V_4 \in P(V^*)$, 
$ V_3 \subset M(V_4)\}$ is 7-dimensional, hence the set of 
quadrics $\Omega$ such that for $X= G_4 \cap \Omega$ \ \ \ $D_1(T)\ne 
\emptyset$ 
is a union of a 7-dimensional set of subvarieties of codimension 8, 
hence has codimension $\ge 1$.  $\square$
\medskip
(4.2) (see below) implies that $\dim F=2$. This means that a condition 
$D_1(T)= \emptyset$ for $X$ is equivalent to a condition of 2-regularity 
of the bundle $T$. Further we shall consider only $X$ that satisfy this
condition.

Properties of determinantal varieties given above imply the following 
\medskip
{\bf Proposition 3.8.} There exist inclusions 
$$F_0 \to G_Z(2,M^*), \; \; F \to P_Z(M^*)$$
hence tautological sheaves on $G_Z(2,M^*)$ and $P_Z(M^*)$ can be restricted 
on $F_0$ and $F$ respectively. Further, $F_0$ is the set of zeros of a 
section of $O_{\eta}(1) \otimes C(2)$ on $G_Z(2,M^*)$, and $F$ is 
the set of zeros of a section of $O_{\eta}(1) \otimes S^2(\tau^*_{3,M})$ on 
$P_Z(M^*)$. We have equalities of sheaves $O_G(2)=O_{\eta}(2)$ (respectively 
 $O_G(1)=O_{\eta}(1)$) on $F_0$ (respectively on $F$). The covering 
$p_F: F \to F_0$ is non-ramified, and the corresponding sheaf is $\sigma =
O(_{\eta}1, _G -1)$. For $f_0 \in F_0$ \ \ \ $l_8(f_0)=P(\tau_{2,M}(f_0))$. 
$\square$
\medskip
If for $V_4 \in P(V^*)$ \ \ \ $Q_G(V_4)$ or $Q_{\Omega }(V_4)$ is already 
a pair of planes then it is clear that $V_4 \in \phi(F_0)$. We denote the set 
of these $V_4$ (and the set of their $\phi$-inverse images on $F_0$) by 
$l_G$ and $l_{\Omega}$ respectively. They are curves on $F_0$, and it is 
clear that $l_G=\bar \phi^*(D_G)$,  $l_{\Omega}=\bar \phi^*(D_{\Omega})$. 
$l_{\Omega}$ is not defined by $X$ uniquely, it depends on a choice of
$\Omega$ such that $X=G_4 \cap \Omega$. (3.3) implies that $l_G=l$.
It is clear that $p_F^{-1}(l)$ is a pair of straight lines, we denote them by
$l_1$ and $l_2$ by such a way that 
$$(c, V_4) \in l_1 \iff \{c=U^* \cap \Omega, \; V_4 \in l\}$$
$$(c, V_4) \in l_2 \iff \{c=\alpha(\psi^{-1}(V_4)) \cap \Omega, \; 
V_4 \in l\}$$
 
It is clear that $p_F \vert _{l_i}: l_i \to l$ are isomorphisms ($i = 1,2$)
and that $r_F(l_1)$ is a point $c_{\Omega} \in F_c$; the corresponding conic 
line $U^* \cap \Omega$ on $X$ is also $c_{\Omega}$. Using properties 
of conic lines on $G$ we see that for $c \in F_c$ \ \ \ $r^{-1}_F(c)$ is 
either 
a point or a straight line on $F$. If $r^{-1}_F(c)$ is a straight line, 
then $\pi(c)$ is a $\beta$-plane on $G$. In this case $\pi(c) \cap H_1 
\cap H_2 \supset c$, hence $\pi(c) \subset G_4$ and $\pi(c) = U^*$, 
$c=c_{\Omega}$. So, $c_{\Omega}$ is the  only point of $F_c$ such that 
$r^{-1}_F(c_{\Omega})$ is not a point. This means that $r_F: F \to F_c$ 
is the blowing up of $c_{\Omega}$. 
\medskip
{\bf Proposition 3.9.} $O_{F_0}(l)=O_{\eta}(1)$,    
$O_{F_0}(l_{\Omega})=O(1,_{\eta}1)$. 
\medskip
{\bf Proof.} This important proposition (it implies smoothness of $F_c$ in
$c_{\Omega}$, see below) can be proved by different ways. Here
is one of them. Let us consider maps of sheaves on $G_Z(2,M^*)$: 
$$0 \to O_{\eta}(-1) \overset{i_{\eta}}\to{\hookrightarrow} O \oplus O(1) 
\overset{i_{1}}\to{\twoheadrightarrow} O(1,_{\eta}1)   \to 0 \eqno{(3.10)}$$

$$O(1) \overset{i_{2}}\to{\hookrightarrow} O \oplus O(1) 
\overset{i_{3}}\to{\to} S^2(M^*) \overset{i_{4}}\to{\twoheadrightarrow} C(2) 
\eqno{(3.11)}$$
where $i_3 = i(Q_{\Omega}) \oplus i(Q_G)$. 
By definition, for $f_0 \in F_0$ \ \ \ $i_4 \circ i_3 \circ i_{\eta}(f_0)=0$, 
and exactness of (3.10) implies that there exists an inclusion $i_5: 
O(1,_{\eta}1)  \to C(2)$ on $F_0$ satisfying a condition 
$i_5 \circ i_1 = i_4 \circ i_3$. It is clear that 
$$ f_0 \in l \iff \im (i_3 \circ i_2(f_0)) = \im (i_3 \circ i_{\eta}(f_0))$$
This implies that for $f_0 \in F_0$ \ \ \ \ \ 
$f_0 \in l \iff i_1 \circ i_2(f_0)=0$, i.e. $l$ is the set of zeros of 
$i_1 \circ i_2$ on $F_0$. This means that $O_{F_0}(l)=O_{\eta}(1)$. 
Replacing in (3.11) $O(1)$ by $O$, we get analogously that 
$O_{F_0}(l_{\Omega})=O(1,_{\eta}1)$.  $\square$
\medskip
We denote the projection $P_Z(M^*) \to Z$ by $\pi$, and its
inversible sheaf by $O_{\pi}(1)$. Further, we denote $c_1(O(1))=H$,
$c_1(O_{\eta}(1))=E$, $c_1(O_{\pi}(1))=P$. Chow ring $A(P_Z(M^*))$ is
generated by $H$, $E$, $P$ satisfying relations 
$$H^5=0, \; E^2=-EH, \; P^4=-3P^3H-5P^2H^2-5PH^3.$$
\medskip
{\bf Proposition 3.13.} Intersection indices on $F_0$ are: $EH=1$, $E^2=-1$.
\medskip
{\bf Proof.} Equality $EH=1$ is obvious, because $\cl (l)=E$, $l$ is 
a straight line in $P(V^*)$, and $H$ is a class of a hyperplane section 
in $P(V^*)$. Equality $E^2=-EH$ is true even on $Z$ and moreover on $F_0$. 
$\square$
\medskip
{\bf Corollary 3.14.} $c_{\Omega} $ is a non-singular point on $F_c$. 
\medskip
{\bf Proof.} $<E^2>_F=-2$. $i_F$ interchanges $l_1$ and $l_2$ and hence 
$<l_1^2>_F = <l_2^2>_F$. Further, $l_1 \cap l_2 = \emptyset$, hence 
$<l_1^2>_F = <l_2^2>_F=-1$. This implies that $c_{\Omega} $ is a 
non-singular point on $F_c$. $\square$
\medskip
Maps of blowing down of $l_1 \cup l_2$ on $F$ and of $l$ on $F_0$ are 
denoted by $r: F \to F_m$ and $r_0: F_0 \to F_{0m}$ respectively. 
There exists an involution $i_{F_m}: F_m \to F_m$ and a projection
$p_{F_m}: F_m \to F_{0m}$ making the following diagram commutative:
$$\matrix & & F & & \overset{i_F}\to{\to} & & F \\
& \overset{r_F}\to{\swarrow} & & \overset{p_F}\to{\searrow} & & 
\overset{p_F}\to{\swarrow} & \\ F_c & & \downarrow & & F_0 & & \downarrow \\
& \searrow & & & \downarrow & & \\
& & F_m & & \overset{i_{F_m}}\to{\to} & & F_m \\
& & & \overset{p_{F_m}}\to{\searrow} & & \overset{p_{F_m}}\to{\swarrow} & \\
& & & & F_{0m} & & \endmatrix $$ 
where vertical maps are $r$, $r_0$, $r$ respectively.

Let us consider a map $r_Z \circ \bar \phi: F_0 \to P^5$. Since $l \subset 
D_{\Omega} = r^*_Z(t)$, we have: $r_Z \circ \bar \phi $ factors through 
$\bar \phi_m: F_{0m} \to P^5$, hence we have a commutative diagram 

$$\matrix F_0 & \overset{\bar \phi}\to{\to} & Z \\
\downarrow & & \downarrow \\
F_{0m} & \overset{\bar \phi_m}\to{\to} & P^5 \endmatrix $$
where vertical maps are $r_0$ and $r_Z$ respectively. Further, we have 
$l=F_0 \underset{Z}\to{\times} D_G$. This implies 
\medskip
{\bf Proposition 3.15.} $\bar \phi_m$ is regular at $c_{0\Omega}$, and 
$\bar \phi_m (c_{0\Omega})=t$. $\square$
\medskip
Since $l_{\Omega}= \bar \phi^* (D_{\Omega})$ and $r_Z(D_{\Omega})$ is a
hyperplane in $P^5$, we have that $\bar \phi_m (l_{\Omega})  $ is a 
hyperplane section of $\bar \phi_m (F_{0m})$. We denote by ${P^5}'$ the set 
of quadrics in $P(V_8)$ that contain $X$; it is the linear envelope of 
$\Omega$ and $q(P(V))$ (the definition of $q(P(V))$ is given in the 
beginning of Section 5). To each $\Omega ' \in {P^5}' - q(P(V))$ we 
can associate its own straight line $l_{\Omega '}$ and hyperplane 
$r_Z(D_{\Omega '})$ in $P^5$. It is easy to see that the map 
$\Omega ' \to r_Z(D_{\Omega '})$  gives us a natural isomorphism ${P^5}' \to 
(P^5)^*$ such that quadrics in $q(P(V))$ correspond to hyperplanes that 
contain $t$. 
\medskip
{\bf Lemma 3.16.} $\deg \phi(l_{\Omega})=40$. 
\medskip
{\bf Proof.} We have $ l_{\Omega}= D_Z(Q_{\Omega})$ on $P(V^*)$. According 
the Porteus formula ([7]), we have: $\deg \phi(l_{\Omega})=4(c_1(M^*)
c_2(M^*)-c_3(M^*))$. Further, $M=8O-5O(1)+O(2)$ in $K_0(P(V^*))$, hence 
$c_t(M^*)=1+3H+5H^2+5H^3$, and $\deg \phi(l_{\Omega})=40$.  $\square$
\medskip
Since $\cl_{F_0}(l_{\Omega})=H+E$, we have $\deg l_{\Omega}=<H+E,H>_{F_0}$, 
and $<H^2>_{F_0}= \deg \phi (F_0) =39$. It is easy to see that
$\deg _{P^5}(\bar \phi_m (F_{0m}))=40$ and
$\cl _Z(\bar \phi (F_0)) = 40 H^3+39 H^2E$ (this formula can be got by
application of the Porteus formula for $T$, but calculations will be longer 
than in Lemma 3.16). 

{\bf Remark.} $<H^2>_{F_0}$ is equal to the number of conics passing through 
a fixed generic point on $X$. Really, let $x \in X$, $x_V$  
the corresponding straight line in $P(V)$ and $\tilde x_V=P(V/x_V)^* \subset 
P(V^*)$ the dual plane in $P(V^*)$, i.e. the set of $P^3$ that contain 
$x_V$. Let $f =(c, V_4)\in F$ and $x \in c$. Then $x_V \subset P(V_4)$, i.e. 
$V_4 \in \tilde x_V$. But $V_4 \in \phi (F_0)$, i.e. if $x \in c$, then 
$\tilde \phi (f) \in \tilde x_V \cap \phi (F_0)$. Clearly the inverse 
is also true, so the number of conics passing through $x$ is equal to 
the number of intersection points of $\phi (V_0)$ and $\tilde x_V $. So, 
there are 39 conics passing through a generic point of $X$. 

We have an equality in $K_0(P_Z(M^*))$:

$$\Omega _F = \Omega _{P_Z(M^*)}-O_{\eta}(-1)\otimes S^2(\tau _{3,M})
\eqno (3.17)$$
We deduce easily from this formula that $\omega_F=O(3,_{\eta}4)$. Further, 
we have an exact sequence (where $i$ is an inclusion $l_1 \cup l_2 \to F$)
$$0 \to r^*(\Omega _{F_m}) \to \Omega _F \to i_*(\Omega _{l_1 \cup l_2} )
\to 0 \eqno(3.18)$$
It implies that $ r^*(\omega _{F_m})=O(3,_{\eta}3)$, i.e. $\omega _{F_m}=
p^*_{F_m}(\bar \phi ^*_m (O_{P^5}(3)))$. Let $\sigma_m$ be a sheaf of order 2 
that corresponds to a two-sheeted covering $p_{F_m}$; it is clear that 
$r_0^*(\sigma _m)=\sigma = O(_G-1, _{\eta}1)$. This means that $\omega _{F_0}=
O(3,_G 1, _{\eta} 3)$ and $\omega _{F_{0m}}=
\bar \phi ^*_m (O_{P^5}(3))\otimes \sigma _m$. Since $\deg \bar \phi _m
(F_{0m})=40$, we get that $c_1(\Omega _{F_{0m}} )^2 = 360 $ and 
$c_1(\Omega _{F_{m}} )^2 = 720$. 

To calculate $c_2(\Omega _{F})$, it is necessary to find $\cl (F)$ in 
$A(P_Z(M^*))$, which is equal $c_6(O_{\eta}(1)\otimes S^2(\tau ^*_{3,M}))$
(see (3.8)). A calculation gives the desired result: 

$$\cl (F)= 80 P^3H^3+240P^2H^4+78P^3H^2E+235P^2H^3E$$
Further, we get (taking into consideration (3.12)) intersection indices of 
divisors on $F$: 
$$H^2=78; \; P^2=-13; \; PH=1; \; EH=2; \; EP=-1$$
We deduce from (3.17) after some calculations: 
$$c_2(\Omega _{F})=4H^2-4P^2+13EH+4EP=386$$
and using (3.18) we get
$$c_2(\Omega _{F_{m}} )=384$$

We define a surface $W \subset X$ as a set of intersection points of 
involutory conics, i.e. 
$$W=\bigcup _{f_0 \in F_0} (\gamma_1(f_0) \cup \gamma_2(f_0))$$
Further, we define a surface $\bar W$ which makes the following diagram 
commutative: 
$$\matrix \bar W & \hookrightarrow & P_{F_0}(\tau_{2,M}) \\
                 & \overset{p_W}\to{\searrow} & \downarrow \\
                 &                           & F_0       \endmatrix $$
where the vertical map is $p$, and $p_W$ is a two-sheeted covering of $F_0$ 
such that for $f_0 \in F_0$ we have 
$$p_W^{-1}(f_0)=\{ \gamma_1(f_0), \gamma_2(f_0)\} \subset p^{-1}(f_0) 
=l_8(f_0)$$
It is clear that there exists a map $\xi : P_{F_0}(\tau_{2,M}) \to P(V_8)$, 
that corresponds to an inclusion of sheaves $\tau_{2,M} \hookrightarrow M 
\hookrightarrow V_8 \otimes O$ on $F_0$. We have: $\xi \vert _W: 
\bar W \to W$ is an epimorphism and $ \xi^*(O_{P(V_8)}(1))=O_P(1)$. 
\medskip
{\bf Lemma 3.19.} Let $f_a=(c_a, c_a', V_{4a}) \in F_0-l$, 
$f_b=(c_b, c_b', V_{4b}) \in F_0$ and $\gamma_1(f_a)=\gamma_1(f_b)$. Then
$j=j(f_a, f_b)=P(M(V_{4a}) \cap M(V_{4b}))$ is a straight line in $G_4$
which is tangent to $\Omega$.
\medskip
{\bf Proof.} Since $P(V_{4a} \cap V_{4b})$  is a plane, we have 
$P(\lambda ^2(V_{4a} \cap V_{4b})) \subset G$ and 
$j=P(\lambda ^2(V_{4a} \cap V_{4b})) \cap H_1 \cap H_2$ is a linear subspace 
in $G_4$. Since $f_a \not\in l$, $j$ is not a plane. If $j$ is a point, then 
$j=\gamma_1(f_a)=\gamma_1(f_b) $ and $T_{P(M(V_{4a}))}(j)$, 
$T_{P(M(V_{4b}))}(j)$ are linearly independent in $T_{P(V_{8})}(j)$. But 
$\dim T_X(j) \cap  T_{P(M(V_{4a}))}(j)  =2$, because this space coincide 
with $T_{Q_G(V_{4a})}(j)$ or $T_{Q_{\Omega}(V_{4a})}(j)$. Analogously, 
$\dim T_X(j) \cap  T_{P(M(V_{4b}))}(j)  =2$, and linear independence of these 
spaces implies that $\dim T_X(j) \ge 4$ that contradicts to non-singularity. 
This means that $j$ is a straight line. It is easy to see that it is 
tangential to $\Omega $ at $\gamma_1(f_a)$. $\square$
\medskip
{\bf Lemma 3.20.} Let $\gamma$ be a straight line on $G_4$ which is tangential 
to $\Omega$ at a point $t_0$. Let $\pi(j)=\cup_{t \in j} t_V$ be 
a plane in $P(V)$ and $\widetilde{\pi(j)}$ the dual straight line in 
$P(V^*)$, i.e. the set of subspaces $P(V_4) \in P(V^*)$ that contain $\pi(j)$. 
If $(c, c', V_4) \in F_0$ and $V_4 \in \widetilde{\pi(j)}$, then $t_0 \in 
\{\gamma_1(c), \gamma_2(c)\}$ or $j \subset X$. 
\medskip
{\bf Proof.} The above conditions imply that $j \subset Q_G(V_4)$. It is
clear that $j \cap \Omega = \{j \cap c, j \cap c'\}$. Since $j$ is tangential
to $\Omega$ at a point $t_0$, then: 
\medskip
If $j \not\subset X$, then $t_0=j \cap c=j \cap c'$, 
\medskip
and hence $t_0 \in c \cap c'$. $\square$
\medskip
{\bf Proposition 3.21.} For a generic point $t \in W$ \ \ \ 
$(\xi \vert _{\bar W})^{-1}(t)$ is one point. 
\medskip
{\bf Proof.} Let for a point $t \in W$ \ \ \ 
$(\xi \vert _{\bar W})^{-1}(t)$ is more than one point. $f_a, \; f_b \in p_W
(\xi \vert _{\bar W})^{-1}(t)$ and $j=j(f_a,  f_b)$. Then a straight line 
$\widetilde{\pi(j)}$ meets $\phi(F_0)$ at points $\phi(f_a), \phi (f_b)$. 
But a variety $B \subset G(2, V^*)$ of bisecants of $\phi(F_0)$ is 
4-dimensional, while a variety $B'$ of straight lines on $G_4$ which are 
tangential to $\Omega$ is 3-dimensional. Let $j$ be such a line. We can 
associate it a straight line $\widetilde{\pi(j)}$ in $P(V^*)$. These lines 
form a threefold $B' \subset G(2, V^*)$. Intersection of $B$ and $B'$ in 
$G(2, V^*)$ is a curve, hence the set of points $t \in W$ such that   
$(\xi \vert _{\bar W})^{-1}(t)$ is more than one point is 1-dimensional.
$\square$
\medskip
{\bf Proposition 3.22.} $O_X(W)=O_X(21)$. 
\medskip
{\bf Proof.} $\forall f_0 \in F_0$ \ \ \ $\{\gamma_1(f_0),\gamma_2(f_0) \}$ is 
a quadric in $P(\tau_{2,M}(f_0))$, and these quadrics form a bundle $\Gamma$
on $(F_0, \tau_{2,M})$. Let us find $L(\Gamma)$. There are maps of sheaves 
on $F_0$: 
$$O_{\eta}(-1) \overset{i_{\eta}}\to{\to} O \oplus O(1) \to S^2(M^*) \to 
S^2(\tau_{2,M}^*)$$
Their composition is 0, hence there exists a map 
$$ O(1, _{\eta}1)=\Coker i_{\eta} \to S^2(M^*) \to 
S^2(\tau_{2,M}^*)$$
It is easy to see that this map is $i(\Gamma)$ and $L(\Gamma)=O(1, _{\eta}1)$.
This implies that 
$$O_{P(\tau_{2,M})}(\bar W)=O(_P2,-1, _{\eta}-1)\eqno(3.23)$$
We denote $c_1(O_{P(V_8)}(1))$ by $R$. We have 
$$\deg W =<W,R,R>_{P(V_8)}=<\bar W,R,R>_{P(\tau_{2,M})}$$
In $A(P_{F_0}(\tau_{2,M}))$ we have an equality $R^2=-c_1R-c_2$, where $c_1
=c_1(\tau_{2,M}) = -3H-E$, $c_2
=c_2(\tau_{2,M}) $. (3.23) implies that $\cl \bar W = 2R-H-E$. This means 
that
$$<\bar W,R,R>_{P(\tau_{2,M})} = (2R^3-HR^2-ER^2)_{P(\tau_{2,M})}$$ $$=
((2c^2_1-2c_2+Hc_1+Ec_1)R+2c_1c_2+Hc_2+Ec_2)_{P(\tau_{2,M})} $$ $$=
(2c^2_1-2c_2+Hc_1+Ec_1)_{F_0}=210.$$ Since $\deg X= 10$, we have 
$O_X(W)=O_X(21)$.   $\square$
\medskip
Let us compare properties of a threefold $X$ under consideration and 
another threefold $X'$ which is a two-sheeted covering of $P^3$ with 
ramification
at a quartic surface $W$ (see [11], [12], [49]). $c_1(\Omega)^2$ and 
$c_2(\Omega) $ of their Fano surfaces $F(X')$ and $F_m(X)$ coincide. 
Since irregularities coincide as well: 
$$h^{1,0}(F_m(X)) = h^{1,0}(F(X')) = 10, \; h^{1,0}(F_{0m}(X)) = 
h^{1,0}(F_0(X')) = 0$$
we have that their Hodge numbers (as well as their decomposition on 
symmetric and antisymmetric 
part with respect to $i_{F_m(X)}$, $i_{F(X')}$) also coincide. 
Moreover, we have the following properties of inclusions $F_m(X) 
\hookrightarrow P^5$, $F(X') \hookrightarrow P^5$: $\deg F_m(X) =
\deg F(X') =40$, $\omega _{F_{0m}(X)}=\omega_{F_0(X')}=O_{p^5}(3) 
\otimes \sigma$. Nevertheless, some properties of $F_m(X)$ and $F(X')$ 
are different. For example, we have for $X'$: 

(a) There 12 conics passing through a generic point $x \in X'$;

(b) $\xi \vert_{\bar W}: \bar W \to W$ is a six-sheeted covering;

(c) image of $F(X')$ in $P^5$ is contained in the Pl\"ucker quadric 
(apparently this is not true for the image of $F_m(X)$ in $P^5$); 

(d) There are no exceptional points and straight lines  on $F(X')$. 

[16] contains a conjecture that $X$ and $X'$ can be birationally isomorphic. 
Differences of properties of their Fano surfaces given above show that this 
is few likely. Nevertheless, there exists a Fano threefold $X''$ ([1]) whose 
properties are in some sense "intermediate" between properties of $X$ and 
$X'$. 

Namely, let $G_3$ be a $G(2,5) \cap P^3$. $X''$ is a two-sheeted covering of 
$G_3$ with ramification in $W=G_3 \cap \Omega$, where $\Omega$ is a quartic 
surface in $P(\lambda ^2 (V))$. There exists a variation in a smooth family 
of $X$ in $X''$. Let $K$ be a cone over $G_4$ with vertex $t$ and $X_4=K 
\cap \Omega$. Then a section of $X_4$ by a hyperplane which does not contain 
$t$ is isomorphic to a $X$ and a section of $X_4$ by a hyperplane which 
contains $t$ 
is isomorphic to a $X''$. From another side, some properties of $X''$ are 
analogous to the ones of $X'$. For example, $F_0(X'')$ (an analog of
$F_0(X')$ ) is isomorphic to the set of conics on $G_3$ which are bitangent
to $W$. 
\bigskip
\bigskip
{\bf Section 4. Tangent bundle theorem.}
\bigskip
Firstly we shall prove this theorem "locally". Recall that there exist maps 
$F \to F_0 \to G_Z(2,M^*)$ and $F \to P_Z(M^*)$, so the corresponding 
tautological sheaves $\tau _{2,m}$ and $\tau _{3,m}$ are defined on $F$. 
We denote $l = l_1 \cup l_2 \hookrightarrow F$. 
\medskip
{\bf Proposition 4.1.} Let $c \in F-l$. Then there exists a natural
isomorphism $\Omega_F(c) \to \tau _{2,m}^*(c)$.
\medskip
Proof of this proposition is a logically necessary step of the proof of the 
Tangent Bundle Theorem for $F$, because it implies the following 
\medskip
{\bf Corollary 4.2.} $F$ is smooth, $\dim F =2$, $T$ is a regular bundle. 
$\square$ 
\medskip
Steps of the proof of Proposition 4.1 correspond to ones of the proof of 
Theorem 4.14 below. The proof of (4.14) is more complicated only because
of necessity to consider exceptional straight lines $l_1$, $l_2$.
\medskip
{\bf Proof of 4.1.} We introduce the following notation. Let $\alpha: 
P_{P(V^*)}(M) \to P(V^*)$ is the projectivization of $M$ and $O_{\alpha}(1)$ 
is the tautological sheaf on $P_{P(V^*)}(M)$. There exists an inclusion of 
sheaves $M \to V_8 \otimes O$ on 
$P(V^*)$ and hence maps 

$$\matrix 
  P(M(V_4))   & \phantom{aaaaaaaaaaaaaa} X & \hookrightarrow & G_4  & 
\hookrightarrow & G \\
\downarrow    &    &                 & \downarrow &                 & 
\downarrow \\ 
P_{P(V^*)}(M) &\hookrightarrow P_{P(V^*)}(V_8 \otimes O) & 
\overset{\pr}\to{\to} & P(V_8) & \hookrightarrow & P(\lambda ^2(V)) 
\endmatrix $$
where vertical maps are inclusions and $P_{P(V^*)}(V_8 \otimes O) = P(V^*) 
\times P(V_8)$. 

We denote by $O_{\alpha}(1)$ the sheaf $O_{P(\lambda ^2(V))}(1)$ and its 
inverse images in all
obgects of this diagram, as well as in subvarieties of $P(M(V_4))$ (for a 
generic $V_4$). 

Let $c = (c, V_4) \in F-l$. We consider diagrams
$$\matrix c & \to & \pi(c) \\  \downarrow & & \downarrow \\ Q_G(V_4) & \to & 
P(M(V_4)) \endmatrix 
\eqno{(4.3)}$$
$$\matrix c & \to & X \\  \downarrow & & \downarrow \\ Q_G(V_4) & \to & G_4 
\endmatrix 
\eqno{(4.4)}$$ 
According deformation theory, $T_F(c)=H^0(c, N(c,X))$, if 
$$H^1(c,N(c,X))=0\eqno{(4.5)}$$
Later it will be shown that condition (4.5) is satisfied, hence we denote 
$ H^0(c, N(c,X))$ by $T_F(c)$. Analogous notations will be used in other 
cases where we shall use 
deformation theory. According Serre duality for $c$ we have: 
$$ H^0(c, N(c,X))= H^1(c, N(c,X)^* \otimes \omega _c)^*= 
H^1(c, N(c,X)^* \otimes \det N(c,X)^* \otimes \omega _c)^*$$ 
$$= H^1(c, N(c,X) \otimes \omega _X)^*$$
and analogously $H^1(c, N(c,X))= H^0(c, N(c,X) \otimes \omega _X)^*$ , i.e. 
$$\Omega_F(c)=H^1(c, N(c,X) \otimes \omega _X)^*\eqno{(4.6)}$$
if $H^0(c, N(c,X) \otimes \omega _X)=0$. 

Since $c \not\in l$ we have that (4.3) is a fibred product, and hence 
$$ N(c,Q_G(V_4))=N(\pi (c), P(M(V_4))) = O_{\alpha}(1)$$
Exact sequences of normal sheaves for both pairs of inclusions in (4.4) are
the following: 
$$ 0 \to N(c,X) \to N(c, G_4) \to N(X, G_4) \vert _c \to 0 \eqno{(4.7)}$$
$$ 0 \to N(c,Q_G(V_4)) \to N(c, G_4) \to N(Q_G(V_4), G_4) \vert _c \to 0 
\eqno{(4.8)}$$
Let us consider the composite map 
$$O_{\alpha}(1) = N(c,Q_G(V_4)) \to N(c, G_4) \to N(X, G_4) \vert _c = 
O_{\alpha}(2)\eqno{(4.9)}$$
Further, let us consider a diagram 
$$\matrix c & \hookrightarrow & c \cup c' & \to & X \\  & & \downarrow & & 
\downarrow \\ & & Q_G(V_4) & \to & G_4 \endmatrix $$
where $c' = i_f(c)$. We have 
$c \cup c' = X \underset{G_4}\to{\times} Q_G(V_4)$, so 
$\tilde N(c \cup c', Q_G(V_4)) \vert _c \to N(X,G_4) \vert _c $ is an 
isomorphism (here and below $\tilde N (Y_1, Y_2)$ means --- for an inclusion 
of varieties $Y_1 \subset Y_2$ --- a locally free sheaf on $Y_2$ such that 
$Y_1$ is the set of zeros of a section of it). 
\medskip
{\bf (4.10).} There is a map $ N(c,Q_{\Omega}(V_4)) \to 
\tilde N(c \cup c', Q_{G}(V_4)) \vert _c$ which is a multiplication by a 
section of a sheaf 
$O_c(c \cap c') = O_{\alpha}(1)$ which corresponds to the divisor  $c \cap c'$
on $c$.
\medskip
It is clear that $N(Q_G(V_4), G_4)= \tau_{2,V}^* \vert _{Q_G(V_4)}$. 
Multiplying (4.7) and (4.8) by $\omega _X=O_{\alpha}(-1)$ we get 
$$0 \to N(c,X) \otimes \omega _X \to N(c, G_4) \otimes O_{\alpha}(-1)   
\to   O_{\alpha}(1)   \to  0 
\eqno{(4.11)}$$
$$0 \to  O   \to N(c, G_4) \otimes O_{\alpha}(-1)   \to   \tau _{2,V} 
\vert _c   \to  0 
\eqno{(4.12)}$$
According (4.10), the composite map 
$$\phi _c: O \to N(c, G_4) \otimes O_{\alpha}(-1)   \to   O_{\alpha}(1)$$
corresponds to the divisor $c \cap c'$ on $c$. Let us consider the long exact 
cohomology sequences for (4.11) and (4.12). Since $c \not\in l$, we have 
$\tau _{2,V} \vert _c = 2 O_c(-1)$ and
$H^i(c, \tau _{2,V} \vert _c)=0$. This means that the long exact cohomology 
sequence for (4.11) is the following: 
$$0 \to H^0( N(c,X) \otimes \omega _X) \to H^0(O_c) 
\overset{H^0(\phi_c)}\to{\to} H^0(O_{\alpha}(1))    \to    \Omega_F(c)    
\to   0      \eqno{(4.13)}$$
We have $H^0(c,  O_{\alpha}(1))=  \tau _{3,M}^* (c)$  (see (4.43) below) and 
$H^0(\phi_c)$ is dual to an epimorhpism
$$\tau_{3,M}(c) \to \tau_{3,M}(c) / \tau_{2,M}(c)$$
This means that $H^0(\phi_c)$ is an inclusion, $H^0( N(c,X) \otimes 
\omega _X)=0$ and 
$\Omega _F(c)=\tau_{2,M}^*(c)$   $\square$
\medskip
{\bf Theorem 4.14.}  (The Tangent Bundle Theorem). There exists an exact 
sequence of sheaves on 
$F$: 
$$0 \to \tau_{2,M}^* \to r^*(\Omega _{F_m}) \to i_*(O_l \oplus O_l) \to 0$$
where $i$ is an inclusion $l \hookrightarrow F$. 
\medskip
{\bf Proof. } Let us consider a Fano family $D_F(c) \hookrightarrow 
F \times X$ defined as follows: 
$(c,x) \in D_F(c) \iff x \in c$, and let 
$$\matrix D_F(c) & \hookrightarrow & f \times X \\ & 
\overset{\beta}\to{\searrow} & \downarrow \\ & & F \endmatrix $$ be a diagram 
of natural projections on $F$, the vertical map is $\alpha$. According the 
deformation theory, $r_F^*(T_{F_c})=
\beta_*(N(D_F(c), F \times X))$, if $ \beta_{*1}(N(D_F(c), F \times X))=0$. 
Like in (4.6), the relative Serre duality for $\beta$ implies that 
$$r_F^*(\Omega_{F_c})=\beta_{*1}(N(D_F(c), F \times X) \otimes 
O_{\alpha}(-1))$$ if
$\beta_{*}(N(D_F(c), F \times X) \otimes O_{\alpha}(-1))=0$.

Let us define the following subvarieties of $P_F(M)$ (here $(f,x)$ is an 
element of 
$P_F(M) \subset F \times P(\lambda ^2(V))$ and $f=(c, V_4)$ ): 

$$D_F(G)=\{(f,x) \vert x \in Q_G(V_4)\}$$

$$D_F(P^2)= P_F(\tau_{3,M})=\{(f,x) \vert x \in \pi (c)  \}$$

$$D_l(P^2)= P_l(\tau_{3,M})=\{(f,x) \vert f \in l \; \hbox{ and } \; 
x \in \pi (c)  \}$$

Recall that for $f=(c,V_4) \in F$ we have $Q_G(V_4) \cap \pi (c)=c$ if $f 
\not\in l$ and 
$Q_G(V_4) \cap \pi (c)=\pi (c)$ if $f \in l$. Hence there is a diagram which 
is a fibred product: 
$$\matrix D_l(P^2)\cup D_F(c) & \to & D_F(P^2) \\ \downarrow & & \downarrow 
\\ D_F(G) & \to 
& P_F(M) \endmatrix \eqno(4.15)$$ 
Normal sheaves of inclusions of this diagram are the following. Since 
$N(l,F)=O_{\eta}(1)$, we have 
$$N(D_l(P^2), D_F(P^2))=O_{\eta}(1)\eqno{(4.16)}$$
and since $M/ \tau_{3,M}=O_{\pi}(1)$ (recall that the projectivization map
$P_Z(M^*) \to Z$ is denoted by $\pi$ ) we have 
$$N(D_F(P^2), P_F(M))=O(_{\alpha}1,_{\pi}1)\eqno{(4.17)}$$
Let us consider a diagram
$$\matrix D_F(c) & \to & F \times X \\ \downarrow & & \downarrow \\ D_F(G) & 
\to 
& F \times G_4 \endmatrix \eqno(4.18)$$
We have 
$$N(F \times X, F \times G_4)=O_{\alpha}(2)\eqno{(4.19)}$$
$$\tilde N(D_F(G), F \times G_4)=O(1) \otimes \tau_{2,V}^* \eqno{(4.20)}$$
where (4.20) is obtained by restriction to $F \times G_4$ of the corresponding 
equality for a general flag variety $D_{P(V_n^*)}(G) \subset P(V_n^*) \times 
G(k, V_n)$ defined as follows: 
$$ D_{P(V_n^*)}(G) =\{(V_{n-1},g) \vert V_{n-1} \in P(V_n^*), g \in G(k, V_n), 
g_{V_n} \subset V_{n-1} \}$$
We have $N(D_{P(V_n^*)}(G), P(V_n^*) \times G(k, V_n)) = O_{P(V_n^*)}(1) 
\otimes \tau_{k,V_n}^*$.
\medskip
{\bf Lemma 4.21.} $\Sing D_F(G)=S$, where $S \subset F \times Q_G(V_4)$ is 
defined as follows: 
$$S=\{(f,x) \vert f \in l, x \in \{\gamma_1(f), \gamma_2(f)\} \}$$
($f \in F$, \ \ \ $x \in Q_G(V_4)$). \ $S$ is a double curve on $D_F(G)$.
\medskip
{\bf Proof.} Explicit calculations in coordinates of the equation of the 
tangent space of points 
of $D_F(G)$, and of the equation of the tangent cone of points of $S$. 
We omit these calculations. $\square$
\medskip
Let us consider the blowing up of (4.15) along $S$ (recall that the blowing 
up along $S$ of a variety $Y \supset S$ is denoted by $\tilde Y \overset {r} 
\to {\to} Y$):
$$\matrix \widetilde{ D_l(P^2)} \cup \widetilde{ D_F(c)} & \to & 
\widetilde{ D_F(P^2)} \\ 
\downarrow & & \downarrow \\ \widetilde{ D_F(G)} & \to 
& \widetilde{ P_F(M)} \endmatrix \eqno(4.22)$$
(4.22) is also a fibred product. (4.16), (4.17) imply that normal sheaves for 
inclusions are the following: 
$$N(\widetilde{ D_l(P^2)}, \widetilde{ D_F(P^2)} )=O(_{\eta}1, _r 1)
\eqno{(4.23)}$$
$$N(\widetilde{ D_F(P^2)}, \widetilde{ P_F(M)} )=O(_{\alpha}1,_{\pi}1, _r 1)
\eqno{(4.24)}$$
where $O_r(1)$ is the tautological sheaf of the map $r$. (4.23), (4.24) imply 
$$N(\widetilde{ D_F(c)}, \widetilde{ D_F(G)}) = 
O(_{\alpha}1,_{\pi}1, _{\eta}-1) \eqno{(4.25)}$$
\medskip
{\bf Lemma 4.26.} Let us denote $N=N(\widetilde{ D_F(G)}, \widetilde{F \times
G_4})$. Then $N$ appears in the exact sequence (4.31) below. 
\medskip
{\bf Proof.} Let us consider a diagram of restictions of blowings up on $S$ 
$$\matrix \tilde S_1 & \hookrightarrow & \tilde S_2 & \hookrightarrow & 
\tilde S_3 & 
\hookrightarrow & \tilde S_4 \\ 
\downarrow & & \downarrow & & \downarrow & & \downarrow \\
S & = & S & = & S & = & S \endmatrix$$
where $ \tilde S_1$,   $ \tilde S_2$,   $ \tilde S_3$,   $ \tilde S_4$ are 
inverse images of $S$ in 
$\widetilde{ D_F(c)}$,  $\widetilde{ D_F(G)}$, $\widetilde{ P_F(M)}$, 
$\widetilde{ F \times G_4}$  respectively; projections (vertical arrows) are 
denoted by $r_1$, $r_2$, $r_3$, $r_4$ respectively. For $s \in S$ \ \ \ 
$r_1^{-1}(s) = P^1$,   $r_2^{-1}(s)$
 is a 2-dimensional quadric, $r_3^{-1}(s) = P^3$ which is the linear envelope 
of $r_2^{-1}(s) $, and $r_4^{-1}(s) =  P^4$. 

(4.20) and properties of blowings up imply that     
$$\tilde N(\tilde S_3 \cup \widetilde{ D_F(G)}, \widetilde{ F \times G_4} = 
O(1, _r1) \otimes \tau_{2,V}^* \eqno{(4.27)}$$
If we have a diagram of inclusions
$$ \matrix Y' =Y_1 \cap Y_2 & \overset {i} \to{\to} & Y_1 \\
\downarrow & & \downarrow \\ 
Y_2 & \to & Y \endmatrix \eqno{(4.28)}$$
where all objects are smooth varieties, $\codim _Y Y_i=k$, $\codim _{Y_i}Y'=1$ 
($i=1,2$), then there exists an exact sequence of sheaves on $Y_1$
$$0 \to N(Y_1, Y) \to \tilde N(Y_1 \cup Y_2, Y) \vert _{Y_1}  \to i_* (N(Y', 
Y_1) \otimes N(Y', Y_2)) \to 0 \eqno{(4.29)}$$
We let in (4.28): $Y=\widetilde{F \times G_4}$, $Y_1=\widetilde{D_F(G)}$, 
$Y_2= \tilde S_3$, $Y'= \tilde S_2$. In this case (4.29) is the 
following: 
$$ 0 \to N \to O(1, _r1) \otimes \tau_{2,V}^* \to i_* (N(\tilde S_2, \tilde 
S_3) \otimes N(\tilde S_2, \widetilde{D_F(G)})) \to 0 \eqno{(4.30)} $$
We have $N(\tilde S_2, \widetilde{D_F(G)})=O_r(-1)$, and since for $s \in S$ \ 
\ \ $r_2^{-1}(s)$ is a quadric in $r_3^{-1}(s)=P^3$, we have $ N(\tilde S_2, 
\tilde S_3)=O_r(2) \otimes r_2^*(L)$ for some sheaf $L$ on $S$. It is 
possible to show that $L=\tilde N(D_F(G), P_F(M)) \vert _S = 
O(-1, _{\alpha}2)$, but we do not need it. So, (4.30) can be written as 
follows: 
$$ 0 \to N \to O(1, _r1) \otimes \tau_{2,V}^* \to i_* (O_r(1) \otimes 
r_2^* (L)) \to 0 \eqno{(4.31)} $$
This is an exact sequence of sheaves on $\widetilde{D_F(G)}$ that we need. 
$\square$
\medskip
Now we need an exact sequence of normal sheaves for a pair of inclusions 
$\widetilde{D_F(c)}  \hookrightarrow \widetilde{D_F(G)} \hookrightarrow 
\widetilde{F \times G_4}$: 
$$ 0 \to O(_{\alpha}1,_{\pi}1, _{\eta}-1) \to r^* N(D_F(c), F \times G_4) 
\otimes O_r(1) \to N \vert _{\widetilde{D_F(c)}} \to 0 \eqno{(4.32})$$
Let us restrict (4.31) to $\widetilde{D_F(c)}$. Since $\tilde S_2$ and 
$\widetilde{D_F(c)}$ are divisors in $\widetilde{D_F(G)}$,  $ \tilde S_2 
\cap \widetilde{D_F(c)} = \tilde S_1$ and $\codim _{\widetilde{D_F(G)}} 
\tilde S_1=2$, we have: 
$\Tor _1 ^{\widetilde{D_F(g)}}(i_*O_{\tilde S_2}, i_* 
O_{\widetilde{D_F(c)}} )=0$, and the restriction of (4.31) is the following: 
$$ 0 \to N \vert _{\widetilde{D_F(c)}} \to O(1, _r1) \otimes \tau_{2,V}^* \to 
i_* O_r(1) \otimes L \to 0 \eqno{(4.33)}$$
where $i$ is an inclusion $\tilde S_1 \to \widetilde{D_F(c)}$.

Multiplying (4.32) and (4.33) by $O_r(-1)$ we get 
$$0 \to O(_{\alpha}1,_{\pi}1, _{\eta}-1, _r -1) \to r^* N(D_F(c), F \times 
G_4) \to N \otimes O_r(-1) \vert _{\widetilde{D_F(c)}} \to 0 \eqno{(4.34})$$
$$ 0 \to N \otimes O_r(-1) \vert _{\widetilde{D_F(c)}} \to O(1) \otimes 
\tau_{2,V}^* \to i_* L \to 0 \eqno{(4.35)}$$
Now we consider direct images of the above sequences to $D_F(c)$. To do it,
we consider an exact sequence of sheaves on $\widetilde{D_F(c)}$:
$$0 \to O \to O_r(-1) \to i_*O_r(-1) \vert _{\tilde S_1} \to 0 $$
We get that $r_*O_r(-1)=O$, $r_{*i}O_r(-1)=0$ for $i > 0$. Long exact 
$r_*$-direct 
images sequences of (4.34), (4.35) are the following (here we denote 
$N'= r_* (N \otimes O_r(-1) \vert _{\widetilde{D_F(c)}})$ ): 
$$0 \to O(_{\alpha}1,_{\pi}1, _{\eta}-1 ) \to N(D_F(c), F \times G_4) \to N' 
\to 0 \eqno{(4.36})$$
We get: $ r_{*1} (N \otimes O_r(-1) \vert _{\widetilde{D_F(c)}})=0$, and 
$$ 0 \to N' \to O(1) \otimes \tau_{2,V}^* \to i_* L \to 0 \eqno{(4.37)}$$

An exact sequence of normal sheaves for a pair of inclusions $ {D_F(c)}  
\hookrightarrow F \times X \hookrightarrow {F \times G_4}$ is the following: 
$$ 0 \to N(D_F(c), F \times X) \to N(D_F(c), F \times G_4) \to O_{\alpha}(2) 
\to 0 \eqno{(4.38})$$
Multiplying (4.36) --- (4.38) by $O_{\alpha}(-1)$ and taking the long exact 
$\beta_*$-direct image sequence we get (here we denote $B = N(D_F(c), F 
\times G_4) \otimes O_{\alpha}(-1)$ ): 
$$0 \to O(_{\pi}1, _{\eta}-1 ) \to \beta_*(B) \to \beta _*(N' \otimes 
O_{\alpha}(-1)) \to 0 \eqno{(4.39})$$
$$0 \to \beta_{*1}(B) \to \beta _{*1}(N' \otimes O_{\alpha}(-1)) \to 0 
\eqno{(4.40})$$
$$0 \to \beta _{*}(N' \otimes O_{\alpha}(-1)) \to
O(1) \otimes \beta _{*} \tau_{2,V} \to
\beta _{*} i_* (L \otimes O_{\alpha}(-1)) \to 
\beta _{*1}(N' \otimes O_{\alpha}(-1)) \to $$ $$
O(1) \otimes \beta _{*1} \tau_{2,V} \to 
\beta _{*1} i_* (L \otimes O_{\alpha}(-1)) \to 0 \eqno{(4.41)}$$
$$ 0 \to \beta_{*}(N(D_F(c), F \times X) \otimes O_{\alpha}(-1))  \to 
\beta_{*}(B) \to 
\beta_{*}O_{\alpha}(1) \to 
r_F^*(\Omega_{F_c}) \to $$ $$
\beta_{*1}(B) \to 
\beta_{*1}O_{\alpha}(1) \to 0 \eqno{(4.42})$$
It is easy to see that the restriction of $\tau_{2,V}$ on a conic line $c 
\subset G$ is equal to $2O_c(-1)$ (resp. $O_c \oplus O_c(-2)$) if $c$ is not 
contained in a $\alpha$-plane (resp. if $c$ is contained in a $\alpha$-plane). 
There exists the following property of the map $\beta: D_F(c) \to F$: 
\medskip
For $f \in F$ we have: $\beta^{-1}(f)$ is a conic line on an $\alpha$-plane 
$\iff f \in l_2$. 
\medskip
Hence, only for $f \in l_2$ cohomology of a sheaf $\tau_{2,V} 
\vert _{\beta^{-1}(f)}$ are non-zero. This means that $\beta _{*} 
\tau_{2,V}=0$ and $\beta _{*1} \tau_{2,V}$ is a sheaf whose support is $l_2$,
and it is invertible on $l_2$. (4.41) implies that $\beta _{*}(N' \otimes
O_{\alpha}(-1))=0$, and (4.39) implies that $\beta_*(B)$ is isomorphic to 
$ O(_{\pi}1, _{\eta}-1 )$.

Since we have a commutative diagram 
$$\matrix S & \overset{i}\to{\hookrightarrow} & D_F(c) \\
\downarrow & & \downarrow \\
l & \overset{i}\to{\hookrightarrow} & F \endmatrix $$
where vertical maps are $\beta_S$, $\beta$ respectively, we see that 
$\beta _{*k} i_* (L \otimes O_{\alpha}(-1)) = i_* \beta _{S*k} (L \otimes 
O_{\alpha}(-1))$. Since $\beta_S$ is a two-sheeted covering, we have 
$\beta _{*1} i_* (L \otimes O_{\alpha}(-1))
=0$ and $\beta _{*} i_* (L \otimes O_{\alpha}(-1))$ is a locally free sheaf 
of rank 2 on $l$. Let us consider a diagram
$$ \matrix D_F(c) & \overset{i_c}\to{\hookrightarrow} & D_F(P^2) & = & 
P_F(\tau_{3,M}) \\
& \overset{\beta}\to{\searrow} & \downarrow & & \\
& & F & & \endmatrix$$
where the vertical map is $\alpha$. Since for any $f \in F$ we have an 
equality of sheaves on 
$\alpha ^{-1}(f)$: $O_{\alpha ^{-1}(f)}(\beta^{-1}(f))=O_{\alpha}(2)$, we 
have globally on $D_F(P^2)$:
$O_{D_F(P^2)}(D_F(c))=\alpha^* L_F \otimes O_{\alpha}(2)$, where $L_F$ is
some invertible sheaf on $F$. There is an exact sequence of sheaves on 
$D_F(P^2)$: 
$$ 0 \to L_F^{-1} \otimes O_{\alpha}(-1) \to O_{\alpha}(1) \to 
i_{c*}(O_{\alpha}(1) \vert _{D_F(c)}) \to 0 $$
Considering its long exact $\alpha_*$-direct image sequence we see that 
$$\beta_{*}O_{\alpha}(1)=\tau_{3,M}^* \eqno{(4.43)}$$
and that $\beta_{*1}(O_{\alpha}(1))=0$. 

Now let us consider a map $i_B: \beta_{*}(B) \to \beta_{*}O_{\alpha}(1)$ 
which is isomorphic to 
$O(_{\pi}1, _{\eta}-1) \to \tau_{3,M}^*$. Its fibre at a point $c=(c, V_4) 
\in F - l$ is $i_B(c)=H^0(\phi_c)$ because of (4.13), i.e. $i_B(c)$ is an 
inclusion 
$(\tau_{3,M}(c)/\tau_{2,M}(c))^* \to \tau_{3,M}^*(c)$. Tautological exact 
sequences on $P_Z(M^*)$ and on $G_Z(2, M^*)$ and an equality of sheaves on 
$F$: $O_{\eta}(1) \otimes \det \tau_{2,M}^* = O$ imply an equality 
$(\tau_{3,M}/\tau_{2,M})^* = O(_{\pi}1, _{\eta}-1)$. So, there are 2 maps of 
sheaves 
$$(\tau_{3,M}/\tau_{2,M})^* = O(_{\pi}1, _{\eta}-1) \to \tau_{3,M}^*$$
namely $i_B$ and the natural inclusion 
$(\tau_{3,M}/\tau_{2,M})^*  \hookrightarrow \tau_{3,M}^*$ whose cokernel is 
$\tau_{2,M}^*$. For all $f \in F-l$ their fibres at $f$ coincide, hence the
maps coincide as well. So, the exact sequence (4.42) (resp. (4.41), taking 
into consideration (4.40) ) can be rewritten as follows: 
$$ 0 \to \tau_{2,M}^* \to r^*(\Omega_{F_c}) \to \beta_{*1}(B) \to 0 
\eqno{(4.44)}$$
respectively
$$0 \to \beta _{*} i_* (L \otimes O_{\alpha}(-1)) \to \beta _{*1}(B) \to 
O(1) \otimes \beta _{*1} \tau_{2,V} \to 0 \eqno{(4.45)}$$
where $\beta _{*} i_* (L \otimes O_{\alpha}(-1))$ is a locally free sheaf of 
rank 2 on $l$ and $ O(1) \otimes \beta _{*1} \tau_{2,V} $ is an invertible 
sheaf on $l_2$. 

Let us prove that $ O(1) \otimes \beta _{*1} \tau_{2,V} = \Omega _{l_2}$ and 
the composition map 
$$r^*(\Omega_{F_c})  \to \beta _{*1}(B) \to 
O(1) \otimes \beta _{*1} \tau_{2,V}$$
coincides with the natural epimorphism $r^*(\Omega_{F_c})  \to 
i_{2*}(\Omega_{l_2})$ associated to an inclusion $i_2: l_2 \to F_c$. 

We define respectively $D_{l_2}(c)$, $D_{l_2}(P^2)$, $D_{l_2}(G)= D_{l_2}(P^2) 
\cup l_2 \times U^*$ as restrictions on $l_2$ of $D_F(c)$, $D_F(P^2)$, 
$D_F(G)$ respectively. We find $\Omega_{l_2}$ by using deformation theory and 
a fact that $l_2$ is a variety of $\alpha$-planes on $G_4$. Let
$$\matrix D_{l_2}(P^2) & \hookrightarrow & l_2 \times G_4 \\
& \searrow & \downarrow \\
& & l_2 \endmatrix $$
be the corresponding family (both maps to $l_2$ are $\alpha$); then we have:
\medskip
$T_{l_2}= \alpha_{*}(N(D_{l_2}(P^2), l_2 \times G_4))$, if 
$$\alpha_{*1}(N(D_{l_2}(P^2), l_2 \times G_4))=0\eqno{(4.46)}$$
Serre duality implies that $\Omega_{l_2}=\alpha_{*2}(N(D_{l_2}(P^2), l_2 
\times G_4) \otimes \omega_{G_4})$. Now we apply (4.29) to the diagram 
$$ \matrix P_{l_2}(\tau_{2,M}) & \overset{i}\to{\to} & D_{l_2}(P^2) \\
\downarrow & & \downarrow \\
l_2 \times U^* & \to & l_2 \times G_4 \endmatrix $$
We get an exact sequence of sheaves on $D_{l_2}(P^2)$: 
$$ 0 \to N(D_{l_2}(P^2), l_2 \times G_4) \to O(1) \otimes \tau_{2,V}^*  
\to i_*(N_1 \otimes N_2) \to 0 $$
where $N_1 = N(P_{l_2}(\tau_{2,M}), D_{l_2}(P^2))$, 
$N_2 = N(P_{l_2}(\tau_{2,M}), l_2 \times U^*)$. We multiply it by 
$\omega_{G_4} = O_{\alpha}(-3)$ and consider the long exact 
$\alpha_{*}$-direct 
image sequence. Since $ \forall x \in l_2$ we have $\tau_{2,V}^* 
\vert _{\alpha ^{-1}(x)} = O \oplus O_{\alpha}(1)$, then 
$\alpha_{*i}(O(1,_{\alpha}-3) \otimes \tau_{2,V}^*  )$ is 0 for $i=0,1$ and 
one-dimensional for $i=2$. Restriction of both sheaves $N_1$ and $N_2$ to 
$\tau_{2,M}(x)$ (i.e. the fibre of the projection $ P_{l_2}(\tau_{2,M}) \to 
l_2$ at $x$ ) is $O_{\alpha}(1)$, so 
$N_1 \otimes N_2 \otimes O_{\alpha}(-3)$ at this fibre is equal to 
$ O_{\alpha}(-1)$, and its cohomology is 0. This implies that the condition 
(4.46) is satisfied and 
$$\Omega_{l_2} \to \alpha_{*2}(O(1,_{\alpha}-3) \otimes \tau_{2,V}^*  )
\eqno{(4.47)}$$
is an isomorphism. 
It is clear that a divisor $D_{l_2}(c) \overset{i}\to{\hookrightarrow}
D_{l_2}(P^2)$ corresponds to a sheaf $ O_{\alpha}( 2)$, so we consider an 
exact sequence on $ D_{l_2}(P^2)$: 
$$ 0 \to O_{\alpha}( -2) \to O \to i_*(O_{ D_{l_2}(c)}) \to 0 $$
We multiply it by $O(1) \otimes \tau_{2,V} = O(1,_{\alpha}-1) \otimes 
\tau_{2,V}^*$ and consider a long exact $\alpha_{*}$-direct image sequence: 
$$ 0 \to \alpha_{*1}( O(1) \otimes \tau_{2,V} \vert _{D_{l_2}(c)} ) \to  
\alpha_{*2}(O(1,_{\alpha}-3) \otimes \tau_{2,V}^* ) \to 0 \eqno{(4.48)}$$
(4.47) and (4.48) imply existence of an isomorphism that we need. It is clear 
that it is 
compatible with the map $r_F^*(\Omega_{F_c}) \to i_{2*}(\Omega_{l_2})$. The 
kernel of this map is $r^*(\Omega_{F_m})$, where $r: F \to F_m$, and the 
exact sequences(4.44), (4.45)
can be rewritten as follows: 
$$0 \to \tau_{2,M}^* \to r^*(\Omega_{F_m}) \to \beta _{*} i_* (L \otimes
O_{\alpha}(-1)) \to 0 \eqno{(4.49)}$$
It is clear that $ r^*(\Omega_{F_m}) \vert _l = 2O_l$, and hence the 
restriction of (4.49) on $l$ gives us an exact sequence 
$$ 2O_l \to \beta _{*} i_* (L \otimes O_{\alpha}(-1)) \vert _l \to 0 
\eqno{(4.50)}$$
Since $\beta _{*} i_* (L \otimes O_{\alpha}(-1)) \vert _l $ is a locally free 
sheaf of rank 2, (4.50) imply that it is isomorphic to $2O_l$ (equality 
$2O_l = \beta _{*} i_* (L \otimes O_{\alpha}(-1))$ can be proved also by 
means of explicit calculations).
 So, (4.49) can be rewritten as follows: 
$$0 \to \tau_{2,M}^* \to r^*(\Omega_{F_m}) \to 2O_l \to 0 \eqno{(4.51)}$$
Q.E.D. $\square$
\medskip
It is clear that (4.51) is isomorphic the the exact sequence associated to 
the divisor $l \subset F$ 
$$ 0 \to O_{\eta}(-1) \to O \to i_*(O_l) \to 0 $$
multiplied by $ r^*(\Omega_{F_m})$. So, $\tau_{2,M}^* = r^*(\Omega_{F_m}) 
\otimes O_{\eta}(-1)$. 
\medskip
{\bf Remark.} (4.51) can be used for a calculation of Chern classes of 
$\Omega_F$, because calculation of the class of $\tau_{2,M}^*$ in $K_0(F)$ 
is much easier than the one of $\Omega_F$.
\bigskip
\bigskip
{\bf Section 5. Special intersection of 3 quadrics in $P^6$.}
\bigskip
In this section $X$ will mean a special intersection of 3 quadrics in $P^6$ 
(see below for its definition). 

Firstly we continue to describe some properties of $G(2,5)$. For a variety 
$Y \subset P^n = P(V_{n+1})$ we denote by $TP_Y(t) \subset P^n$ (resp. 
$TV_Y(t) \subset V_{n+1}$) the projective (resp. vector) tangent space to 
$Y$ at a point $t \in Y$, i.e. $TP_
Y(t)  = P(TV_Y(t))$.

For any $n$ \ \ \ $G(2, V_n)$ is an intersection of quadrics in 
$P(\lambda^2V_n)$, and the set of these quadrics is isomorphic to $G(n-4, 
V_n)$. Namely, for $V_{n-4} \subset V_n$ the corresponding quadric is 
$\Pl (V_n/V_{n-4}) \in P(S^2 \lambda^2 V_n/V_{n-4})^*$. For $n=5$ we denote 
the corresponding inclusion by $q: P(V) \to P(S^2 \lambda^2 V^*)$; it is 
linear. All quadrics $q(P(V))$ are cones whose vertex has dimension 3. 

Let $t \in G$, $v \in P(V)$. Then $t \in \Sing q(v) \iff v \in t_V$. 
\medskip
{\bf 5.1.} Let $TV_G(t) = V_7 \subset \lambda^2(V)$. Then $V_7$ is a member 
of exact sequences 
$$0 \to \lambda^2 t_V \to V_7 \to t_V \otimes V/t_V \to 0 $$
$$0 \to V_7 \to \lambda^2 V \to \lambda^2 (V/t_V) \to 0 $$
$\forall v \in t_V$ we have $\Sing q(v) \subset TP_G(t) \subset q(v)$. The 
same inclusions are valid for projections from $t$. 
\medskip
{\bf 5.2.} Let $c \subset G$ be a conic line such that $\pi(c) \not\subset G$. 
Then there exists the only $V_4 \in P(V^*)$ such that $v \in P(V_4) \iff 
\pi(c) \subset q(v)$. It is easy to see that $c \subset G(2,V_4)$. 

Let $t \in G$; for any variety $Y$ we denote by $Y_t$ the projection of $Y$
from $t$. Since $\deg G = 5$, we have: $\deg G_t=4$, for $v \in t_V$ \ \ \ 
$(q(v))_t$ is a quadric containing $G_t$. This means that $G_t$ is a complete 
intersection of 2 quadrics. They generate a straight line $P^1$ which is 
isomorphic to $t_V$, all quadrics belonging to this line are cones whose 
vertex is 2-dimensional, and they contain $(TP_G(t))_t$. 

Let $t \in G_4$ and $t_V \cap U = \emptyset$; then $(G_4)_t$ is a complete 
intersection of two quadrics in $P^6$, they generate a straight line $P^1$ 
which is isomorphic to $t_V$, all quadrics belonging to this line are cones 
whose vertex is a point, they
 contain $P^3_T = (TP_G(t) \cap H_1 \cap H_2)_t$, and their vertices are 
contained in $P^3_T$ and form a normcubic $\bar c_3$ in it. The set of 
straight lines on $G_4$ passing through $t$ is canonically isomorphic to
$\bar c_3$; a projection from $t$ sends any such a line to the corresponding 
point of $\bar c_3$. 

Let us recall now properties of $X^3_8$ --- a generic intersection of 3 
quadrics in $P^6 = P(V_7)$ (see [21], [28], [47], [50]). Let $P^2$ be the 
linear envelope in $P(S^2V_7^*)$ of these 3 quadrics, and for $h \in P^2$ 
let $Q(h)$ be the corresponding quadric. The Hesse curve $H_7 = H_7(X^3_8)$ 
is defined as follows: $h \in H_7 \iff Q(h)$ is a cone. The set of planes on 
a non-singular 4-dimensional quadric $Q$ has 2 connected components, we denote
this set of components by $\Gamma (Q)$. Planes belonging to one component 
either coincide or their intersection is one point; intersection of planes 
belonging to different components is either empty or a straight line, and 
analogously 
for $P^3$ lying on a cone over such a quadric. We denote a two-sheeted 
covering $p = p_{H_7} : \hat H_7 \to H_7$ as follows: $p^{-1}(h) = 
\Gamma(Q(h))$.

For any plane curve $H$ of degree $n$ and its two-sheeted covering 
$p: \hat H \to H$ whose involution is denoted by $i_H: \hat H \to \hat H$, 
we define a set $\bar S = \bar S(p)$ of effective divisors on $\hat H$ as 
follows: 
\medskip
$d \in \bar S \iff O_H(p_*(d))=O_{P^2}(1) \vert _H \iff p_*(d)$ is a sum of
$n$ points $H \cap l$ where $l$ is a straight line in $P^2$. 
\medskip
$\bar S$ is non-connected, $\bar S = S \cup S'$, and if $d \in \bar S$, 
$d = \sum _{i=1}^n \hat h_i$ for $\hat h_i \in \hat H$, then $d$ and $d - 
\hat h_1 + i_H(\hat h_1)$ belong to different components of $\bar S$. 
Particularly, if $n$ is even then $S$ has an involution 
$i_S$: $i_S(\sum \hat h_i) = \sum i_H(\hat h_i)$. $S$ is irreducible ([50]). 

There exists an isomorphism $\phi_7: F(X^3_8) \to S = S(p_{H_7})$, where
$ F(X^3_8)$ is a set of conics in $X^3_8$. $\phi_7$ is defined as follows. 
For a conic $c \subset X^3_8$ we define firstly a straight line $l(c) 
\subset P^2$ as follows: $h \in l(c) 
\iff Q(h) \supset \pi(c)$. For $h \in l(c) \cap H_7$ \ \ \ $\pi(c)$ defines 
an element $\gamma(c,h) \in \Gamma (Q(h)) = p^{-1}(h) \subset \hat H_7$. We 
let 
$$\phi_7(c) = \sum_{ h \in l(c) \cap H_7} \gamma(c,h)$$
A fact that $\phi_7$ is an isomorphism follows easily from the following 
proposition ([47]): 
\medskip
{\bf 5.3.} Let $l \subset P(S^2 V_7^*)$ be a straight line in a space of 
quadrics in $P^6$, \ \ \ $h_i \in l$ 
$(i = 1, \dots , 7)$ are such that $Q(h_i)$ are cones. Then $\cap_{h \in l}
Q(h)$ contains $2^6$ planes. Let $\pi$ be such a plane. It defines a 
component in $\Gamma(Q(h_i))$ denoted by $\gamma(\pi, h_i)$. Let us consider 
$\gamma(\pi, h_i)$ as a function in $i$, it maps a set $\{1, \dots , 7\}$ to 
$\cup^7_{i=1}\Gamma(Q(h_i))$. Then these functions are different for 
different $\pi$, and values of 2 such functions are different at even numbers 
of elements of $\{1, \dots , 7\}$.
\medskip
The composition map $S \to \Div \hat H_7 \to \Pic \hat H_7$ (defined up to a 
shift) factors through the map $\Pr (\hat H_7/ H_7) \hookrightarrow \Pic \hat
H_7$. 
\medskip
{\bf Theorem 5.4.} ([50]) The corresponding map $\Alb (S) \to 
\Pr (\hat H_7/ H_7)$ is an isogeny. 
$\square$
\medskip
{\bf Theorem 5.5.} ([21]) There exists an isomorphism $\Pr (\hat H_7/ H_7) 
\to J^3(X^3_8)$. 
$\square$
\medskip
We denote the Abel - Jacobi map for conics by $\Phi$. We get a diagram where 
the left vertical map is $2\Phi$: 
$$\matrix \Alb F & \overset{\Alb(\phi_7)}\to{\to} & \Alb S \\
\downarrow & & \downarrow \\
J^3(X^3_8) & \gets & \Pr (\hat H_7/ H_7) \endmatrix $$
\medskip
{\bf Theorem 5.6.} ([28]) It is commutative. 
$\square$
\medskip
Let us consider now properties of a threefold $G_4 \cap \Omega$ which is 
smooth everythere except one point $t$ which is an ordinary double point. 
Recall that in this section it will be  
denoted by $X$. The set of conics on $X$ will be  
denoted by $F=F(X)$. It is well-known
that the projection of $X$ from $t$ is an intersection of 3 quadrics in 
$P^6$; this fact is used for a proof of non-rationality of a "generic" 
$G(2,5) \cap H_1 \cap H_2 \cap \Omega$ ([26]). 
\medskip
{\bf Lemma 5.7.} There exists a cone $\Omega ' \in P(S^2V_8^*)$ such that 
$t$ is its vertex and such that $X = G_4 \cap \Omega '$. 
\medskip
{\bf Proof.} We deprojectivize the map $q: P(V) \to P(S^2 \lambda^2V^*)$ and 
we 
restrict it to $V_8$. We get a map $q: V \to S^2V_8$. Let $T: V \to V_8^*$ be 
a map defined as follows: $T(v)$ is an equation of the polar line of $q(v)$ 
with respect to $t$. Then $\im T \subset (V_8/TV_{G_4}(t))^*$ and $\Ker T$ is 
a set of points $v \in V$ such that $t \in \Sing q(v)$, i.e. $\Ker t = t_V$, 
hence $T$ is an epimorphism on $(V_8/TV_{G_4}(t))^*$. Since $t$ is a singular
point in $G_4 \cap \Omega$, we have $TV_{\Omega}(t) \in (V_8/TV_{G_4}(t))^*$ 
and $\exists v \in V$ such that $T(v)= TV_{\Omega}(t) $. This means that the 
equation of the tangent hyperplane at $t$ of the quadric $\Omega - q(v)$ 
is 0, i.e. $\Omega ' = \Omega - q(v)$ is the desired cone. $\square$
\medskip
Further we shall suppose that $\Omega$ is such a cone. We see that $X_t$ is 
an intersection of 3 quadrics: $\Omega _t$ and $(q(t_V))_t$. As earlier, let 
$P^2$ be the plane that they generate, and we denote by $t_P$ the straight 
line that corresponds to quadrics $(q(t_V))_t$, or, in formal notations,
$h \in t_P \iff Q(h) \in (q(t_V) \vert _{H_1 \cap H_2})_t$. 

$t_P$ is a component of $H_7(X_t)$, another component is a curve of degree 6 
denoted by $H_6=H_6(X)$. It is clear that for a generic $X$ it is 
non-singular, later we shall consider only these $X$. Let $t_P \cap 
H_6=\{R_1, \dots , R_6\}$. $Q(R_i)$ are ordinary cones, like all other 
quadrics corresponding to points of $t_P$ and $H_6$. Hence there are a 
two-sheeted non-ramified covering $p = p_{H_6}: \hat H_6 \to H_6$ and a 
corresponding surface $S = S(p_{H_6})$ together with an involution $i_S$. 
Vertices of
 $Q(R_i)$ are projections of 6 straight lines on $X$ passing through $t$, and 
they 
are the only singular points on $X_t$. These vertices are intersection points 
of $\bar c_3$ and $\Omega _t$ as well.
\medskip
{\bf Proposition 5.8.} There exists a birational isomorphism $\phi_6: F(X) 
\to S(p_{H_6})$. 
\medskip
{\bf Proof.}  Let $c=(c, V_4) \in F$. If $t \not\in \pi(c)$ then the 
projection $c_t \subset X_t$ is a conic. We denote by $\chi$ a rational map 
$(c, V_4) \mapsto c_t$ from $F(X)$ to $F(X_t)$. $c_t$ defines a straight line 
$l(c) = l(c_t) \subset P^2$. We denote 
$ l(c) \cap t_P = h_0(c)$, $l(c) \cap H_6 = \{h_1(c), \dots, h_6(c)\}$. If
$c \notin l_1 \cup l_2$, then $l(c) \ne t_P$ and (5.2) implies that 
$Q(h_0(c)) = q(t_V \cap V_4)$. 
\medskip
{\bf Lemma 5.9.} $\Sing Q(h_0(c)) = P(M(V_4))_t \cap P^3_T$.
\medskip
{\bf Proof.} This statement is formulated on $G(2,V)$ as follows. For some 
given 
$t \in G$, $V_4 \subset V$ we set $s=t_V \cap V_4$. Then the linear envelope 
of $t$ and $P(\lambda^2(V_4) \cap TP_G(t)$ coincide with $\Sing q(s)$. 
Really, let $s_0$, $t_1$ (resp. $s_0$, $v_1$, $v_2$, $v_3$) be a basis of 
$t_V$ (resp. $V_4$). Then (considering vector spaces) we have 

a basis of $t \subset \lambda^2V$ is $s_0 \land t_1$

a basis of $TV_G(t)$ is $s_0 \land t_1$, $s_0 \land v_i$, $t_1 \land v_i$

a basis of $\lambda^2 V_4$ is $s_0 \land v_i$, $v_i \land v_j$

and a basis of $\Sing q(s)$ is $s_0 \land t_1$, $s_0 \land v_i$ (see (5.1); 
$i,j = 1,2,3$). This implies the proposition for $G$. Intersecting it with 
$H_1$, $H_2$ and projecting from $t$ we get the desired. $\square$
\medskip
$\forall h \in t_P$ the quadric $Q(h)$ contains $P^3_T$, we denote the 
corresponding element of $\Gamma(Q(h))$ by $\beta(h)$. 
\medskip
{\bf Lemma 5.10.} $\beta(h_0(c)) \ne \gamma(c_t, h_0(c))$. 
\medskip
{\bf Proof.} Since $\pi(c) \subset P(M(V_4))$, we have: $\pi(c_t)$ and 
$(P^3_T)_{\Sing Q(h_0(c))}$ have the empty intersection after the projection 
from ${\Sing Q(h_0(c))}$. $\square$
\medskip
We define $\phi_6(c)=\sum_{i=1}^6 \gamma(c_t, h_i(c))$. (5.10) and (5.3) imply 
that at a generic point $\phi_6$ is an inclusion. Since $S$ is irreducible and 
$\dim F = \dim S$ we have that $\phi_6$ is a birational isomorphism. $\square$
\medskip
{\bf Lemma 5.11.} $\phi_6$ commutes with involutions $i_F$, $i_S$. 
\medskip
{\bf Proof.} Let $(c, V_4)$ and $(c',V_4) \in F$ be 2 involutory conics. 
Then $\pi(c) \cap \pi(c')$ is a straight line, hence (if $t \notin P(M(V_4))$, 
i.e. if $V_4 \not\supset t_V$) \ \ \ $\pi(c_t) \cap \pi(c'_t)$ is also a 
straight line. For $h \in l(c)$ \ \ 
\ $Q(h)$ contains $c'_t$ and $\pi(c_t) \cap \pi(c'_t)$, i.e. a curve of degree 
3 on $\pi(c'_t)$ and hence the whole $\pi(c'_t)$, i.e. $l(c)=l(c')$, 
$h_i(c)=h_i(c')$. The intersection of projections of $\pi(c_t) $ and 
$ \pi(c'_t)$ from $\Sing Q(h_i(c))$ ($i = 
1, \dots , 6$) is a straight line, and if 
$$\Sing Q(h_i(c)) \notin P(M(V_4))_t$$
then $\gamma (c_t, h_i(c)) \ne \gamma (c'_t, h_i(c'))$. But else 
$$ P(M(V_4))_t \subset Q(h_i(c))  \Rightarrow \Sing Q(h_0(c)) \in Q(h_i(c)).$$
If $h_i(c) \notin \{R_1, \dots, R_6\}$ then $h_i(c) \notin t_P$ and $\Sing 
Q(h_0(c))$ belongs to 3 basis quadrics of $P^2$ and hence to $X_t$, i.e. 
$h_0(c) \in \{R_1, \dots, R_6\}$. This is not true at an open set, hence at 
this set $c$ and $c'$ are mapped to involutory elements of $S$, hence the 
last is true at all $F$. $\square$
\medskip
The above description of $\phi_6$ does not indicate images of points 
$(c, V_4) \in F$ for which $V_4 \supset t_V$, as well as images of straight 
lines $l_1, l_2 \subset F$. The set of spaces $V_4$ which contain $t_V$ is 
isomorphic to $P^2_T= P((V/t_V)^*) \hookrightarrow P(V^*)$. For this $V_4$ 
\ \ \ $Q_{\Omega}(V_4)$ is a cone in $P(M(V_4))$ whose vertex is $t$, and 
$Q_G(V_4))$ is a generic quadric containing $t$. Let $(c, V_4) \in F$, $V_4 
\in P^2_T$ and 
$$Q= T( \bar \phi(c)) = \pi(c) \cup \pi(c') \in < Q_{\Omega}
(V_4), Q_G(V_4) >$$ 
is a pair of planes. If $Q \notin Q_{\Omega}(V_4) $, then 
$ Q_{\Omega}(V_4) \cap Q_G(V_4) = c \cup c' = Q \cap Q_{\Omega}(V_4) $. Let 
$t \in \pi(c)$, then $ Q_{\Omega}(V_4) \cap \pi(c)$ is a pair of straight 
lines passing 
through $t$, and $ Q_{\Omega}(V_4) \cap \pi(c')$ is a generic conic. There are
6 straight lines on $X$ passing through $t$, i.e. there are only 
$\binom 62$ = 15 such 
spaces $V_4$ while there are 30 conics. It is possible to prove (although we 
do not need it) that $\phi_6$ maps these conics (as well as lines $l_1, l_2 
\subset F$) to 32 points $S \cap ((p_{H_6})^{-1} (\sum _i R_i))$. 

Let now $Q= Q_{\Omega}(V_4) $. We denote the set of these $V_4 \in P^2_T$ by 
$K_0$, i.e. $K_0 = l_{\Omega } \cap P^2_T$. We multiply an inclusion of 
vector spaces $(t) \to V_8$ by $O_{P^2_T}$; it is clear that the obtained 
map $t \otimes O_{P^2_T} \to V_8
 \otimes O_{P^2_T} $ factors through the inclusion $M \vert _{P^2_T} \to V_8 
\otimes O_{P^2_T} $, and we denote $\Coker (t \otimes O_{P^2_T} \to M 
\vert _{P^2_T} )$ by $M_t$. Since $\Omega$ is a cone whose vertex is $t$, 
we have: a bundle of quadrics $Q_{\Omega}$ restricted to $P^2_T \subset
P(V^*)$ induces a bundle of quadrics $ Q_{\Omega_t} $ on $(P^2_T, M_t)$ 
whose fibre at $V_4 \in P^2_T$ is a projection $Q_{\Omega}(V_4)$ from its 
vertex. It is clear that $K_0=D_2(Q_{\Omega_t})$, hence $K_0$ is a curve on 
$P^2_T$, and according the Porteus formula ([7]) we have $\deg K_0 = 2 c_1 
(M_t^*) = 6$. 

Since $K_0 \subset F_0(X)$, there exists a two-sheeted non-ramified covering 
$p_K: K \to K_0$ whose involution is denoted by $i_K : K \to K$, and $p_K = 
p_F \vert _K$, $i_K = i_F \vert _K$. Involutory conics corresponding to
points of $K$ have 2 points of
 intersection, one of which is $t$. For a generic cone $\Omega$ whose vertex 
is $t$, we have: $K_0$ is non-singular and $K$ is connected. 

Let $d: \tilde F \to F$ be the desingularization map and 
$\tilde K = d^{-1} (K)$. There exists an involution $i_{\tilde F}$ on $\tilde 
F$ such that $i_F \circ d = d \circ i_{\tilde F}$. There exist rational maps 
$\tilde \phi _6: \tilde F \to S$ and $\tilde 
\chi : \tilde F \to F(X_t)$. 
\medskip
{\bf Proposition 5.12.} $\tilde K$ is non-connected, $\tilde K = K_1 \cup K_2$ 
and $d_i = d \vert _{K_i} : K_i \to K$ are isomorphisms ($i=1,2$). For 
$f \in \tilde K$ we have: $d(f)$ is a conic on $X$ passing through $t$,
$\tilde \chi (f)$ is a pair of straight lines on $X_t$, one of them is 
contained in a quadric $Q_2 = P^3_T \cap \Omega _t$. 
\medskip
{\bf Proof.} Since $S$ is smooth and $\tilde K$ is not a union of straight 
lines, we 
have that $\tilde \phi_6 (\tilde K)$ is a curve on $S$. Let us find which 
points $s \in S$ do not belong to $\phi _6 (F - K)$. Let $s = \sum _{i=1}^6 
\hat h_i$, \ \ \ $\hat h_i \in \hat
 H_6$, \ \ \ $p_{h_6} (\hat h_i) = h_i \in H_6$, and $l$ a straight line in 
$P^2$ such that 
$\{h_i\} = l \cap H_6$. We can suppose that $l \ne t_P$ (this is true for 
almost all points $s \in S$). We let $h_0 = l \cap t_P$ and $\hat h_0 
\in \Gamma (Q(h _0))$, \ \ \ $\hat h_0 \notin \beta (h_0)$. According 
(5.3), there exists a plane $\pi = \pi (s)$ on $\cap_{h \in l}  Q(h)$ that 
corresponds to the set $\{ \hat h_i \}$ ($i = 0, \dots , 6$). Really, their 
images under the projection from $\Sing Q(h_0)$ must intersect namely 
by this way. Let $\pi \cap P^3_T = \emptyset$; then $c = \pi \cap 
X_t$ is a conic. Let us consider the 2-dimensional cone containing $c$ and 
having a vertex $t \in P(V_8) \supset X$. Since $X$ is an intersection of 
quadrics, any straight line either intersects $X$ at most at 2 points, or is 
contained in $X$. Since $c \subset X_t$, any ruling line of the above cone 
intersects $X$ at one point (the other intersection point is $t$). This point 
does not coincide with $t$ because $c \cap P^3_T=\emptyset$. So, the set of 
these intersection points is a conic, and $s \in S$ is its image under the map
$\phi _6$. (If we chose a plane $\pi$ on $\cap _{h \in l} Q(h)$ that 
corresponds to a set $\beta (h_0), \hat h_1, \dots , \hat h_6$, then $\pi 
\cap P^3_T$ is one point, and instead of a conic we should get a normcubic on 
$X$ passing throught $t$). 

So, if $s = \tilde \phi_6(c)$ for $c \in \tilde K$, then $\pi(s) \cap P^3_T = 
P^1$. It is clear that 
$\pi (s) = \pi (\tilde \chi (c))$ and $\forall h \in P^2 - l(c)$ we have 
$\tilde \chi(c) = \pi(s) \cap Q(h)$. Choosing $h \in t_P$ we get that
$\tilde \chi(c) \supset \pi(s) \cap P^3_T = P^1$. This means that $\tilde 
\chi (c)$ is a pair of straight lines, one of which lies on a quadric $Q_2 = 
\Omega _t \cap P^3_T$, and other intersects it. 

Inversely, let $l_X$ be any straight line on $X_t$ that intersects $Q_2$. 
There are 2 
ruling lines $u_1$, $u_2$ of $Q_2$ that pass through a point $l_X \cap Q_2$, 
hence we get 2 conics $c_i(l_X)= l_X \cup u_i$ ($i=1,2$). We can associate 
them:

(a) 2 different straight lines $l(c_i(l_X)) \subset P^2$;

(b) points $h_{ij} = h_j(c_i(l_X)) \in H_6$;

(c) elements $\gamma (c_i(l_X), h_{ij}) \in \hat H_6$ 

such that the element 
$$s_i(l_X) = \sum _{j=1}^6 \gamma (c_i(l_X), h_{ij})\in S$$ 
and $\pi (s_i(l_X)) \cap P^3_T = u_i$. Let us find 
$\tilde \phi _6^{-1}(s_i(l_X))$. A cone over $c_i(l_X)$ having a vertex $t$ 
in $P(V_8)$ is a plane $<l_X, t>$, and each straight line
 in this plane passing through $t$ meets $X$ exactly at one point 
(except $t$). Further, one straight line in this plane (it is the inverse 
image of $l_X \cap Q_2$ under the projection from 
$t$) is tangential to $X$ at $t$. So, $<l_X, t> \cap X$ is a conic $c_X(l_X)
\in K$, and $\tilde \phi _6^{-1}(s_i(l_X)) \in 
d^{-1}(c_X(l_X))$. This means that $\forall c \in K$ \ \ \ $d^{-1}(c)$ is 2 
points, i.e. $\tilde K \to K$ is a non-ramified two-sheeted covering, and for 
$c \in \tilde K$ we have $\tilde \chi(c) = d(c)_t \cup u(c)$, where $u(c)$ is 
a ruling line of $Q_2$. 
Since there are 2 systems of rulings on $Q_2$, we have: $\tilde K = K_1 \cup 
K_2$ is non-connected, namely $c \in K_i$ iff $u(c)$ is a ruling line of type 
$i$. It is clear that restrictions $d_i = d \vert _{K_i} : K_i \to K$ are 
isomorphisms ($i=1,2$). $\square$ 
\medskip
{\bf Remark. }  Let $\tilde r: \tilde F \to \tilde F_m$ be a map of blowing 
down of $l_1 \cup l_2 \subset \tilde F$ to points. There exists an isomorphism 
$i: \tilde F_m \to S$ such that $\tilde \phi_6=i \circ \tilde r$.
\medskip
Let us denote the isomorphism $d_2^{-1} \circ d_1: K_1 \to K_2$ by $\alpha$, 
and let $i_{K_1} =  d_1^{-1} \circ i_K \circ d_1$ be an involution of $K_1$. 
\medskip
{\bf Lemma 5.13.} For $c \in K_1$ we have $i_{\tilde F}(c)=
\alpha (i_{K_1}(c))$. 
\medskip
{\bf Proof.} $i_{\tilde F}(c)$ is equal to $i_{K_1}(c)$ or 
$\alpha (i_{K_1}(c))$. We have: the intersection of $\tilde \chi (c)$ and 
$\tilde \chi (i_{K_1}(c))$ is 1 point, while
$\forall f \in \tilde F$ the intersection of $\tilde \chi (f)$ and 
$\tilde \chi (i_{\tilde F}(f))$ is 2 points. $\square$ 
\medskip
We denote by $\tilde X$ the desingularization of $X$. Any point of $\tilde F$ 
defines a conic on $\tilde X$ by an obvious way (for $c \in \tilde F$ the 
corresponding conic is  
the inverse image of $\tilde \chi(c)$ under the desingularization map 
$\tilde X \to X$). So, there 
exists the Abel - Jacobi map $\Phi: \tilde F \to J^3(\tilde X)$. 

The next 3 propositions are analogs for $X_t$ of the theorems (5.4) - (5.6) 
for a generic $X^3_8$, and their proofs follow the proofs of these theorems. 
\medskip
{\bf Proposition 5.14.} The map $\Alb S(p_{H_6}) \to \Pr \hat H_6 / H_6$ is an
isogeny. 
\medskip
{\bf Proof.} [50], Proposition 1.3.3 implies that this map is an epimorphism, 
hence it is sufficient to prove that $h^{1,0}(S) \le 9$. The proof of ([50], 
Proposition 3.4.4) implies that: 
\medskip
(1) There exist plane curves $C_i$ of degree 6 ($i=0, \dots , 6$), their 
non-ramified two-sheeted coverings $p_i: \hat C_i \to C_i$ and points 
$P_1, \dots , P_6 \in P^2$ in a general position such that:

(a) $p_0 = p_{H_6}$; 

(b) $\Sing C_i = \{P_1, \dots , P_i\}$ and these points are ordinary double 
points on $C_i$; 
\medskip
(2) $\forall i = 0, \dots , 5$ there exists a commutative diagram 

$$\matrix \hat D_i & \to & D_i \\ & \searrow & \downarrow \\ & & T_i 
\endmatrix $$ 
where $T_i$ is an open subset in $P^1$, $D_i$ is a smooth surface, the fibre 
of the covering $\hat D_i \to D_i$ at a point $t \in T_i$ is a non-ramified 
two-sheeted covering of smooth curves $(\hat D_i)_t \to (D_i)_t$ for all $t$
 except $t=0$, where $(D_i)_t
$ has one double point and $(\hat D_i)_t$ has 2 double points which are its 
inverse images. Further, for $t=t_0$ the above fibre is isomorphic to the 
normalization of the covering $p_i$ denoted by $\tilde p_i: \tilde 
{\hat C}_i \to \tilde C_i$, and for $t
=0$ it is isomorphic to the blowing up of the covering $p_{i+1}$ at  $P_1 
\cup \dots \cup P_i$ for $C_{i+1}$ and their $p_{i+1}$-inverse images for 
$\hat C_{i+1}$. 
\medskip
(3) There exist maps $D_i \to T_i \times P^2$ over $T_i$ whose fibres at 
points $t_0$ and $0 \in T_i$ are ordinary maps $\tilde C_i \twoheadrightarrow 
C_i \hookrightarrow P^2$ and $(\tilde C_{i+1})_{ P_1 \cup \dots \cup P_i} 
\twoheadrightarrow C_{i+1}\hookrightarrow P^2$.
\medskip
Some properties of $C_i$ and coverings $D_i \to T_i$ of purely technical 
nature are omitted, their complete list is given in [50], 2.2.1. 

These conditions and ([50], 2.2.7) imply that 
$$ \forall i \; \; h^{1,0}(S(\tilde p_i)) \le h^{1,0}(S(\tilde p_{i+1}))+1$$
If  points $P_1, \dots , P_6$ do not belong to a conic then the linear 
system of the sheaf $O_{P^2}(1) \vert _{\tilde C_6}$ is non-special. ([50], 
1.3.8) implies that $ h^{1,0}(S(\tilde p_6)) = g(\tilde C_6)-1 =3$ and hence
$h^{1,0}(S(p_{H_6}))\le 9$. $\square$
\medskip
{\bf Proposition 5.15.} There exists an isomorphism $I: \Pr \hat H_6/H_6 \to 
J^3(\tilde X)$. 
\medskip
{\bf Proof.} Let $\hat h \in \hat H_6$ and $p(\hat h)= h$. Let us consider the 
cone $Q(h)$. The set of spaces $P^3$ of type $\hat h$ on $Q(h)$ is isomorphic 
to $P^2$, so for all these $P^3$ curves $P^3 \cap X$ are birationally 
equivalent. We denote such a curve by $C_4(\hat h)$, they can be lifted on
$\tilde X$ and these lifts define a cylinder map $\Phi: H_1(\hat H_6, \n Z) 
\to H_3(\tilde X, \n Z)$ and an Abel - Jacobi map $\hat I: J(\hat H_6) \to 
J^3( \tilde X)$. 

Let us choose a generic straight line $l \subset X$; $\forall \hat h \in 
\hat H_6$ there exists the only $P^3$ of type $\hat h$ such that $P^3 
\supset l$. We have 
$$P^3 \cap X = C_4(\hat h) = l \cup C_3(\hat h, l)$$
where $C_3(\hat h, l)$ is a normcubic in $P^3$ meeting $l$ at 2 points. A 
system of curves $C_3$ defines the same cylinder map $\Phi: H_1(\hat H_6, 
\n Z) \to H_3(\tilde X, \n Z)$. Let us apply a formula ([20], p. 44) to 
$\Phi$: 
$$\forall \gamma_1, \gamma_2 \in H_1(\hat H_6, \n Z) \;\; <\Phi(\gamma_1),
\Phi(\gamma_2)> = k <\gamma_1, \gamma_2> + <i(\gamma_1), \gamma_2>$$
where $k \in \n Z$ and $i \in \End H_1(\hat H_6, \n Z)$ corresponds to a 
correspondence $C$ on $\hat H_6 
\times \hat H_6$ defined as follows: $(\hat h_1, \hat h_2) \in C \iff 
C_3(\hat h_1, l) \cap C_3(\hat h_2, l) \ne \emptyset$. Further, [20], p. 96 
implies that $i = H_1(i_H)$. 

From another side, we have $g(H_6)=10$, $h^{2,1}(\tilde X)=9$ (because 
$h^{2,1}(X)=10$ and because $\tilde X$ is a blowing up of a degenerated fibre
of a Lefschetz pencil). This implies that $\dim \Ker \Phi \ge h_1(\hat H_6) 
- h_3 (\tilde X) = 20$, i.e. for a 20-dimensional space of cycles $\gamma$ 
we have $i(\gamma) + k \gamma=0$. This can be only for $k=-1$, hence $\Ker 
\Phi = J(H_6)$ and $\Pr \hat H_6/H_6 \to J^3(\tilde X)$ is an isomorphism. 
$\square$
\medskip
{\bf Proposition 5.16.} The diagram 
$$\matrix \Alb (\tilde F) & \overset{\Alb (\tilde \phi_6)}\to{\to} & 
\Alb (S) \\
\downarrow & & \downarrow \\
J^3(\tilde X) & \overset{I}\to{\gets} & \Pr (\hat H_6/H_6) \endmatrix $$
where the left vertical map is $2 \Phi$, is commutative. 
\medskip
{\bf Proof.} Let $c \in \tilde F$, $l(c) \cap t_P=h_0$, $l(c) \cap H_6= 
\{h_1, \dots, h_6\}$ and $\gamma(c_t, h_i)= \hat h_i$. Then $\phi_6(c)=
\sum_{i=1}^6 \hat h_i$ and it is sufficient to prove that codimension 2 
cycles $2c - \sum_{i=1}^6 C_4(\hat h_i)$
 on $\tilde X$ are linearly equivalent for all $c \in \tilde F$. Let $X_4 = 
\cap_{h \in l(c)} Q(h)$ be an intersection of 2 quadrics and $P^3_i$ the 
only $P^3$ of type $\hat h_i$ containing $\pi(c_t)$. Then $X_4 \cap P^3_i$ 
is a quadric containing $\pi(c_t)$, i.e. $X_4 \cap P^3_i$ is a pair of planes 
$\pi(c_t) \cup L_i$. According ([47], Section 3, Lemma 6) there is a linear
equivalence ($\equiv$) on $X_4$ of 2 cycles: 
$$\sum_{i=0}^6 L_i + 5 \pi(c_t) \equiv 3(P^4 \cap X_4)$$
Intersecting this equivalence with $X$ we get: 
$$\sum_{i=0}^6 L_i \cap X + 6c \equiv 2c - (c + L_0 \cap X) + const $$
But $c + L_0 \cap X = const$, because $c \cup (L_0 \cap X) = P^3_0 \cap X$ 
and $P^3_0$ varies along $t_P= P^1$. Equality $C_4(\hat h_i)= (L_i \cap X) 
\cup c)$ implies that $\sum_{i=1}^6 C_4(\hat h_i) -2c \equiv const$. $\square$
\bigskip
\bigskip
{\bf Section 6. Abel - Jacobi map of $F(X)$ is an isogeny.}
\bigskip
Here $X$ will mean again a Fano threefold of genus 6. Let us prove that 
$h^{1,0}(F(X))= \dim J^3(X)=10$. We consider $X_4' = G \cap H_1 \cap \Omega$
and a Lefschetz pencil $f': X_4' \to P^1$, we restrict it to a neighbourhood 
of zero $\Delta \subset P^1$, 
so the only singular fibre ${f'}^{-1}(t)$, \ \ \  $t \in \Delta$, is a 
threefold $X_0={f'}^{-1}(0)$ having the only double point. We denote 
$F_3'=\cup_{t \in \Delta} F({f'}^{-1}(t))$, and $f'_F: F_3' \to \Delta$ is 
the projection. 
\medskip
{\bf Lemma 6.1.} For a generic $f': X_4' \to P^1$ \ \ \ $F_3'$ is non-singular.
\medskip
{\bf Proof.} Non-singularity of $F_3'$ is not obvious only at points of 
$K=F(X_0) \subset F_3'$. For a given $X_0$ it is possible to choose $X_4' 
\subset P^8 = P(H_1)$ such that $X_0 = X_4' \cap P^7$, where $P^7 \subset 
TP_{X_4'}(t)$ and such that all points of $K$ are non-singular on $F(X_4')$ 
--- the variety of conics on $X_4'$. A generic $P^6 \subset P^7$ defines a
Lefschetz pencil $f': X_4' \to P^1$; we denote by $U(P^6)$ a set of conics 
$c \subset P^8$ such that $\dim (\pi(c) \cap P^6) \ge 1$. We have $F_3'= 
U(P^6) \cap F(X_4')$. Let $c \in K$. Let us look which conditions must be 
imposed on $P^6 \subset P^7$ in order to achieve non-singularity of $F_3' =
F_3'(P^6)$ at points of $c$. For any variety $Y \ni c$ we denote $TV_Y(c)$ by 
$\bar Y$, and we denote by $U_7$ (resp. $U_8$) varieties of conics in $P^7$ 
(resp. $P^8$). We get a diagram of tangent spaces at $c$: 
$$ \matrix \overline{F(X_0)} & & \hookrightarrow & & \overline{F(X_4')} \\
 \downarrow & & & & \downarrow \\
\bar U_7 & \hookrightarrow & \overline{U(P^6)} & \hookrightarrow & \bar U_8 
\endmatrix $$
(vertical maps are inclusions). Here we have: $\overline{F(X_0)} = 
\overline{F(X_4')} \cap \bar U_7$, $\codim _{\bar U_8}\bar U_7 = 3$, $\dim 
\overline{F(X_4')} = 5$, $\dim \overline{F(X_0)} = 3$, because $c$ is a 
smooth point on $F(X_0)$. Condition of smoothness of $F_3'(P^6)$ at $c$ is 
equivalent to a condition that $\dim \overline{U(P^6)} \cap \overline{F(X_4')} 
= 3$, i.e. that $\overline{U(P^6)} \not\subset <\overline{F(X_4')}, \bar 
U_7>$. But $\overline{U(P^6)} \in P(\bar U_8/ \bar U_7)=P^2$, and a map $P^6 
\mapsto P(\overline{U(P^6)}/ \bar U_7) \in P^2$ is a linear projection from 
$P^{7*}$ to $P^2$ whose center $P^4 \subset P^{7*}$ is dual to $\pi(c) 
\subset P^7$. The inverse image of $<\overline{F(X_4')}, \bar U_7>$ under 
this projection is a hyperplane in $P^{7*}$ containing $P^4$, i.e. the set
of $P^6 \subset P^7$ which contain a fixed point in $\pi(c)$. It is easy to 
see that this point is $t$. So, for a non-singularity of $F_3'$ it is 
sufficient to choose a generic $P^6$ which does not contain 
$t$. $\square$
\medskip
It is clear that for any $X$ we can choose $f': X_4' \to P^1$ such that:

1. $X={f'}^{-1}(t_0')$ for $t_0' \in \Delta$;

2. The special fibre satisfies properties of general position of Section 5 
(smoothness of $H_6$, $K_0$);

3. $F_3'$ is non-singular. 
\medskip
{\bf 6.2.} To transform $f'$, $f'_F$ to a desired degeneration we use a 
standard 
method (see, for example, [30]). Let us consider blowings up of $X_4'$ and of 
$F_3'$ along singular points of the fibre over 0, we denote the corresponding 
maps by $\tilde f': \tilde X_4' \to \Delta$, $\tilde f'_F : \tilde F_3' \to
\Delta$. Their fibres at 0 are $\tilde X_0 \cup P^3$, $\widetilde{F(X_0)} 
\cup R$ respectively, where $P^3$, $R$ are exceptional divisors of blowing up 
of $X_4'$ and $F_3'$, and $R$ is a fibration over $K=K(F(X_0))$ whose fibres 
are isomorphic to $P^1$. It is clear that $\tilde X_0 \cap P^3=Q$ is a 
non-singular quadric in $P^3$ and $\widetilde{F(X_0)} \cap R = K_1 \cup K_2$. 
Components of the fibre at 0 of $P^3$ and $R$ are double, that's why we 
consider a map $\phi_2: \Delta \to \Delta$: $t \mapsto t^2$ and a diagram 
which is a fibred product: 
$$\matrix X_4 & \to & \tilde X_4' \\
\downarrow & & \downarrow \\
\Delta & \overset{\phi_2}\to{\to} & \Delta \endmatrix $$
where the right vertical map is $\tilde f'$. We denote the left vertical map 
by $f$. The fibre of $f$ at 
$t_0 = (t_0')^{\frac12}$ (resp. at $t=0$) is $X$ (resp. $X(0)= \tilde X_0
\cup \bar P_3$), where $\bar P_3$ is a two-sheeted covering of $P^3$ ramified 
at $Q$, so $\tilde X_0 \cap \bar P_3=Q$, and both components have multiplicity 
1. Analogously we define a map $f_F: F_3 \to \Delta$ having a singular fibre 
$F(0) = F \cup \bar R$, where $F=\widetilde{F(X_0)}$, and $\bar R \to R$ is a 
two-sheeted covering ramified at $K_1 \cup K_2
$, so $F \cap \bar R = K_1 \cup K_2$, and composition maps $K_i \hookrightarrow
 \bar R \to R \to K$ coincide with $d_i$ ($i=1,2$). 

Clemens - Schmidt exact sequence ([8], [40]) for $f_F$ is the following:

$$0 \to H^1(F(0)) \overset{\nu}\to{\to} H^1(F(X)) \overset{N}\to{\to} 
H^1(F(X)) \eqno{(6.3)}$$
There exists a mixed Hodge structure on $H^1(F(0))$, $H^1(F(X))$, and $\nu$ 
and $N$ are strict morphisms with respect to these structures of degrees (0,0) 
and (-1, -1) respectively. We shall need only weight filtrations; their 
definition in our case is the following ([8]): for $H^1(F(X))$ we have
$$ 0 \subset W_0 \subset W_1 \subset W_2= H^1(F(X))$$ 
where $W_0 = \im N$, $W_1= \Ker N$. For $H^0(F(0))$ we have
$$0 \subset W_0 \subset W_1=H^1(F(0))$$
where $W_0=H^1(\Pi(F(0)))$ and $\Pi(F(0))$ is a CW-complex whose vertices 
correspond to components of $F(0)$ and edges correspond to irreducible 
components of intersections of these components. We have $$W_1/W_0 =
\Ker H^1(F) \oplus H^1 ( \bar R) \to H^1(F \cap \bar R)$$ for a natural 
restriction map. In our case $\Pi(F(0)) = S^1$ (2 vertices, 2 edges), i.e. 
$\dim W_0(H^1(F(0)))=1$. Exactness of (6.3) implies that $\dim W_2/W_1 = 
\dim W_0(H^1(F(X))) = \dim W_0(H^1(F(0))) =1$ and $W_1/W_0$ for $H^1(F(X)
)$ and $H^1(F(0))$ coincide. So, the equality $h^1(F(X))=20$ is a corollary of 
$$\dim \Ker H^1(F) \oplus H^1 ( \bar R) \to H^1(F \cap \bar R)=18\eqno{(6.4)}$$
So, we have to prove (6.4).
\medskip
{\bf Lemma 6.5.} Let us identify $K_1$ and $K_2$ via $\alpha$. Then (6.4) 
follows from the following statement: the image of a restriction map 
$$H^1(F) \to H^1(K_1) \oplus H^1(K_2) \eqno{(6.6)}$$
is contained in the diagonal.
\medskip
{\bf Proof.} Since $\bar R \to K$ is a fibration with a fibre $P^1$, we have: 
$H^1(K) \to H^1( \bar R)$ is an isomorphism. This means that the restriction 
map $H^1(\bar R) \to H^1(K_1) \oplus H^1(K_2)$ is an inclusion, and its image 
is the diagonal. According (5.14), 
$h^1(F)=18$, this implies the lemma. $\square$
\medskip
{\bf Lemma 6.7.} (6.6) is equivalent to the following condition
$$ \forall h \in H^1(F) \ \ \ h \vert _K = i^*_K(i^*_F(h) \vert _K)
\eqno{(6.8)}$$
(here and below $K=K_1$). 
\medskip
{\bf Proof.} Follows from commutativity of the diagram 
$$\matrix K_1 & \overset{i_K}\to{\to} & K_1 & \overset{\alpha}\to{\to} & 
K_2 \\
\downarrow & & & & \downarrow \\
F & & \overset{i_F}\to{\to} & & F \endmatrix $$
(vertical maps are natural inclusions). $\square$
\medskip
There exists a composite map defined up to a shift 
$$f: F \overset{\tilde \phi_6}\to{\to} S \to \Pr \hat H^6/H^6$$
\medskip
{\bf Lemma 6.9.} $f \circ i_F = -f$ (up to a shift). 
\medskip
{\bf Proof.} For $c \in F$ let $f(c)=\sum_{i=1}^6 \hat h_i$ for $\hat h^i 
\in \hat H_6$; then $f(i_F(c))=\sum_{i=1}^6 i_H( \hat h_i)$. But 
$O_{\hat H_6}(\sum_{i=1}^6 (\hat h_i + i_H (\hat h_i)))= O_{P^2}(1) 
\vert _{\hat H_6}$, i.e. $ \sum \hat h_i \equiv 
- \sum i_H (\hat h_i) + \const$. $\square$ 
\medskip
This means that $H^1(i_F)$ is a multiplication by --1. This is equivalent to
the following: $H^1(F_0) =0$, where $F_0=\widetilde{F_0(X_0)}$.

Let $E_i \hookrightarrow H^1(K)$ ($i = \pm 1$) be eigenspaces of $H^1(i_K)$ 
with eigenvalues $i$, then $E_1 = H^1(K)$, $E_{-1}=H^1( \Pr K/K_0)$, and (6.8) 
is equivalent to the following: 
\medskip
The image of the restriction map $H^1(F) \to H^1(K)$ is contained in $E_{-1}$, 
or, which is the same, the composition map $\Alb K_0 \to \Alb K \to \Alb F$ is 
0.
\medskip
Taking into consideration (5.15) and (5.16), it is sufficient to prove that 
the map $\Alb K_0 \to J^3(\tilde X_0)$ is 0. Let $k_{10}, k_{20} \in K_0$ and 
$k_1$, $k_2$ be their representatives in $K$. Let us find a cycle $D$ on 
$\tilde X_0$ which is the image
 of $k_{10}- k_{20}$: 
$$D=\tilde \chi(k_1) + \tilde \chi(i_K(k_1)) - \tilde \chi(k_2) - \tilde 
\chi(i_K(k_2)).$$
We have $\forall f \in F$ \ \ $ \Phi (\tilde \chi(f)) + 
\Phi (\tilde \chi(i_F(f))) 
= const$ and $i_F(i_K(k_1)) = \alpha (k_1)$, hence 
$$D=\tilde \chi(k_1) - 
\tilde \chi(\alpha(k_1)) - \tilde \chi(k_2) + \tilde \chi(\alpha(k_2))= 
U(k_1) - U(\alpha(k_1)) + U(\alpha(k_2)) - U(k_2)$$ i.e. $D \equiv 0$. 
\medskip
So, we have proved that $h^{1,0}(F(X))=10$.
\medskip
{\bf Remark.} The above considerations imply that there exists a map 
$\Pr K/K_0 \to J^3(\tilde X_0) \to \Pr \hat H_6/H_6$. For a generic covering 
$p_{H_6}$ existence of this map implies that $p_K = p_{H_6}$ which implies 
in its turn that $p_K (X_0) = p_{H_6}(X_0)$ is true for any $X_0$. 
\medskip
Calculation of $h^{1,0}(F_0(X))$, as well as a proof of connectedness of 
$F(X)$ and of $F_0(X)$ is made by the same method, but easier. Analogously to 
a fibration $f_F : F_3 \to \Delta$ we construct a fibration $f_{0F}: F_{03}
\to \Delta$ having a generic
 fibre $f_{0F}^{-1}(t_0)=F_0(X)$ and a special fibre $f_{0F}^{-1}(0)=F_0(0)= 
F_0 \cup \bar R_0$, where $F_0 = \widetilde{F_0(X_0)}$ and $\bar R_0$ is a 
two-sheeted covering of $R_0$ ($R_0$ is a fibration over $K_0$ whose fibre is 
$P^1$) ramified at $K \subset R_0$, and $F_0 \cap \bar R_0 = K$. There exists 
a weight filtration on $H^1(F_0(0))$: 
$$0 \subset W_0 \subset W_1 = H^1(F_0(0))$$
where $W_0 = H^1(\Pi(F_0(0)))$, $W_1/W_0 = \Ker H^1(F_0) \oplus H^1(R_0) 
\to H^1(K)$. Since $\Pi(F_0(0))$ is a segment, then $W_0=0$. As we have seen 
(Lemma 6.9), $H^1(F_0)=0$. Since the cohomology map 
$H^1(K_0) \to H^1(\bar R_0)$ 
of the fibration $\bar R_0 \to K_0$ is an isomorphism, we have: 
$$\Ker H^1(\bar R_0) \to H^1(K) = \Ker H^1(K_0) \to H^1(K)=0$$ So, 
$H^1(F_0(0))=0$. It follows immediately from the Clemens - Schmidt exact
sequence of the fibration $f_{0F}$ that $H^1(F_0(X))=0$ as well. 

To calculate $H^0(F(X))$ we consider the even Clemens - Schmidt exact sequence 
of the fibration $f_{F}$:
$$0 \to H^0(F(0)) \to H^0(F(X)) \overset{N}\to{\to} H^0(F(X)) $$ 
Connectedness of $F$ implies that $h^0(F(0))=1$ and that the map $N$ is 0. 
This implies 
connectedness of $F(X)$ and hence of $F_0(X)$ as well. 
\medskip
{\bf Proposition 6.10.} The Abel - Jacobi map $\Phi: \Alb F \to J^3(X)$ is an
epimorphism (here $F = F(X)$). 
\medskip
{\bf Proof.} Let $F_1$ be the family of straight lines on $X$, $D \subset F_1 
\times F_1$ the incidence divisor, i.e. $(t_1, t_2) \in D \iff t_1 \ne t_2$ 
and $t_1 \cap t_2 \ne \emptyset$, \ \ \ $p_j: D \to F_1$ a projection to the 
$j$-th component $(j=1,2)$, $N = \deg p_j$, i.e. the number of straight lines 
meeting 
a given straight line, $i: D \to D$ the natural involution (permutation of 
$t_1$, $t_2$) and $p_D: D \to D_0$ the 
factorization of $D$ by $i$. Then $D_0 \overset{d}\to{\hookrightarrow} F$ is 
a set of reducible conics. Let us consider cylinder maps in complex homology 
\medskip
(a) for conics: $\Phi_F:H_1(F, \n C) \to H_3(X, \n C)$;

(b) for straight lines: $\Phi_{F_1}:H_1(F_1, \n C) \to H_3(X, \n C).$
\medskip
It is clear that for $\gamma \in H^1(D)$ we have 
$$\Phi_{F} \circ d_* \circ p_{D*}(\gamma) = \Phi_{F_1} (p_{1*}(\gamma)) + 
\Phi_{F_1} (p_{2*}(\gamma)) \overset{\opr}\to{=} B(\gamma)$$
If $\gamma = p_1^*(\gamma')$, where $\gamma'=H_1(F_1)$, then 
$$B(\gamma) = \Phi_{F_1} (p_{1*}(p_1^*(\gamma'))) + \Phi_{F_1} 
(p_{2*}(p_1^*(\gamma'))) = N \cdot 
\Phi_{F_1} (\gamma') + \Phi_{F_1} (p_{1*}(i^*(p_1^*(\gamma'))))$$
According ([20], Section 4.3, Lemma 6) we have
$$\Phi_{F_1} (p_{1*}(i^*(p_1^*(\gamma'))))=(1-n) \Phi_{F_1} (\gamma') $$
where $n$ is a number of straight lines passing through a generic point of 
$X_4'$, so 
$B(\gamma) = (N+1-n) \Phi_{F_1} (\gamma')$. For $X$ we have $N=11$, $n=6$. 
According to a lemma of Clemens ([20], Section 4.3) $\Phi_{F_1}$ is an 
epimorphism, hence $\gamma' \mapsto
 B(\gamma) = \Phi_F \circ d_* \circ p_D(\gamma)$ is an epimorhism as well. 
This implies that $\Phi_F$ is also an epimorphism, which is equivalent to 
(6.10). $\square$
\medskip
Since $\Alb F$ and $J^3(X)$ have the same dimension, we get 
\medskip
{\bf Theorem 6.11.} The Abel - Jacobi map $\Phi: \Alb F \to J^3(X)$ is an 
isogeny. $\square$
\medskip
Now we give a sketch of a proof of the following theorem: 
\medskip
{\bf Theorem 6.12.} The class of $\Phi(F)$ in $A(J^3(X))$ is 
$\frac{2 \Theta ^8}{8!}$ where $\Theta$ is the class of a Poincar\'e divisor.
\medskip
{\bf Sketch of proof.} For a degeneration $f: X_4 \to \Delta$ (see 6.2) we 
denote $f^{-1}(t)=X_t$ and $F(X_t)=F_t$ for $t \ne 0$. There exists a complex 
analytic fibration $j: J_{11} \to \Delta$ whose fibre at $t \in \Delta^*=
\Delta - \{0\}$ is $j^{-1}(t)
=J^3(X_t)$ ([52], Section 2). Its fibre $j^{-1}(0)=J(0)$ at 0 is a fibration 
over $J^3(\tilde X_0)$ whose fibre is $\n C^*$ ([52], 4.56). There exists a 
compactification of $J(0)$. It is a fibration $\mu: \tilde J_0 \to 
J^3(\tilde X_0)$ whose fibre is $P^1$. It has 2 sections $\alpha_0$, 
$\alpha_{\infty}: J^3(\tilde X_0) \to \tilde J_0$ which map any point 
$t \in J^3(\tilde X_0)$ to points $(0)$, $(\infty) \in P^1 = \mu^{-1}(t)$ 
such that $J(0)= \tilde J_0 - \im \alpha_0 - \im \alpha_{\infty}$. There 
exists also a compactification of $J_{11}$ which is a fibration $j_c: J_{11,c} 
\to \Delta$ whose non-singular fibre is the same as the one of $j$ and whose 
singular fibre is $j_c^{-1}(0)=J_0$. Further, there is a map $\nu: 
\tilde J_0 \to J_0$ which is the normalization of $J_0$. $\nu$ identifies 
(with a shift) $\im \alpha_0$ and $\im \alpha_{\infty}$, outside these
subsets $\nu$ is an inclusion. 

Let us consider now the Abel - Jacobi map. We fix a section $s_0 : \Delta 
\to F_3$ such that $s_0(0) \in F(X_0)-K$, and for any section $s : \Delta \to 
F_3$ such that $s(0) \in F(X_0)-K$ we consider an element $\Phi_s(t)=
\Phi(s(t) - s_0(t)) \in J^3(X_t)$.
 They form a section $\Delta ^* \to J_{11}$ which can be expanded to a 
section $\Phi_s: \Delta \to J_{11}$. We have: $\Phi_s(0)$ depends only on 
$s(0)$ but not on $s$ and $\mu(\Phi_s(0)) = \Phi_{\tilde X_0}(s(0) - s_0(0))$
([52], 4.58) where $\Phi_{\tilde
 X_0}$ is the Abel - Jacobi map of $\tilde X_0$. We denote by 
$\cup _s \Phi_s(t) = F_{\Phi}(t) \subset j^{-1}(t)$ the image of $F_t$ under 
the Abel - Jacobi map shifted by such a way that $s_0(t)$ goes to 0. For 
$t=0$ we have an inclusion $F(X_0) - K \hookrightarrow J(0)$. There exists a 
commutative diagram where the above inclusion is the left vertical map: 
$$\matrix F(X_0)-K & \hookrightarrow & \widetilde{F(X_0)} & 
\twoheadrightarrow & F(X_0) & = & F_{\Phi}(0) \\
\downarrow &                        & \downarrow          &            
                 &\downarrow & & \\
J(0) & \hookrightarrow & \tilde J_0 & \overset{\nu}\to{\twoheadrightarrow} & 
J_0 & & \endmatrix $$
(vertical maps are inclusions) and we get a family of surfaces $(F_{\Phi}
\to \Delta) \subset (j_c: J_{11,c} \to \Delta)$ whose fibre at $t \in \Delta$ 
is $F_{\Phi}(t) = F(X_t)$. Since 
$$\mu(\widetilde{F(X_0)}) = 
\Phi_{\tilde X_0}(\widetilde{F(X_0)}) \subset J^3(\tilde X_0), \ \ \ 
J^3(\tilde X_0)= \Pr(\hat H_6/H_6)\hbox{ and }\widetilde{F(X_0)} =S(p_{H_6})$$
then according ([28]) the homology class of $\Phi_{\tilde X_0}
(\widetilde{F(X_0)})$ in $H_4(J^3(\tilde X_0), \n Z)$ is 
$$\frac{2\Theta^7}{7!} = 2 \sum _{i < j} \gamma_i \times \delta_i \times 
\gamma_j \times \delta_j \eqno{(6.13)}$$
where $\Theta$ is the class of a Poincar\'e divisor and $(\gamma_i, \delta_i)$
is a simplectic basis of $H_1(J^3(\tilde X_0), \n Z)$ with respect to the 
principal polarization associated with $\Theta \subset J^3(\tilde X_0)$. 

There exists a map of topological spaces $(0,0): J^3(X_t) \to J_0$ ([25]) 
such that $(0,0) (F_{\Phi}(t))$ is homologically equivalent to $F_{\Phi}(0) 
\hookrightarrow J_0$. By analogy with proofs of ([25], Lemma 2.6) and ([31], 
Theorem 8.27) it is possible
 to show --- using (6.13) --- that the class of $ F_{\Phi}(t) \subset 
J^3(X_t)$ in $H_4(J^3(X_t), \n Z)$ is 
$$2 \sum _{i < j} \gamma_i \times \delta_i \times \gamma_j \times \delta_j = 
\frac{2\Theta^8}{8!}
\; \; \square$$
\bigskip
\bigskip
{\bf Section 7. Geometric interpretation of the tangent bundle theorem and 
recovering of $X$ by its Fano surface.}\nopagebreak
\bigskip
We denote $V_{10}=H^0(F,\Omega_F)^*$. Here we impose the following condition 
on $F$: 
\medskip
{\bf 7.1.} The natural map $V_{10}^*= H^0(F_m,\Omega_{F_m}) \to 
\Omega_{F_m}(c_{\Omega})$ is an epimorphism.
\medskip
Apparently this condition is true for any $X$ or at least for a generic $X$. 
The theorem on recovering of $X$ by its Fano surface $F(X)$ will be proved 
only for $X$ that satisfy (7.1). 

Let us define a map $B_8: F \to G(2, V_8)$ as follows: $B_8(f) = 
<\gamma_1 (f), \gamma_2 (f)>$. We denote by $B_{10}: F_m \to G(2, V_{10})$ a 
map associated to $\Omega_{F_m}$. (7.1) means that $B_{10}$ is regular at 
$c_{\Omega}$. It is clear that $B_8$ factors through $p_F : F \to F_0$. Since 
$H^0(F_0, \Omega_{F_0})=0$, we have: $B_{10}$ factors through $p_{F_m}: F_m 
\to F_{0m}$, and the corresponding map $F_{0m} \to G(2, V_{10})$ is 
associated to the sheaf  $\Omega_{F_{0m}}\otimes \sigma_m$.
\medskip
{\bf Theorem 7.2.} (Geometric interpretation of the tangent bundle theorem). 
The map $B_{10} \circ r: F \to G(2, V_{10})$ is regular. There exists the only 
epimorphism $p_{10,8}: V_{10} \to V_8$ such that
$$B_8=G(2, p_{10,8}) \circ B_{10} \circ r \eqno{(7.3)}$$
If (7.1) is satisfied then $ p_{10,8}$ is the projection from the straight 
line $B_{10}(c_{\Omega})$. 
\medskip
{\bf Proof.}  We denote $H=H^0(F, \tau_{2,M}^*)^*$ and let $B_H: F \to 
G(2,H)$ be the map defined by $\tau_{2,M}^*$. Let us consider the long exact
cohomology sequence for (4.51): 
$$0 \to H^* \overset{i_1}\to{\to} V_{10}^* \to H^0(2O_{l_1 \cup l_2})$$
It is clear that (7.1) is equivalent to a condition $H=V_8$. Since the linear 
envelope of lines $<\gamma_1 (f), \gamma_2 (f)> = l_8(f) \subset P(V_8)$ is 
the whole $P(V_8)$, there exists an inclusion $i_2: V_8^* = H^0(G(2, V_8), 
\tau_{2,V_8}^* ) \hookrightarrow H^* =H^0(F, \tau_{2,M}^*)$, and we have 
$$B_8=G(2, i_2^*) \circ B_H \eqno{(7.4)}$$
This implies that $\forall f \in F$ \ \ \ $B_H$ is regular at $f$. For 
$f \in F-l$ the restriction of (4.51) to $f$ is 
$$0 \to \tau_{2,M}^* \vert _f  \to r^* \Omega_{F_m} \vert _f \to 0$$
and there exists a commutative diagram 
$$ \matrix & & H^* & \overset{i_1}\to{\hookrightarrow} & V_{10}^* & & \\
& & \downarrow & & \downarrow & & \\
H^0(\tau_{2,M}^* \vert _f) & = & \tau_{2,M}^* (f) & \to & H^0(r^* \Omega_{F_m} 
\vert _f) & = & r^* \Omega_{F_m} (f) \endmatrix \eqno{(7.5)}$$
Since $H^* \to \tau_{2,M}^* (f)$ is an epimorphism, we have that $V_{10}^* 
\to r^* \Omega_{F_m} (f)$ is also an epimorphism, i.e. $B_{10}$ is regular at 
$F - l$. If so, there exists a continuation of $B_{10}$ to $l$, i.e. $B_{10}$ 
is regular on $F$. Commutativity of the diagram dual to (7.5) shows that 
$$B_H(f)=G(2,i_1^*) \circ B_{10}(f) \eqno{(7.6)}$$
We denote now $p_{10,8}=(i_1 \circ i_2)^*$. (7.4) and (7.6) imply that on 
$F-l$ we have
$B_8=G(2, p_{10,8}) \circ B_{10} \circ r$, hence this is true on the whole 
$F$. 

If (7.1) is true then for $f \in l$ we have $B_{10} \circ r(f) = B_{10} 
(c_{\Omega})$ which does not depend on $f$. But $B_8(f)$ for $f \in l$ are 
different straight lines having no commun points. This observation shows that 
(7.3) can be satisfied only if $p_{10,8
}$ is the projection from the straight line $B_{10}(c_{\Omega})$. $\square$
\medskip
{\bf Theorem 7.7.} Here we suppose that (7.1) is satisfied for $F_c$. Then 
$X$ such that $F_c  = F_c(X)$ is defined uniquely up to isomorphism. 
\medskip
{\bf Proof.} Let us consider the map $F_c \to G(2, V_{10})$ associated to
$\Omega_{F_c}$. It factors through $F_{0m}$ and hence it permits us to 
recover uniquely the point $c_{\Omega} \in F_c$, the straight line 
$l_2 \subset F_c$, the involution $i_F$ on $F$ and the straight line 
$B_{10}( c_{\Omega}) \subset P(V_{10})$. Since (7.1) implies that $p_{10,8}$ 
is the projection from $B_{10}(c_{\Omega})$, we can recover as well the map 
$B_8: F_0 \to G(2, V_8)$. Since the straight line $l \subset F_0$ can be 
recovered by $F_c$ and the straight lines $l_0(c)
$ for $c \in l$ are contained in $U^*$ and are tangent to the conic $c^*_u$, 
we have that $U^*$, $c^*_u$ and the isomorphism $l \to c^*_u$ are recovered
by $F_c$. 
\medskip
{\bf Proposition 7.8.} $c_{\Omega} \subset U^*$ can be recovered uniquely by 
$F_c$.
\medskip
{\bf Proof.} We prove 

(a) If $c \in F_0 - l$ and $l_8(c) \cap U^* \ne \emptyset$ then $l_8(c) 
\cap U^* \in c_{\Omega}$. 

(b) The quantity of these points $c \in F_0 - l$ is finite and $> 4$. 

So, having $F_c$ it is possible to construct more than 4 points on 
$ c_{\Omega} $, i.e. to recover $ c_{\Omega} $ itself. 
\medskip
(a) can be proved easily. For $c \in F_0 - l$ we have: $X \cap l_8(c)= G_4 
\cap l_8(c)$. Since $U^* \subset G_4$, we have $l_8(c) \cap U^* \in G_4$, and 
hence $l_8(c) \cap U^* \in X$. But $ c_{\Omega} = X \cap U^*$. 
\medskip
Let us define now a curve $S$ as a set of pairs $(t,c)$ where $t \in 
c_{\Omega}$, $c \in c_u \cap t_V$. We define a straight line $j(t,c)$ as a 
line lying in $\alpha$-plane $\alpha (c)$ (see Section 3) and which is the 
tangent of
$\alpha(c) \cap \Omega $ at $t$. Lemma 3.20 (whose notations are used here) 
implies that either $j(t,c) \subset X$ or $\{$ if for $f_0 \in F_0 - l$  \ \ 
\ $\tilde \phi(f_0) \in \widetilde{\pi(j(t,c))}$ $\}$ \ \ $t \in l_8(f_0)$. 
We define a surface
$$R=\bigcup_{(t,c) \in S} \widetilde{\pi(j(t,c))} \subset P(V^*)$$
It is clear that $R \supset l$. Let $V_4 \in (R - l) \cap (\phi(F_0 - l))$, 
\ \ \ $V_4 = \phi(f_0)$, \ \ \ $V_4 \in \widetilde{\pi(j(t,c))}$. Then either 

(i) $l_8(f_0) \cap c_{\Omega} = \{t\}$ or

(ii) There exists a straight line on $X$ passing through $t$. 

Let us find for how many $V_4 \in (R - l) \cap (\phi(F_0 - l))$ we have the
case (ii). Since the union of all straight lines on $X$ is an intersection of 
$X$ and a hypersurface of degree 10, there are 20 points on $ c_{\Omega}$ 
which belong to a straight line on $X$. Let $t$ and $j(t)$ are one of these 
points and 
lines. We define $c=\cap _{t' \in j(t)} t'_V$; it is clear that $c \in c_u 
\cap t_V$ and $j(t) = j(t,c)$. Since $j(t)$ meets 11 straight lines on $X$, 
we have: a straight line $\widetilde{\pi(j(t,c))} \subset R$ meets 
$\phi(F_0)$ at 11 points. 
But one of them belongs to $l$. We get that for $20 \cdot (11-1)=200$ 
intersection points of $R - l$ and $\phi(F_0 - l)$ there exist corresponding 
lines on $X$.

Let us find now $\# (R - l) \cap (\phi(F_0 - l))$. Let 
$r: \widetilde{P(V^*)} \to P(V^*)$ be the blowing up along $l$, \ \ $\tilde R$ 
and $\tilde F_0$ the proper inverse images of $R$ and $\phi (F_0)$ 
respectively. Then $\# (R - l) \cap (\phi(F_0 - l)) = 
<\tilde R, \tilde F_0> _{\widetilde{P(V^*)}}$. 
\medskip
{\bf Lemma 7.9.} $<\tilde R, \tilde F_0> _{\widetilde{P(V^*)}} = 36d +8$ 
where $d = \deg R$. 
\medskip
{\bf Proof.} $\widetilde{P(V^*)}$ is a variety of pairs $(V_2, V_4)$ where 
$V_2 \subset V_4 \subset V$ and $V_2 \subset U_3$; $r(V_2, V_4) = V_4$. There 
exists an isomorphism
$$\widetilde{P(V^*)} = P_{U^*} (O \oplus O \oplus O(1)) \overset{\mu}\to{\to} 
U^*$$ and $\mu(V_2, V_4) = V_2 \in U^*$. The ring $A(\widetilde{P(V^*)})$ has 
generators $D=\mu^*(c_1(O_{U^*}(1)))$, $\tilde H= c_1(O_{\mu}(1))$ and 
relations $D^3=0$, $\tilde H^3=D\tilde H^2$. Further, we have 
$$r^{-1}(l)=\tilde H - D, \ \ \ r^{-1}(H) = \tilde H,\ \ \ r_*(\tilde H^2) = 
r_*(\tilde HD) = r_*(D^2) = H^2$$ 
where $H = c_1(O(1))$. Let us find the class of $R$ in 
$A(\widetilde{P(V^*)})$. By definition, 
$$\{(V_2, V_4) \in \tilde R \iff P(V_2)=t_V \text{ for some } t \in 
c_{\Omega} \text{ and } P(V_4) \supset \pi(j(t,c))\}$$
where $c \in t_V \cap c_u$. So, $\mu(\tilde R) = c_{\Omega}$, and for $t 
\in c_{\Omega}$ \ \ \ $\mu^{-1}(t) \cap \tilde R$ is a pair of straight lines 
$\widetilde{\pi(j(t,c_i))}$ ($i=1,2$), where $\{c_1, c_2\}= t_V \cap c_u$. 
This means that $\tilde R$ is a divisor in $
\mu^{-1}( c_{\Omega}) $ and $O_{\mu^{-1}( c_{\Omega}) }(\tilde R) = L \otimes 
O_{\mu}(2) $ where $L$ is an invertible sheaf on $c_{\Omega}$. Since 
$\mu^{-1}( c_{\Omega})$ is a divisor on $\widetilde{P(V^*)} $ of degree $2D$, 
we have: $\cl (\tilde R) = 4 D
 \tilde H + nD^2$ where $n$ is defined by the condition $r_*(\tilde R)= dH^2 
= (4 + n) H^2$, i.e. $n=d-4$. 

Let us find now the class of $\tilde F_0$ in $A(\widetilde{P(V^*)})$. Let us 
assume that there exists a locally free sheaf $E$ of rank 2 on $P(V^*)$ such
that $\phi(F_0)$ is the set of zeros of a section of $E$. Then $\tilde F_0$ 
is the set of zeros of a section of $r^*(E) \otimes O_r(1)$ on 
$\widetilde{P(V^*)}$ and 
$$\cl (\tilde F_0) = c_2(r^*(E) \otimes O_r(1))= r^* c_2(E) + r^* c_1(E) 
\cdot c_1(O_r(1)) + c_1(O_r(1))^2.$$ 
We have $c_2(E) = \deg \phi(F_0) = 39H^2$, $c_1(O_r(1)) = D - \tilde H$, and 
$c_1(E)$ can be easily found by consideration of the exact sequence of normal 
sheaves for inclusions $l \hookrightarrow \phi(F_0) \hookrightarrow P(V^*)$: 
$$ 0 \to O_l(-1) \to 3O_l(1) \to E \vert _l \to 0$$
Namely, $\det E \vert _l = O_l(4)$, i.e. $c_1(E) = 4H$ and 
$$\cl \tilde F_0 = 
39 \tilde H^2 + 4\tilde H (D - \tilde H) + (D - \tilde H)^2 = 36 \tilde H^2 + 
2 \tilde H D + D^2$$

Really, $E$ does not exist, but it is clear that the above formalism is 
correct. So, $<\tilde R, \tilde F_0> _{\widetilde{P(V^*)}} = (4 \tilde H D+ 
(d-4)D^2)( 36 \tilde H^2 + 2 \tilde H D + D^2) = 36d +8.$  $\square$
\medskip
This means that there are $36 d + 8 - 200 $ straight lines of type $l_8(f_0)$ 
which 
meet $c_{\Omega}$. Since this number is $\ge 0$, we get that $d \ge 6$ and 
this number is $\ge 24$. So, we can recover $c_{\Omega}$ uniquely by $F_c$. 
$\square$
\medskip
Let us consider now the birational isomorphism $b: G_4 \to P^4$ (see (3.4)) 
which is the restriction to $G_4$ of the projection $p_{85}: P(V_8) \to P^4= 
P(V_8/U_3^*)$. 
\medskip
{\bf Proposition 7.10.} It is possible to recover the normcubic $c_3 \subset 
P^3 \subset P^4$ and the isomorphism $\psi_2: c_u \to c_3$ by $F_c$. 
\medskip
{\bf Proof.} For $f_0 \in F_0$ we let $l_5(f_0)= G(2, p_{85})(l_8(f_0))$, 
i.e. $l_5$ is a map from $F_0$ to $G(2, V_8/U_3^*)$. Since $l$ is a divisor
of $F_0$, then $l_5(f_0)$ is defined also for a generic point $f_0 \in l$. We 
shall consider $l_5(f_0)$ also as a straight line in $P^4$. A natural map of 
straight lines 
$p_{85} \vert _{l_8(f_0)}: l_8(f_0) \to l_5(f_0)$ is defined also for a 
generic point $f_0 \in l$. Analogously, the map $\xi _5 = p_{85} \circ \xi: 
P_{F_0}(\tau_{2,M}) \to P^4$ is defined at a generic point of 
$P_l(\tau_{2,M})$. 

It is clear that in $P_{F_0}(\tau_{2,M})$ we have: $S = \bar W \cap 
P_l(\tau_{2,M})$ (here we identify $c_u$ and $l$ via $\psi$). Let us consider 
a commutative diagram 
$$ \matrix S & \hookrightarrow & \bar W & \overset{\xi \vert _W}\to{\to} & 
W & \hookrightarrow & X & \hookrightarrow & G_4 \\
\downarrow & & \downarrow & & & & & & \downarrow \\
P_l(\tau_{2,M}) & \hookrightarrow & P_{F_0}(\tau_{2,M}) & & & 
\overset{\xi_5}\to{\to} & & & P^4 \endmatrix $$
where the left and middle vertical maps are inclusions, the right vertical map 
is $b$. Let $t_1 \in c_u$, then $ l_8(\psi(t_1)) \subset U^*$, and let 
$\{t', t''\}= l_8(\psi(t_1)) \cap c_{\Omega}$. We denote $\{t_1, t_2\} = 
t'_V \cap c_u$ and $\{t_1, t_3\} = t''_V \cap c_u$. Then $(t', t_1) \in S$, 
$(t', t_2) \in S$. The image of these elements in $G_4$ is $t'$, hence their 
images in $P^4$ belong to a bisecant $<\psi_2(t_1), \psi_2(t_2)>$ of $c_3$. 
Since $ c_{\Omega}$ can be recovered by $F_c$, we have: for $t_1 \in c_u$ we
can find images of points $(t', t_1)$ and $(t', t_2)$ in $P^4$ going along 
the left-lower way of the above diagram, and hence we can find the bisecant 
$<\psi_2(t_1), \psi_2(t_2)>$. Analogously we can find the bisecant 
$<\psi_2(t_1), \psi_2(t_3)>$; their intersection point indicates us the point 
$\psi_2(t_1) \in P^4$. $\square$ 
\medskip
Let us consider an isomorphism $I: P(V_8) \to P(V_8')$ and a map $b^{-1}: 
P^4 \to I(G_4)$ such that $b^{-1} \circ b = I \vert _{G_4}$. Since we can 
recover uniquely $c_3 \subset P^4$, we can also recover $b^{-1}: P^4 \to 
I(G_4)$ and $P(V_8')$ as the linear envelope of $I(G_4)$. 
\medskip
{\bf Proposition 7.11.} Isomorphism $I: P(V_8) \to P(V_8')$ can be recovered
by $F_c$. 
\medskip
{\bf Proof.} Let $\phi(f_0)=V_4$ where $f_0 \in F_0 - l$, and let $\{t\} = 
Q_G(V_4) \cap U^*$. Then (see 3.4) $b(Q_G(V_4)) = P^2$ and since points 
$b(\gamma_1(f_0))$, $b(\gamma_2(f_0))$ belong to $P^2$, we have: $l_5(f_0) 
\subset P^2$. This implies that $b^{-1}(l_5(f_0))$ is a conic on $I(Q_G(V_4))$ 
passing through $I(t)$, $I(\gamma_1(f_0))$, $I(\gamma_2(f_0))$. Since $b$ is 
the  restriction of $p_{85}: P(V_8) \to P^4$ to $G_4$, we have that for $a \in 
l_8(f_0)$ points $I(t)$, $I(a)$ and $b^{-1}(p_{85}(a))$ are collinear.

We can recover an isomorphism $S^2(c_3) \to I(U^*)$ and hence $I 
\vert _{U^*} : U^* \to I(U^*)$ as well. For $f_0$ we can recover a point 
$l_5(f_0) \cap P^3$. There exists only one bisecant of $c_3$ passing through 
this point, and the linear envelope of this bisecant and $l_5(f_0)$ is equal 
clearly $P^2 = b(Q_G(V_4))$. So, we can recover this $P^2$. 

Let now $f_1, f_2 \in F_0 - l$ be such that $l_5(f_1) \cap l_5(f_2) = \{m\}$, 
where $m =m(f_1, f_2) \in P^4$, and $m_1 \in l_8(f_1)$, $m_2 \in l_8(f_2)$ 
are inverse images of $m$: $p_{85}(m_1) = p_{85}(m_2) = m$. Let, further, 
$t_1 = Q_G(\phi(f_1)) \cap U^*$, $t_2 = Q_G(\phi(f_2)) \cap U^*$. Then the 
point $b^{-1}(m) $ belongs to both lines $<I(t_i), I(m_i)>$ ($i=1,2$), hence 
$I^{-1}(b^{-1}(m)) = <t_1, m_1> \cap <t_2, m_2>$, because these lines are
different. If $f_1$, $f_2$ are given, then we can construct
 points $m$, $m_1$, $m_2$, $t_1$, $t_2$, $b^{-1}(m)$, and hence 
$I^{-1}(b^{-1}(m))$. So, there are many points in $P(V_8)$ whose $I$-image 
can be constructed. To prove that they permit to recover uniquely the linear 
map $I$ it is sufficient to prove that all points of type 
$b^{-1}(m(f_1, f_2))$ do not belong to a hyperplane containing $I(U^*)$, or 
--- the same --- to prove that all points of type $m(f_1, f_2)$ do not belong 
to a hyperplane in $P^4$. 
\medskip
{\bf Lemma 7.12.} Let $Y$ be a surface on $G(2,V)$, $d = c_1(\tau_{2,V} 
\vert _Y)^2 - c_2(\tau_{2,V} \vert _Y)$. Then for a generic point $ y \in Y$ 
we have: the straight line $y_V$ meets $d-3$ straight lines of type $y_V'$ 
where $y' \in Y - y$. 
\medskip
{\bf Proof.} Follows easily from a calculation of intersection index on 
$\tilde G_y$. $\square$
\medskip
Chern classes of $l_5^*(\tau_2)$ (this sheaf corresponds to the inclusion 
$ l_5 \times F_0 \to G(2, (V_8/V_3)^*)$ \ ) \ can be easily found by using 
exact sequences (4.51) and 
$$0 \to \tau_2^* \vert _{F_0} \to \tau_{2,M}^* \vert _{F_0} \to 
i_*(\tau_{2,M}^* \vert _l) \to 0$$ where $i: l \to F_0$. We get: $d=168$. 
This means that straight lines of type $l_5(f_0)$ for $f_0 \in F$ do not 
belong to one hypersurface and a generic straight line of this type meets 
165 others. This implies
that all points of type $m(f_1, f_2)$ do not belong to a hyperplane in $P^4$. 
$\square$
\medskip
Now we can recover $G_4 \subset P(V_8)$ by $F_c$, because $G_4 = 
I^{-1}(b^{-1}(P^4))$. It is clear that for $f_0 \in F_0 - l$ 
$$\{\gamma_1(f_0), \gamma_2(f_0)\} = G_4 \cap l_8(f_0)$$ hence we can recover 
$W$ by $F_c$. Since $O_X(W) = O_X(21)$ and X is an intersection of quadrics 
in $P(V_8)$, $X$ is an intersection of all quadrics in $P(V_8)$ which contain 
$W$. So, we can recover uniquely $X$ by $F_c$, Q.E.D. $\square$
\bigskip
{\bf REFERENCES}
\bigskip
1. Gushel N.P. Special Fano threefold of genus 6. In: Constructive algebraic 
geometry. Yaroslavl, 1982, 
N 200, p. 26 - 34

2. Iskovskih V.A. Anticanonical models of algebraic threefolds. In: 
Contemporary problems of 
mathematics. v. 12. (Itogi nauki i tehniki, VINITI). Moscow, 1979, p. 59 - 157

3. Iskovskih V.A. Birational automorphisms of threefolds. In: Contemporary 
problems of
mathematics. v.12. (Itogi nauki i tehniki, VINITI). Moscow, 1979, p. 159 - 236 

4. Iskovskih V.A. On birational automorphisms of threefolds. Dokl. AN SSSR, 
1977, v. 234, N 4, p. 743 - 745

5. Iskovskih V.A. Birational automorphisms of a Fano theefold $V^3_6$. Dokl. 
AN SSSR, 1977, v. 235, N 3, p. 509 - 511 

6. Iskovskih V.A., Manin Yu. I. Quartic threefolds and counterexamples to the
Luroth problem. Math. Sbornik, 1971, v. 86, N 1, p. 140 - 166

7. Jozefiak T., Lascoux A., Pragacz P. Classes of determinantal varieties 
associated to symmetric and skew-symmetric matrices. Izv. AN SSSR. Ser. mat., 
1981, v 45, p. 662 - 673

8. Kulikov Vik. S. Degenerations of K3 surfaces and Enriques surfaces. Izv. 
AN SSSR. Ser. mat., 1977, v. 41, N 5, p. 1008 - 1042

9. Logachev D. Yu. Abel - Jacobi map of a Fano threefold of genus 6 is an 
isogeny. In: Constructive 
algebraic geometry. Yaroslavl, 1982, N 200, p. 67 - 76

10. Logachev D. Yu. Isomorphism of middle Jacobians of some dual sections of 
determinantals. Izv. AN SSSR, Ser. Mat., 1982, v. 46, N 2, p. 269 - 275

11. Tikhomirov A. S. Geometry of Fano surface of a double space $P^3$ 
ramified at a quartic. Izv. AN SSSR, Ser. mat., 1980, v. 44, N 2, p. 415 - 422

12. Tikhomirov A. S. Middle Jacobian of a double space $P^3$ ramified at a 
quartic. Izv. AN SSSR, Ser. mat., 1980, v. 44, N 6, p. 1329 - 1377

13. Tikhomirov A. S. Fano surface of a double cone of Veronese. Izv. AN SSSR, 
Ser. mat., 1981, v. 45, N 5, p. 1121 - 1197

14. Tikhomirov A. S. Singularities of a theta divisor of a double space $P^3$ 
of index 2. Izv. AN SSSR, Ser. mat., 
1982, v. 46, N 5, p. 1062 - 1081

15. Tikhomirov A. S. Special double space $P^3$ of index 2. In: Constructive 
algebraic geometry. Yaroslavl, 1982, N 200, p. 121 - 129 

16. Tyurin A. N. Middle Jacobian of threefolds. In: Contemporary problems of 
mathematics. v. 12. (Itogi nauki i tehniki, VINITI). Moscow, 1979, p. 5 - 57

17. Tyurin A. N. On Fano surface of a non-singular cubic in $P^4$. Izv. AN 
SSSR, Ser. mat., 1970, v. 34, N 6, p. 1200 - 1208

18. Tyurin A. N. Geometry of a Fano surface of a non-singular cubic $F 
\subset P^4$ and Torelli theorems. Izv. AN SSSR, Ser. mat., 1971, v. 35, p. 
498 - 529

19. Tyurin A. N. Geometry of a Poincar\'e divisor of a Prym variety. Izv. AN 
SSSR, Ser. mat., 1975, v. 39, N 5, p. 1003 - 1043

20. Tyurin A. N. Five lectures on threefolds. Uspehi Mat. Nauk, 1972, v. 27, 
N 5, p. 3 - 50

21. Tyurin A. N. On intersection of quadrics. Uspehi Mat. Nauk, 1975, v. 30, 
N 6, p. 51 -99

22. Hashin S. N. Birational automorphisms of a double cone of Veronese of 
dimension 3. Vestnik MGU. Ser. mat., mech. 1983, N 1. 

23. Altman A., Kleiman S. Foundations of theory of Fano schemas. Preprint. 

24. Andreotti A., Mayer A. On period relations for abelian integrals. Ann.
Scuola Norm. Sup. Pisa, 1967, v. 21, N 2, p. 189 - 238

25. Barton C., Clemens C.H. A result on the integral Chow ring of a generic 
principally polarized complex abelian variety of dimension four. Comp. Math., 
1977, v. 34, fasc. 1, p. 49 - 67

26. Beauville A. Vari\'eti\'es de Prym et Jacobiennes interm\'ediaires. Ann. 
Scient. Ecole Norm. Sup., 1977, v. 10, p. 304 - 392

27. Beauville A. Prym varieties and Schottky problem. Invent. Math., 1977, v. 
41, N 2, p. 149 - 196

28. Beauville A. Sous - vari\'eti\'es speciales des vari\'eti\'es de Prym. 
Preprint. 

29. Bombieri E., Swinnerton-Dyer H. On the local zeta function of a cubic 
threefold. Ann. Scuola Norm. Sup. Pisa, 1967, v. 21, p. 1 - 29

30. Clemens C. H. Degeneration of K\"ahler manifolds. Duke Maht. J., 1977, 
v. 44, N 2, p. 215 - 290

31. Clemens C. H. Double solids. Preprint.

32. Clemens C. H., Griffiths Ph. The intermediate Jacobian of a cubic 
threefold. Ann. of Math., 1972, v. 95, N 2, p. 281 - 356

33. Collino A., Murre J. P., Welters G. E. On the surface of conics on a 
quartic threefold. Rend. Sem. Math. Univ. Politecn. Torino. 1980, v. 38, N 1. 

34. Fano G. Sul sistema $\infty ^2$ di rette contenuto in una variet\`a 
cubica generale dello spazio a quattro dimensioni. Atti R. Acc. Sc. Torino,
1904, v. 39, p. 778 - 792

35. Fano G. Osservazioni sopra alcune variet\`a non razionali aventi tutti i 
generi nulli.. Atti R. Acc. Sc. Torino, 1915, v. 50, p. 1067 - 1072

36. Fano G. Sulle seczioni spaziali della variet\`a Grassmanniana delle rette 
spazio a cinque dimensioni. Rend. R. Accad. Lincei. 1930, v. 11, N 6, p. 
329 - 356

37. Griffiths Ph. On periods of integrals on algebraic manifolds I. Amer. J. 
Math. 1968, v. 90, p. 568 - 626

38. Griffiths Ph. On periods of integrals on algebraic manifolds II. Amer. J.
Math. 1968, v. 90, p. 805 - 865

39. Griffiths Ph., Harris J. Principles of algebraic geometry. New York, 1977

40. Griffiths Ph., Schmid W. Recent developments in Hodge theory: a 
discussion of techniques and results. In: Proc. Intern. Colloq. on discrete 
subgroups of Lie groups. Bombay, Oxford Univ. Press, 1973

41. Hironaka H. Resolution of singularities of an algebraic variety over a
field of characteristic zero. I, II. Ann. of Math., 1964, v. 79, p. 109 - 326

42. Lieberman D. Higher Picard varieties. Amer. J. Math., 1968, v. 90, N 4, 
p. 1165 - 1199

43. Mumford D. Prym varieties I. Contrib. Anal. New York - London, 1974, p. 
325 - 350

44. Murre J. P. Algebraic equivalence modulo rational equivalence on a cubic 
threfold. Comp. Math., 1972, v. 25, p. 161 - 206

45. Murre J. P. Reduction of a proof of the non-rationality of a non-singular 
cubic threefold to a result of Mumford. Comp. Math., 1973, v. 27, p. 63 - 82

46. Murre J. P. Some results on a cubic threefold. In: Lecture Notes in 
Math., 1974, v. 412.

47. Reid M. The complete intersection of two or more quadrics. Ph.D. thesis. 
Cambridge Univ., 1972

48. Todd J. A. The locus representing the lines of four-dimensional space. 
Proc. London Math. Soc., 1930, v. 30, p. 513 - 550

49. Welters G. E. The Fano surface of lines on a double $P^3$ with 4th order 
discriminant locus I. Univ. Utrecht. Preprint N 123, 60 p. 

50. Welters G. E. Divisor varieties, Prym varieties and a conjecture of 
Tyurin. Univ. Utrecht. Preprint N 139, 57 p. 

51. Welters G. E. Abel - Jacobi isogenies for certain types of Fano 
threefolds. Ph.D. thesis. Amsterdam, 1981

52. Zucker S. Generalized intermediate Jacobians and the theorem on normal 
functions. Invent. Math., 1976, v. 33, N 3, p. 185 - 222

\enddocument